\documentclass[a4paper,12pt]{article}
\usepackage[utf8]{inputenc}

\usepackage[spanish]{babel}
\spanishdecimal{.}

\usepackage[inline]{enumitem}   

\usepackage{mathrsfs}
\usepackage{graphicx} 
\usepackage{graphics} 
\usepackage{caption}
\usepackage[x11names,table]{xcolor}
\usepackage{amsmath,amssymb,amsfonts,latexsym,cancel,amsthm}

\usepackage{geometry}
\usepackage{multirow}
\usepackage{multicol,caption}
\usepackage{float}

\usepackage{subfigure} 

\usepackage{rotating}
\usepackage{adjustbox}
\usepackage[graphicx]{realboxes}
\usepackage[justification=centering]{caption}
\newtheorem{ejemplo}{Ejemplo}

\usepackage{amsthm}

\title{Aplicando diferencias finitas para resolver ecuaciones y sistemas de 
       ecuaciones diferenciales parciales sobre dominios planos irregulares 
       simplemente conexos y no conexos}
\author{Miriam Sosa-Díaz\\
  \and
   Faustino Sanchez Gardu\~{n}o\\
   \texttt{faustinos403@gmail}}

\begin{document}




\maketitle

\begin{abstract} 
Basandonos en los fundamentos del
{\it m\'{e}todo de exhausti\'{o}n, las f\'{o}rmulas 
de diferencias finitas} y {\it los m\'{e}todo de 
diferencias finitas} deducimos un nuevo m\'{e}todo
para obtener soluciones num\'{e}ricas de ecuaciones
y sistemas de ecuaciones diferenciales parciales 
lineales y no lineales sobre regiones irregulares
simplemente conexas y no conexas. Despu\'{e}s 
aplicamos este m\'{e}todo para obtener soluciones
num\'{e}ricas de algunos problemas lineales y no 
lineales sobre regiones irregulares con condiciones
de Dirichlet en la frontera.

\end{abstract}

\section{Introducci\'{o}n}

                En la ciencia e ingenier\'{i}a la resoluci\'{o}n de problemas
                que involucran ecuaciones diferenciales parciales lineales 
                y no lineales es determinante en la toma de decisiones.                
                
                El an\'{a}lisis de estos problemas desde un punto de vista 
                te\'{o}rico muchas veces, no es suficiente para conocer por 
                completo la soluci\'{o}n. Por lo que tenemos que recurrir 
                a los m\'{e}todos num\'{e}ricos, que bajo ciertas condiciones
                nos pueden proveer de buenas soluciones num\'{e}ricas. 
                
                M\'{a}s all\'{a} de s\'{o}lo aplicar un m\'{e}todo para encontrar
                soluciones num\'{e}ricas, lo ideal es realizar el an\'{a}lisis
                num\'{e}rico que corresponde al problema completo. De esta 
                manera los resultados que obtengamos ser\'{a}n m\'{a}s confiables.
               
                Usualmente la elecci\'{o}n de un m\'{e}todo para resolver un
                problema de este tipo depende de las caracter\'{i}sticas de 
                la ecuaci\'{o}n o del sistema de ecuaciones diferenciales 
                parciales. Un factor determinante y clave en la elecci\'{o}n
                del m\'{e}todo es el tipo de regi\'{o}n $\Omega$ donde se desea 
                resolver el problema. Por ejemplo, 
                \begin{enumerate}
                       \item si $\Omega$ es un cuadrado o un rect\'{a}ngulo se suele 
                             usar diferencias finitas,
                       \item si $\Omega$ es una regi\'{o}n irregular se suele usar
                             elemento finito u otro m\'{e}todo. 
                \end{enumerate}                                                 
                
                Para aplicar cualquiera de los m\'{e}todos antes mencionados, debemos 
                discretizar la regi\'{o}n $\Omega$, es decir, debemos asociar una 
                malla a la regi\'{o}n $\Omega$.
                
                En el caso de diferencias finitas, a la regi\'{o}n $\Omega$ se le asocia
                una malla estructurada, mientras que en elemento finito la malla es no 
                estructurada. \'{E}sta es una de las caracter\'{i}stica que determinar\'{a} 
                la estructura de la discretizaci\'{o}n de la ecuaci\'{o}n o del sistema 
                de ecuaciones diferenciales parciales que se quiere resolver. 
                
                Con la finalidad de establecer un nuevo m\'{e}todo para resolver
                ecuaciones diferenciales parciales sobre regiones irreguales
                simplemente conexas y no conexas veamos, de manera muy breve,
                c\'{o}mo se han abordado algunos problemas matem\'{a}ticos desde
                la antig\"{u}edad.
                
                No vamos a encontara la solucion del problema, en vez de eso 
                lo que haremos es construir un conjunto finito de puntos
                que son aproximaciones a la solucion
                \textquestiondown como podemos construir un conjuntos finito de 
                puntos  ?
                que verifiquen aproximadamente una ecuacion diferencial
                
                Los objetivos del analisis numerico
                Encontrar una aproximacion de la solucion para un problema dado
                Determinar una cota para la exactitud de la aproximacion
            
\section{Fundamentos}
                
                Uno de los m\'{e}todos m\'{a}s antiguos que se ha aplicado para la 
                resoluci\'{o}n de problemas, tales como el c\'{a}lculo de \'{a}reas,
                per\'{i}metros, volumenes, entre otros, ha sido el m\'{e}todo de
                exhausti\'{o}n.
                
                En \cite{Ribnikov}, el autor nos explica que en el m\'{e}todo de 
                exhausti\'{o}n debemos de construir una sucesi\'{o}n de soluciones
                aproximadas tales que \'{e}stas convergan a la soluci\'{o}n del
                problema en cuesti\'{o}n. El autor expone el m\'{e}todo de la 
                siguiente manera.

                \begin{center}
                \begin{minipage}{0.9\linewidth}
                
                {\it  \textquotedblleft                   
                   Uno de los m\'{a}s antiguos m\'{e}todos de este genero es el 
                   m\'{e}todos de exhausti\'{o}n. Su creaci\'{o}n se le atribuye 
                   a Eudoxo. Ejemplos de su utilizaci\'{o}n est\'{a}n expuestos 
                   en el libro duod\'{e}cimo de los \textquotedblleft 
                   Elementos\textquotedblright de Euclides y en una serie de obras  
                   de Arqu\'{i}mides. El m\'{e}todo de exhausti\'{o}n se aplicaba 
                   al c\'{a}lculo de las \'{a}reas de figuras, vol\'{u}menes de 
                   cuerpos, longitud de cuerdas, b\'{u}squeda de la subtangentes
                   de las curvas, etc. La esencia matem\'{a}tica del m\'{e}todo 
                   (se entiende que en forma algo diferente a la forma de 
                   exposici\'{o}n de los griegos antiguos) consiste en la sucesi\'{o}n
                   de las siguientes operaciones: }
                   
                \end{minipage}
                \end{center}       
                
                \begin{figure}[h]
                \centering
                \begin{tabular}{||c||}
                \hline                
                \includegraphics[width=6cm,height=5cm]{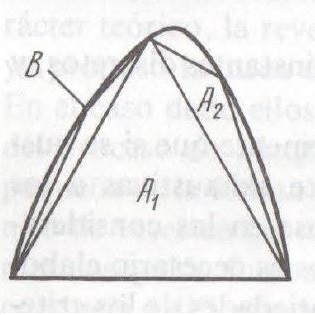} \\             
                \hline
                \end{tabular}
                \caption{Imagen tomada de \cite{Ribnikov}.}
                \label{fig:tri}
                \end{figure}                                                
                   
                \begin{center}
                \begin{minipage}{0.9\linewidth}
                {\it                    
                                      
                   \begin{enumerate}
                   
                          \item si es necesario, por ejemplo cuadrar la figura B 
                                (figura \ref{fig:tri}), entonces como primer paso 
                                 en esta figura se inscribe una sucesi\'{o}n de 
                                 otras figuras $A_1$, $A_2$, $\ldots$, $A_n$, $\ldots$
                                 cuyas \'{a}reas crecen mon\'{o}tonamente y para cada 
                                 figura de esta sucesi\'{o}n pueden ser determinadas.
                                 
                          \item las figuras $A_k$ $(k=1,2, \ldots)$ se eligen de tal 
                                forma que la diferencia positiva $B-A_k$ pueda ser hecha
                                tan peque\~{n}a como se quiera;
                                
                          \item del hecho de la existencia y construccion de las figuras
                                descritas se hace la deducci\'{o}n sobre la acotaci\'{o}n 
                                superior de la sucesi\'{o}n de las figuras inscritas
                                \textquotedblleft agotadoras\textquotedblright ;
                                
                          \item impl\'{i}citamente, por lo general con ayuda de otras 
                                consideraciones te\'{o}ricas y pr\'{a}cticas, se busca $A$, 
                                es decir, el l\'{i}mite de la sucesi\'{o}n de las figuras 
                                inscritas;
                                
                          \item se demuestra, para cada problema por separado, que $A=B$,
                                esto es, que el l\'{i}mite de la sucesi\'{o}n de las figuras 
                                inscritas es igual al \'{a}rea $B$. La demostraci\'{o}n, 
                                como regla, se realiza por reducci\'{o}n al absurdo. Sea
                                $A \neq B$. Entonces $B>A$ \'{o} $B<A$. Supongamos que $B>A$,
                                elijamos un elemento $A_n$ de la sucesi\'{o}n tal que
                                $B-A_n < B-A$. Esto es posible para cada cualquier diferencia
                                fijada $B-A$. Pero entonces debe ser $A_n > A$ y esto es 
                                imposible ya que en realidad $A > A_n$ para todo $n$ finito. 
                                La suposici\'{o}n contraria $(B<A)$ tambi\'{e}n conduce 
                                a una contradicci\'{o}n, puesto que puede elegirse $A_n$ tal 
                                que $A-A_n < A-B$. Pero entonces debe obtenerse que $A_n>B$ 
                                y esto es imposible.           
 
                   \end{enumerate}
                 }
                   
                \end{minipage}
                \end{center}       
                   
                \begin{center}
                \begin{minipage}{0.9\linewidth}
                {\it                                       
 
                   Con el m\'{e}todo de exhausti\'{o}n se demuestra, de esta forma 
                   la unicidad del l\'{i}mite. En combinaci\'{o}n con otros m\'{e}todos
                   es \'{u}til para la b\'{u}squeda de l\'{i}mites. Sin embargo,
                   este m\'{e}todo no puede dar la soluci\'{o}n de la cuesti\'{o}n
                   sobre la existencia del l\'{i}mite.                
                   \textquotedblright (pp. $75-76$).
                }

                \end{minipage}
                \end{center}       
                
                Por otro lado, al inicio del cap\'{i}tulo de integraci\'{o}n 
                en \cite{Spivak}, el autor realiza la siguiente reflexi\'{o}n con 
                respecto al problema del c\'{a}lculo de \'{a}reas.

                \begin{center}
                \begin{minipage}{0.9\linewidth}
                \vspace{5pt}                        
                
                {\it  \textquotedblleft                                      
                   En geometr\'{i}a elemental se deducen f\'{o}rmulas para las \'{a}reas 
                   de muchas figuras planas, pero un poco de reflexi\'{o}n hace ver que 
                   raramente se da una definici\'{o}n aceptable de \'{a}rea. El \'{a}rea 
                   de una regi\'{o}n se define a veces como el n\'{u}mero de cuadrados 
                   de lado unidad que caben en la regi\'{o}n. Pero esta definici\'{o}n es
                   totalmente inadecuada para todas las regiones con excepci\'{o}n de las
                   m\'{a}s simples. Por ejemplo, el c\'{i}rculo de radio $1$ tiene por
                   \'{a}rea el n\'{u}mero irracional $\pi$, pero no est\'{a} claro en
                   absoluto cu\'{a}l es el significado de $\ll$ $ \pi$ cuadrados $\gg$. 
                   Incluso si consideramos un c\'{i}rculo de radio $\frac{1}{\sqrt{\pi}}$ 
                   cuya \'{a}rea es $1$, resulta dif\'{i}cil explicar de qu\'{e} manera 
                   un cuadrado unidad puede llenar este c\'{i}rculo, ya que no parece 
                   posible dividir el cuadrado unidad en pedazos que puedan ser 
                   yuxtapuestos de manera que former un c\'{i}rculo.\textquotedblright (p. $317$).
                }  
                                
                \vspace{5pt}
                \end{minipage}
                \end{center}                    
                
                Finalmente en \cite{Apostol} se le\'{e} lo siguiente
                
                \begin{center}
                \begin{minipage}{0.9\linewidth}
                \vspace{5pt}                        
                
                {\it  \textquotedblleft 
                   Definici\'{o}n. Sea $S$ un subconjunto de un intervalo $[a,b]$ en $E_n$. 
                   Para toda partici\'{o}n $P$ de $[a,b]$ definimos $ \underline {J}(P,S)$ como 
                   la suma de las medidas de aquellos subintervalos de $P$ que \'{u}nicamente 
                   contienen puntos interiores de $S$ y $\overline{J}(P,S)$ como la suma de las
                   medidas de los subintervalos de $P$ que contienen puntos de $S \cup \partial S$.
                   (Recordar las Definiciones $3-24$ y $8-38$.) Los n\'{u}meros
                   
                   \begin{eqnarray}
                          \underline{c}(S) = sup \{ \underline {J}(P,S) | P \in \mathcal{P} [a,b]  \}, \\
                          \overline{c}(S) = inf   \{   \overline{J}(P,S)  |  P \in \mathcal{P} [a,b] \}
                   \end{eqnarray}

                   se llaman, respectivamente, el contenido de Jordan interior y exterior
                   ($n-$dimensional) de S. Se dice que el conjunto $S$ es medible-Jordan 
                   si $\underline{c}(S) =  \overline{c}(S) $, en cuyo caso este valor 
                   com\'{u}n se llama el contenido de Jordan de $S$, representado por $c(S)$. 
                }   
                   
                \end{minipage}
                \end{center}      
                
                \begin{figure}[t]
                \centering
                \begin{tabular}{||c||}
                \hline                
                \includegraphics[width=6cm,height=5cm]{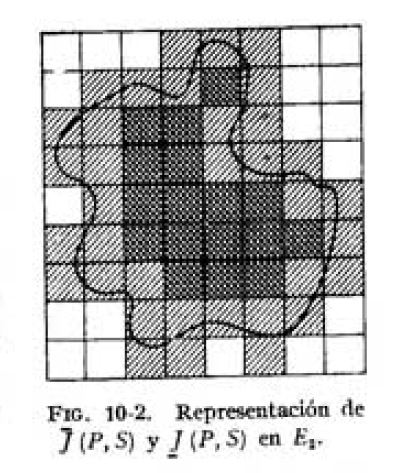} \\             
                \hline
                \end{tabular}
                \caption{Representacion de $\overline{J}(P,S)$ y $\underline {J}(P,S)$ en $E_2$.}
                \label{fig:exa}
                \end{figure}     
                   
                \begin{center}
                \begin{minipage}{0.9\linewidth}
                {\it                    
                   
                   Nota. Es f\'{a}cil verificar que $\underline{c}(S)$ y $\overline{c}(S)$ 
                   dependen tan s\'{o}lo de $S$ y no del intervalo $[a,b]$ que contiene $S$. 
                   Tambi\'{e}n, $0 \leq \underline{c}(S) \leq \overline{c}(S)$. Para 
                   conjuntos finitos tenemos que $ \underline{c}(S) = \overline{c}(S)  =0$. 
                   El contenido exterior de Jordan $n$-dimensional de un conjunto acotado
                   $r$-dimensional es cero si $r<n$. Los conjuntos elementales $S$ son
                   medibles-Jordan y $c(S) = \mu(S)$. (Hablando con precisi\'{o}n, 
                   deber\'{i}amos poner $c_n$ mejor que $c$ para poner de manifiesto que el
                   contenido es $n$-dimensional, pero no se originan confusiones al omitir 
                   esta alusi\'{o}n a $n$.) Se dice que los conjuntos medibles-Jordan $S$ en
                   $E_2$ tienen \'{a}rea $c(S)$. En este caso, las sumas $\underline{c}(S)$
                   y $\overline{c}(S)$ representan aproximaciones del \'{a}rea desde el 
                   $\ll$ interior $\gg$ y desde el $\ll$ exterior $\gg$, respectivamente. 
                   Esto est\'{a} representado en la figura \ref{fig:exa}, en la cual los 
                   rect\'{a}ngulos sombreados ligeramente est\'{a}n considerados en 
                   $\overline{c}(S) $, y los rect\'{a}ngulos de sombreado m\'{a}s intenso 
                   en $\underline{c}(S)$. Para los conjuntos en $E_2$, $c(S)$ se llama el 
                   volumen de $S$. \textquotedblright (p. $246$).
                }  
                                
                \vspace{5pt}
                \end{minipage}
                \end{center}

\section{Discretizaci\'{o}n de un dominio irregular} \label{Metodo}

                Con el objetivo de resolver problemas que involucren la 
                soluci\'{o}n de una ecuaci\'{o}n \'{o} de un sistema de 
                ecuaciones diferenciales parciales lineales y no lineales
                sobre un dominio irregular vamos a construir una sucesi\'{o}n 
                de mallas uniformes y estructuradas tal que converga al dominio
                irregular. \'{E}sto nos permitir\'{a} aplicar el m\'{e}todo  
                de diferencias finitas sobre cada una de las mallas para
                construir una sucesi\'{o}n de soluciones.
                
                Para \'{e}sto consideremos cualquier dominio irregular $\Omega$
                y hagamos lo siguiente.                                
                                                                
                \begin{enumerate}
                        \item Construir un rect\'{a}ngulo $R$ en el plano tal que $\Omega \subset R$.                         
                        \item Discretizar el rect\'{a}ngulo $R$ por una malla uniforme $G_R$ 
                              de $mx \times my$, eligiendo $\Delta x $ y $\Delta y $ \'{o}  
                              $mx$ y $my$. Nos referiremos a los puntos de $G_R$ como
                              $(j \Delta x, k \Delta y)$ para $j=0,1,\ldots, mx$ 
                              y $k=0,1,\ldots, my$.
                        \item Intersectar el dominio $\Omega$ con la malla $G_R$. A partir de 
                              \'{e}sto determinamos los siguientes conjuntos.
                              \begin{enumerate}
                                     \item El conjunto de puntos que est\'{a}n en el interior
                                           de $\Omega$.
                                     \item El conjunto de puntos que est\'{a}n en la frontera 
                                           de $\Omega$.
                                     \item El conjunto de puntos que est\'{a}n en el complemento 
                                           de $\Omega$
                              \end{enumerate}
                                
                              Debe de quedar claro que para definir los conjuntos anteriores
                              estamos considerando \'{u}nicamente la posici\'{o}n de los 
                              puntos, pero con respecto al dominio, es decir, si \'{e}stos 
                              est\'{a}n en el interior, en la frontera \'{o} en el complemento
                              del dominio $\Omega$. 
                              
                              As\'{i} mismo, tambi\'{e}n debe de quedar claro que la manera 
                              en que definimos los conjuntos anteriores difiere de la que
                              se aborda en la definici\'{o}n del contenido de Jordan, ya que
                              en \'{e}sta, el conjunto de puntos se construye en base a la
                              posici\'{o}n de los rect\'{a}ngulos.                                                            
                              
                              Hay que tener presente que estamos interesados en determinar
                              una soluci\'{o}n aproximada de una ecuaci\'{o}n o de un 
                              sistema de ecuaciones en cada uno de los puntos interiores de 
                              $\Omega$ y en el caso de que consideraramos la definici\'{o}n del 
                              contenido de Jordan habr\'{i}a puntos interiores que 
                              quedar\'{i}an excluidos y creemos que \'{e}sto afectar\'{i}a 
                              la exactitud de la aproximaci\'{o}n por la forma en la que se 
                              har\'{a}.
                              
                              Ahora pasemos a definir como vamos a considerar a 
                              la frontera, para esto pasemos al siguiente punto.
                                
                        \item Definamos los siguientes conjuntos                                 
                              \begin{enumerate}                                        
                                     \item El conjunto de puntos $\overline {\Omega}$ 
                                           determinado por todos los puntos que est\'{a}n 
                                           en el interior de $\Omega$ y los puntos que
                                           determinan su frontera, {\it pero con respecto a
                                           la malla}.   
                                              
                                           Estos \'{u}ltimos puntos pueden \'{o} no estar  
                                           en la frontera del dominio $\Omega$.
                                              
                                     \item El conjunto de puntos $\underline{\Omega}$ 
                                           que est\'{a}n \'{u}nicamente en el interior 
                                           de $\Omega$. En este caso hay puntos interiores
                                           que pasan a ser puntos frontera. 
                              \end{enumerate} 
                              
                              Una vez que ya hemos construido $\overline {\Omega}$
                              y $\underline{\Omega}$ es importante verificar que los puntos
                              interiores de ambas discretizaciones, $\overline {\Omega}$ 
                              y $\underline{\Omega}$, est\'{e}n conectados. Ya que en caso 
                              contrario estar\'{i}amos asociando una discretizaci\'{o}n 
                              discontinua que corresponder\'{i}a a un dominio no continuo, 
                              a un dominio continuo, lo cual no es correcto. De hecho 
                              sino verificamos est\'{a} propiedad impl\'{i}citamente
                              estar\'{i}amos cambiando el problema a resolver.
                                                            
                              Ahora bien, por la forma en que se construyeron $\overline {\Omega}$ 
                              y $\underline{\Omega}$, resulta claro que la colecci\'{o}n
                              de todos los puntos de cada uno de estos dos conjuntos forman
                              una malla uniforme y estructurada. De hecho, 
                              {\it diremos que $\overline {\Omega}$ y $\underline{\Omega}$
                              aproximan al dominio $\Omega$ mediante una extrapolaci\'{o}n  
                              y una interpolaci\'{o}n de una malla uniforme y estructurada,
                              respectivamente}.       
                              
                              N\'{o}tese que por para cada par de valores $\Delta x $ 
                              y $\Delta y $ \'{o} $mx$ y $my$ estamos generando una
                              malla por extrapolaci\'{o}n $\overline {\Omega}$ y otra 
                              por interpolaci\'{o}n $\underline{\Omega}$ que aproximan 
                              al dominio $\Omega$ y que denotaremos como
                              $\overline {\Omega}_{m,n}$ y $\underline{\Omega}_{m,n}$
                              respectivamente, pero si no se origina confusi\'{o}n 
                              usualmente omitiremos los sub\'{i}ndices $m,n$.                                                             
                                                                                                                                                      
                              De esto se sigue que, estamos generando dos 
                              sucesiones de mallas uniformes y estructuradas. Este 
                              hecho da lugar a que tendremos dos sucesiones que 
                              aproximan a la soluci\'{o}n exacta del problema, como
                              veremos a continuaci\'{o}n. 
                              
                        \item Debido a la estructura de $\overline {\Omega}$ y de $\underline{\Omega}$,
                              resulta claro que podemos aplicar el m\'{e}todo de diferencias
                              finitas, ver por ejemplo \cite{Jwth} y \cite{Jwt2}, o alg\'{u}n otro 
                              m\'{e}todo adecuado que nos permita encontar una soluci\'{o}n
                              aproximada sobre cada punto de $\overline {\Omega}$ \'{o} de
                              $\underline{\Omega}$. Para \'{e}sto simplemente aplicamos el
                              m\'{e}todo elegido                                                                                                              
                              \begin{enumerate}
                                     \item al conjunto de puntos
                                           $(j \Delta x, k \Delta y) \in \overline{\Omega}$, 
                                             
                                     \item al conjunto de puntos 
                                           $(j \Delta x, k \Delta y) \in \underline {\Omega}$,                                               
                              \end{enumerate} 
                              y discretizamos la variable temporal $t$ en caso de que \'{e}sta 
                              est\'{e} presente.                                                                                                                                                     
                              
                              As\'{i} c\'{o}mo por cada par de valores de  $\Delta x $ 
                              y $\Delta y $ \'{o} $mx$ y $my$ estamos generando una malla
                              por extrapolaci\'{o}n y otra por interpolaci\'{o}n se sigue 
                              que por cada malla estamos generando una soluci\'{o}n: una                               
                              por {\it extrapolaci\'{o}n} y otra por {\it interpolaci\'{o}n}.  
                                                                                                                                                                                    
                \end{enumerate}

                Con el objetivo de ver que tan bien funciona el m\'{e}todo que
                acabamos de describir vamos a ver un ejemplo.

                      

                \begin{ejemplo} \label{ejemplo1} Consideremos el siguiente problema
                      con condiciones de Dirichlet en la frontera
                      \begin{eqnarray}
                       \label{Parabola:eq}
                            \nabla^2 a & = & 4,\hspace*{0.3cm}\text{para}\hspace*{0.3cm}(x,y)\in\Omega, \\
                                     a & = & f(x,y),\hspace*{0.3cm} \text{para} \hspace*{0.3cm} (x,y) \in   \partial  \Omega , \nonumber
                      \end{eqnarray}    
                      donde $a=a(x,y)$, $\Omega$ es el dominio determinado por 
                      la par\'{a}bola $y=x^2$, la recta $x=1$ y el eje horizontal,
                      $f(x,y) = (x+y)^2$ y $\partial \Omega$ es la frontera de $\Omega$.
                      
               
                                           
                      \begin{enumerate}
                            \item Vamos a resolver este problema usando el m\'{e}todo 
                                  de la secci\'{o}n \ref{Metodo}.
                            \item Luego vamos a comparar las soluciones obtenidas 
                                  con la soluci\'{o}n exacta
                                  \begin{equation*}
                                        a(x,y) = (x+y)^2,
                                  \end{equation*}
                                  sobre $\Omega$.
                            \item Finalmente vamos a comparar el m\'{i}nimo, el 
                                  m\'{a}ximo y sus posiciones con los valores
                                  exactos. En este caso el valor m\'{i}nimo y 
                                  m\'{a}ximo exacto es $0$ y $4$, los cuales se
                                  localizan en $(0,0)$ y en $(1,1)$, respectivamente.
                      \end{enumerate}                                             
                      
                      Primero que nada notemos que $\Omega$ es un dominio irregular,
                      as\'{i} que de acuerdo al m\'{e}todo \ref{Metodo}, para resolver 
                      este problema aplicando diferencias finitas debemos de considerar
                      un rect\'{a}ngulo $R$ tal que $\Omega \subset R$, en este caso es 
                      suficiente considerar
                      \begin{equation*}
                            R= [-0.2,1.2] \times [ -0.2, 1.2].
                      \end{equation*}

                      Ahora vamos a discretizar $R$, como usualmente se hace,
                      considerando $\Delta x = 0.1 $ y $\Delta y = 0.1$ y
                      despu\'{e}s intersectamos la malla obtenida con el dominio
                      $\Omega$. \'{E}sto nos permite determinar el conjunto de 
                      puntos interiores del dominio $\Omega$, los cuales se 
                      pueden ver en la figura \ref{fig:puntosInterioresA}.
                      
                      \begin{figure}[t]
                      \centering
                      \begin{tabular}{||c||}
                      \hline                
                      \includegraphics[width=6cm,height=5cm]{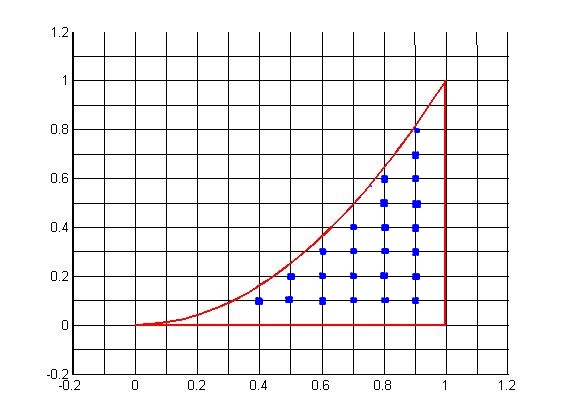} \\             
                      \hline
                      \end{tabular}
                      \caption{Puntos interiores de un dominio irregular $\Omega$.}
                      \label{fig:puntosInterioresA}
                      \end{figure}   
                      
                      A partir de este conjunto de puntos vamos a construir dos mallas, 
                      $\overline {\Omega}$ y $\underline{\Omega}$, que aproximan al 
                      dominio $\Omega$.
                      
                      Para construir la primera malla, $\overline {\Omega}$, primero 
                      consideramos todos los puntos interiores que ya hab\'{i}amos 
                      determinado y despu\'{e}s agregamos todos los puntos necesarios,
                      con respecto a la malla, para construir la frontera de \'{e}stos.
                      Estos nuevos puntos pueden estar en la frontera del dominio
                      original.
                      
                      Mientras que para construir la segunda malla, 
                      $\underline{\Omega}$, volvemos a considerar de nuevo el conjunto
                      de todos los puntos interiores que ya hab\'{i}amos determinado 
                      y en este caso la frontera son aquellos puntos que no est\'{a}n
                      en el interior de este conjunto (de puntos interiores), es decir 
                      algunos puntos interiores pasan a ser puntos frontera.
                      
                      La estructura final de las dos mallas uniformes y estructuradas, 
                      $\overline {\Omega}$ y $\underline{\Omega}$, se pueden visualizar
                      en la figura \ref{fig:pFronteraUtilizable}, respectivamente. En 
                      \'{e}sta los puntos de color azul representan los puntos interiores
                      de cada malla, los puntos de color rosa son los puntos que se 
                      agregagon para formar a la frontera de $\overline {\Omega}$ y los
                      puntos azules rodeados de color rosa son los puntos interiores que 
                      pasaron a formar la frontera de $\underline{\Omega}$. En cualquiera 
                      de los dos casos vamos a considerar que el valor de cada punto que
                      forma la frontera es conocido y est\'{a} dado por la segunda 
                      igualdad de (\ref{Parabola:eq}).

                      \begin{figure}[t] 
                      \centering                
                      \begin{tabular}{||c | c||}     
                      \hline                                                
                      \includegraphics[width=6.0cm,height=5.0cm]{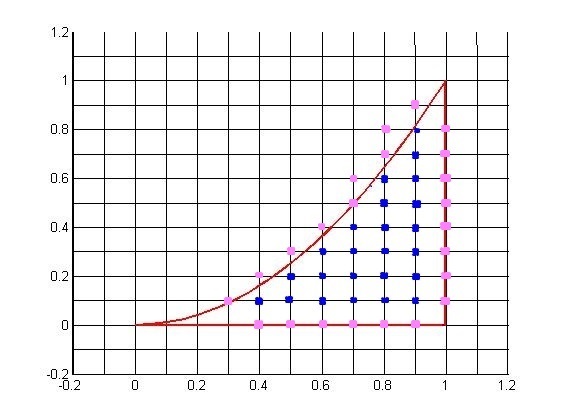} & \includegraphics[width=6.0cm,height=5.0cm]{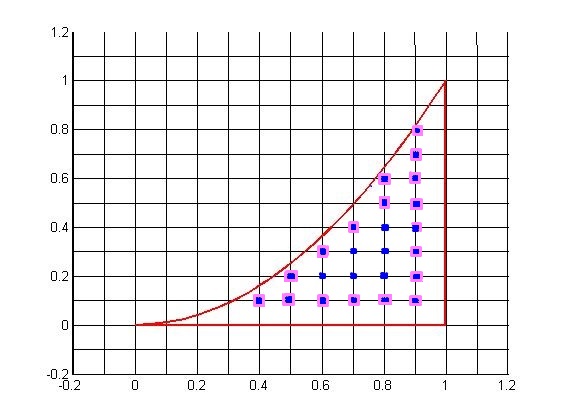} \\                  
                      \hline
                      \end{tabular}
                      \caption{Representaci\'{o}n de las mallas $\overline {\Omega}$ 
                               y $\underline{\Omega}$, que aproximan al dominio $\Omega$, 
                               respectivamente.}
                      \label{fig:pFronteraUtilizable}
                      \end{figure}   
                
                      De manera inmediata notamos que para resolver el problema
                      (\ref{Parabola:eq}) s\'{o}lo tenemos que aplicar el m\'{e}todo  
                      de diferencias finitas. Esto lo haremos en dos partes,
                      primero para $\overline {\Omega}$ y luego para $\underline{\Omega}$.
                      
                      Denotemos con $u_{i j}$ la soluci\'{o}n aproximada de $a_{i j}$
                      en el punto $(i \Delta x, j \Delta y) \in \overline {\Omega}$, 
                      entonces la discretizaci\'{o}n del problema (\ref{Parabola:eq})
                      est\'{a} dada por
                      \begin{equation}
                      \label{Parabola:disc}   
                             \displaystyle\frac{u_{i+1 j} - 2 u_{i j} + u_{i-1 j}}{\Delta x^2} + \displaystyle\frac{u_{i j+1} - 2 u_{i j} + u_{i j-1}}{\Delta y^2}  =  4,                       
                      \end{equation}                                
                      ver por ejemplo \cite{Jwt2}. Mientras 
                      \begin{equation}
                      \label{Parabola:fron} 
                             u_{i j}  =  f (i \Delta x, j \Delta y), \hspace*{0.3cm} \text{para} \hspace*{0.3cm} (i \Delta x, j \Delta y) \in   \partial  \overline {\Omega},
                      \end{equation}
                      denota la discretizaci\'{o}n de la frontera.      
                                                                                                              
                      Esto nos lleva a que tenemos que resolver un sistema de
                      ecuaciones lineales de la forma $Au=f$, donde $A$ es una
                      matriz invertible de $24\times24$. 
                      
                      Para resolver este sistema de ecuaciones lineales vamos
                      a aplicar el m\'{e}todo iterativo de Gauss-Seidel. 
                      
                      Una vez que hayamos encontrado una soluci\'{o}n aproximada 
                      estaremos interesados en saber que tan cercana est\'{a} 
                      la soluci\'{o}n aproximada $u$ de la soluci\'{o}n exacta $a$. 
                      Para esto vamos a calcular  
                      \begin{equation}
                       \label{errorOver}
                            error = \| a - u \|_ \infty. 
                      \end{equation}

                      \begin{table}[t!]
                       \centering
                      \begin{tabular}{|c|c|c|c|c|c|}
                      \hline
                      \multirow{2}{1cm}{\hspace*{0.3cm} $i$} & \multirow{2}{1cm}{\hspace*{0.3cm} $j$} & 
                      \multirow{2}{1.2cm}{$(x_i,y_j)$}  & \multirow{2}{1.6cm}{ $a(x_i,y_j)$} &
                      \multirow{2}{1.6cm}{$u(x_i,y_j)$} & \multirow{2}{1.6cm}{$v(x_i,y_j)$}   \\  
                            &           &            &            &             &           \\ \hline
                      \hline
                      
  $  4 $ & $  7 $ & $ (0.4000, 0.1000) $ & $ 0.2500 $ & $ 0.2500 $ &              \\

  $  4 $ & $  8 $ & $ (0.5000, 0.1000) $ & $ 0.3600 $ & $ 0.3600 $ &              \\

  $  4 $ & $  9 $ & $ (0.6000, 0.1000) $ & $ 0.4900 $ & $ 0.4900 $ &              \\

  $  4 $ & $ 10 $ & $ (0.7000, 0.1000) $ & $ 0.6400 $ & $ 0.6400 $ &              \\

  $  4 $ & $ 11 $ & $ (0.8000, 0.1000) $ & $ 0.8100 $ & $ 0.8100 $ &              \\

  $  4 $ & $ 12 $ & $ (0.9000, 0.1000) $ & $ 1.0000 $ & $ 1.0000 $ &              \\

  $  5 $ & $  8 $ & $ (0.5000, 0.2000) $ & $ 0.4900 $ & $ 0.4900 $ &              \\

  $  5 $ & $  9 $ & $ (0.6000, 0.2000) $ & $ 0.6400 $ & $ 0.6400 $ & $ 0.6400 $   \\

  $  5 $ & $ 10 $ & $ (0.7000, 0.2000) $ & $ 0.8100 $ & $ 0.8100 $ & $ 0.8100 $   \\

  $  5 $ & $ 11 $ & $ (0.8000, 0.2000) $ & $ 1.0000 $ & $ 1.0000 $ & $ 1.0000 $   \\

  $  5 $ & $ 12 $ & $ (0.9000, 0.2000) $ & $ 1.2100 $ & $ 1.2100 $ &              \\

  $  6 $ & $  9 $ & $ (0.6000, 0.3000) $ & $ 0.8100 $ & $ 0.8100 $ &              \\

  $  6 $ & $ 10 $ & $ (0.7000, 0.3000) $ & $ 1.0000 $ & $ 1.0000 $ & $ 1.0000 $   \\

  $  6 $ & $ 11 $ & $ (0.8000, 0.3000) $ & $ 1.2100 $ & $ 1.2100 $ & $ 1.2100 $   \\

  $  6 $ & $ 12 $ & $ (0.9000, 0.3000) $ & $ 1.4400 $ & $ 1.4400 $ &              \\

  $  7 $ & $ 10 $ & $ (0.7000, 0.4000) $ & $ 1.2100 $ & $ 1.2100 $ &              \\

  $  7 $ & $ 11 $ & $ (0.8000, 0.4000) $ & $ 1.4400 $ & $ 1.4400 $ & $ 1.4400 $   \\

  $  7 $ & $ 12 $ & $ (0.9000, 0.4000) $ & $ 1.6900 $ & $ 1.6900 $ &              \\

  $  8 $ & $ 11 $ & $ (0.8000, 0.5000) $ & $ 1.6900 $ & $ 1.6900 $ &              \\

  $  8 $ & $ 12 $ & $ (0.9000, 0.5000) $ & $ 1.9600 $ & $ 1.9600 $ &              \\

  $  9 $ & $ 11 $ & $ (0.8000, 0.6000) $ & $ 1.9600 $ & $ 1.9600 $ &              \\

  $  9 $ & $ 12 $ & $ (0.9000, 0.6000) $ & $ 2.2500 $ & $ 2.2500 $ &              \\

  $ 10 $ & $ 12 $ & $ (0.9000, 0.7000) $ & $ 2.5600 $ & $ 2.5600 $ &              \\

  $ 11 $ & $ 12 $ & $ (0.9000, 0.8000) $ & $ 2.8900 $ & $ 2.8900 $ &              \\   
                      
                      \hline
                      \end{tabular}
                      \caption{Comparaci\'{o}n de las soluciones de la ecuaci\'{o}n (\ref{Parabola:eq})
                               para $\Delta x =0.1 $ y $\Delta y = 0.1$ sobre $\overline {\Omega}$ 
                               y $\underline{\Omega}$.}  
                      \label{parabola:Valor1}  
                      \end{table}

                      \begin{sidewaysfigure}
                      \centering
                      \begin{tabular}{|c|c|c|}
                      \hline
                      \verb+ Solucion + & \verb+ Solucion  + &  \verb+ Solucion +  \\ 
                      \verb+ Exacta + & \verb+ por exterpolacion + & \verb+ por interpolacion + \\ \hline
                      \includegraphics[width=7cm,height=6cm]{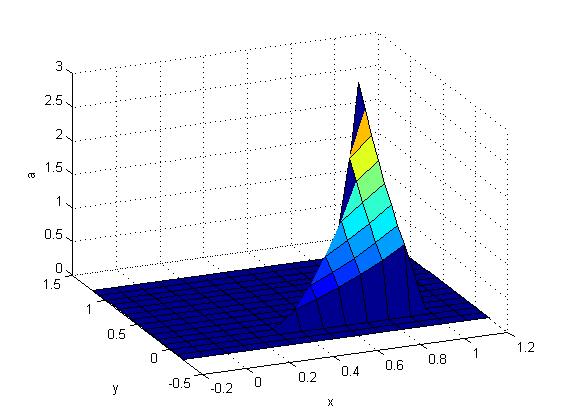}   & \includegraphics[width=7cm,height=6cm]{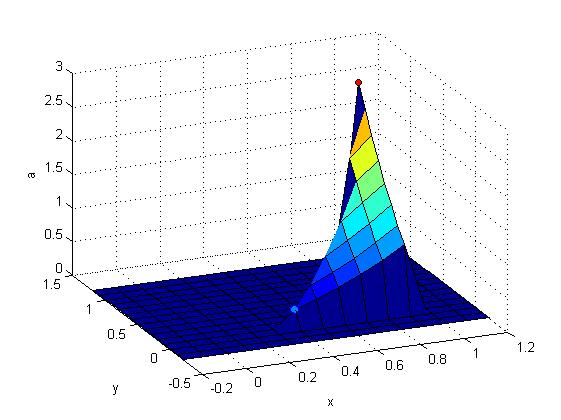} &\includegraphics[width=7cm,height=6cm]{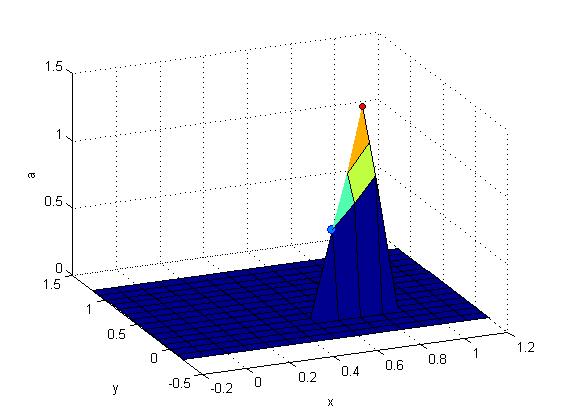} \\                  
                      \includegraphics[width=7cm,height=6cm]{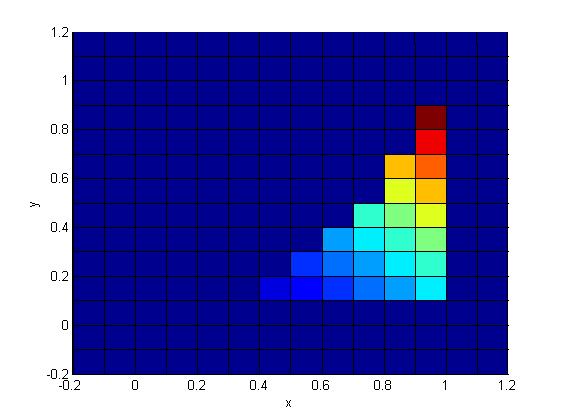}   & \includegraphics[width=7cm,height=6cm]{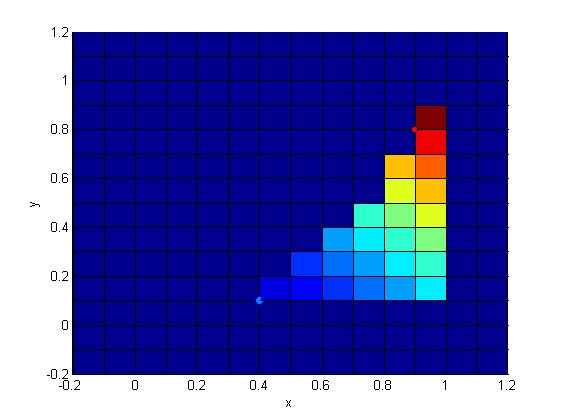} & \includegraphics[width=7cm,height=6cm]{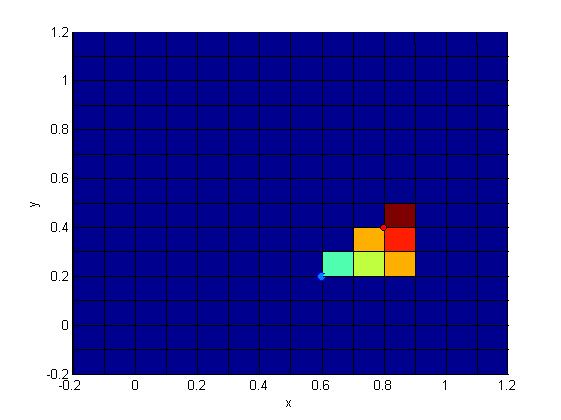} \\                             
                      \hline
                      \end{tabular}
                      \caption{Soluciones y sus proyecciones de la ecuaci\'{o}n (\ref{Parabola:eq}).}
                      \label{parabola:fig1}
                      \end{sidewaysfigure}
                      
                      En la tabla \ref{parabola:Valor1} se muestran los resultados 
                      obtenidos junto con los valores correctos, en \'{e}sta podemos
                      observar que el valor
                      \begin{equation*}
                       min = 0.2500 \hspace{0.3cm} \text{y} \hspace{0.3cm} max = 2.8900,                           
                      \end{equation*}
                      los cuales se localizan en
                      \begin{equation*}
                        (0.4000, 0.1000) \hspace{0.3cm} \text{y} \hspace{0.3cm} (0.9000, 0.8000),
                      \end{equation*}
                      respectivamente, y el $error = 8.8818e-16$. 
                                                     
                      La busqueda del $ min $ y del $ max $ \'{u}nicamente se 
                      llev\'{o} a cabo en el conjunto de puntos interiores 
                      de $\overline {\Omega}$. De esta manera evitamos que 
                      el $ min $ y el $ max $ no est\'{e}n en dominio $\Omega$, 
                      ya que hay puntos frontera de $\overline {\Omega}$ 
                      que no pertenecen al dominio $\Omega$. En particular 
                      notemos que el $ min $ y el $ max $ est\'{a}n en la 
                      frontera del interior de $\overline {\Omega}$. 
                      
                      Es importante mencionar que los resultados que se muestran 
                      en la tabla \ref{parabola:Valor1} son los resultados que se
                      muestran en {\it Matlab} usando formato {\it format short}, 
                      formato que ser\'{a} utilizado a lo largo de este documento.
                      
                      Ahora pasemos a resolver el problema (\ref{Parabola:eq}) 
                      sobre $\underline{\Omega}$, para esto denotemos con $v_{i j}$ 
                      la soluci\'{o}n aproximada de $a_{i j}$ en el punto
                      $(i \Delta x, j \Delta y) \in \underline{\Omega}$, entonces 
                      la discretizaci\'{o}n del problema (\ref{Parabola:eq})
                      es 
                      \begin{equation}
                      \label{Parabola:discv}   
                             \displaystyle\frac{v_{i+1 j} - 2 v_{i j} + v_{i-1 j}}{\Delta x^2} + \displaystyle\frac{v_{i j+1} - 2 v_{i j} + v_{i j-1}}{\Delta y^2}  =  4,                       
                      \end{equation}                                
                      ver por ejemplo \cite{Jwt2}. Mientras 
                      \begin{equation}
                      \label{Parabola:fronv} 
                             v_{i j}  =  f (i \Delta x, j \Delta y), \hspace*{0.3cm} \text{para} \hspace*{0.3cm} (i \Delta x, j \Delta y) \in \partial \underline{\Omega},
                      \end{equation}
                      denota la discretizaci\'{o}n de la frontera.                                             
                      
                      De manera inmediata observamos que esta discretizaci\'{o}n 
                      es an\'{a}loga a la discretizaci\'{o}n dada por (\ref{Parabola:disc})
                      y (\ref{Parabola:fron}). 
                      
                      Al igual que en el caso anterior, la discretizaci\'{o}n
                      (\ref{Parabola:discv}) y (\ref{Parabola:fronv}) nos lleva a
                      que tenemos que resolver un sistema de ecuaciones lineales
                      de la forma $Av=f$, s\'{o}lo que en este caso $A$ es una
                      matriz invertible de $6\times6$. 
                      
                      Antes de continuar observemos que no todos los puntos que  
                      se encuentran en la frontera se utilizan para determinar 
                      los valores desconocidos $v$, para ver \'{e}sto es suficiente  
                      con observar la segunda gr\'{a}fica de la figura 
                      \ref{fig:pFronteraUtilizable}.                      
                      
                      Para ser consistentes con lo que estamos haciendo vamos a 
                      resolver este sistema de ecuaciones lineales aplicando el
                      m\'{e}todo iterativo de Gauss-Seidel y para conocer que  
                      tan cercana est\'{a} la soluci\'{o}n aproximada $v$ de la
                      soluci\'{o}n exacta $a$, calculamos
                      \begin{equation*}
                            error = \| a - v \|_ \infty. 
                      \end{equation*}                                          
                      
                      Los resultados obtenidos se muestran en la \'{u}ltima columna
                      de la tabla \ref{parabola:Valor1}, en \'{e}sta podemos ver que
                      el valor
                      \begin{equation*}
                       min = 0.6400 \hspace{0.3cm} \text{y} \hspace{0.3cm} max = 1.4400,                           
                      \end{equation*}
                      los cuales se localizan en
                      \begin{equation*}
                        (0.6000, 0.2000) \hspace{0.3cm} \text{y} \hspace{0.3cm} (0.8000, 0.4000),
                      \end{equation*}
                      respectivamente, y el $error = 3.3307e-16$. 
                      
                      La busqueda del $ min $ y del $ max $, al igual que en el
                      caso anterior, se llev\'{o} acabo s\'{o}lo en los puntos
                      interiores de $\underline{\Omega}$, \'{e}sto con el fin de ser
                      consistentes con las comparaciones. Una vez m\'{a}s notemos
                      que el $ min $ y el $ max $ est\'{a}n en la frontera del 
                      interior de $\overline {\Omega}$.

                      Es importante mencionar que en este caso s\'{i} valdr\'{i}a 
                      la pena considerar tambi\'{e}n a la frontera de $\underline{\Omega}$
                      en la busqueda del $ min $ y del $ max $ ya que todos los
                      puntos de $\underline{\Omega}$ pertenecen al dominio $\Omega$.

                      De los resultados anteriores n\'{o}tese que el m\'{i}nimo y 
                      el m\'{a}ximo de $\overline {\Omega}$ y de $\underline{\Omega}$ 
                      tal vez podr\'{i}an no considerarse buenas aproximaciones a 
                      los valores exactos del dominio $\Omega$. Esto puede deberse 
                      a que 
                      \begin{enumerate}
                            \item $\Delta x $ y $\Delta y$ son muy grandes,
                                  
                            \item $\overline {\Omega}$ y $\underline{\Omega}$
                                  no cubre de manera adecuada al dominio $\Omega$,                                                          
                                  
                            \item el m\'{i}nimo y el m\'{a}ximo de $\Omega$ no 
                                  pertenecen a $\overline {\Omega}$ ni a $\underline{\Omega}$. 
                      \end{enumerate}
                      
                      Pero es importante notar que tanto estos valores como  
                      los dem\'{a}s valores obtenidos se aproximan muy bien 
                      a los valores exactos, este hecho puede corroborarse 
                      con tan s\'{o}lo comparar uno a uno los resultados de 
                      la tabla \ref{parabola:Valor1}. 
                      
                      Finalmente, en la figura \ref{parabola:fig1} se puede 
                      observar la representaci\'{o}n gr\'{a}fica de la 
                      soluci\'{o}n exacta evaluada en los puntos de la 
                      discretizaci\'{o}n de $\Omega$ y la soluci\'{o}n 
                      obtenida sobre $\overline{\Omega}$ y $\underline{\Omega}$ 
                      y sus respectivas proyecciones sobre el plano $xy$.                       
                      La escala del espectro de colores en {\it Matlab} que 
                      estamos usando en la gr\'{a}fica de la soluciones y 
                      que seguiremos usando de aqu\'{i} en adelante es como 
                      sigue. Con el color rojo indicaremos el valor m\'{a}ximo, 
                      mientras que con el color azul indicaremos el valor 
                      m\'{i}nino. 
                                                                                  
                \end{ejemplo}                                                               
                
                Hasta este momento hemos visto un ejemplo en donde resolvemos 
                una ecuaci\'{o}n diferencial parcial lineal sobre un dominio 
                irregular usando diferencias finitas. Para esto construimos 
                dos mallas que aproximan al dominio irregular con $\Delta x =0.1 $ 
                y $\Delta y = 0.1$, y es natural 
                preguntarnos \textquestiondown qu\'{e} sucede con el 
                comportamiento cualitativo de la soluci\'{o}n aproximada 
                cuando se eligen otros valores para $\Delta x $ y $\Delta y $?
                Para responder esta pregunta empecemos por analizar el
                siguiente ejemplo

                \begin{ejemplo} Consideremos de nuevo el problema del ejemplo anterior, el 
                       cual reescribimos aqu\'{i}
                      \begin{eqnarray}
                       \label{Parabola:eq2}
                            \nabla^2 a & = & 4,\hspace*{0.3cm}\text{para}\hspace*{0.3cm}(x,y)\in\Omega, \\
                                     a & = & (x+y)^2,\hspace*{0.3cm} \text{para} \hspace*{0.3cm} (x,y) \in   \partial  \Omega , \nonumber
                      \end{eqnarray}    
                      donde $a=a(x,y)$, $\Omega$ es el dominio determinado por 
                      la par\'{a}bola $y=x^2$, la recta $x=1$ y el eje horizontal
                      y $\partial \Omega$ es la frontera de $\Omega$.
                      Vamos a resolver este problema usando diferencias finitas
                      para diferentes valores de $\Delta x $ y $\Delta y $.                                              
                      
                      Recuerde que ya se hab\'{i}a mencionado que en el proceso de 
                      ir variando $\Delta x $ y $\Delta y $ se van a ir generando 
                      dos sucesiones de dominios, $\overline{\Omega}$ 
                      y $\underline{\Omega}$ y sobre cada uno de estos dominios 
                      vamos a resolver el problema (\ref{Parabola:eq2}) aplicando 
                      el mismo procedimiento que en el ejemplo anterior. 
                      
                      Despu\'{e}s de haber realizado varias pruebas, los resultados 
                      que observamos al resolver el problema (\ref{Parabola:eq2})
                      con diferentes valores de $\Delta x $ y $\Delta y $ son 
                      los siguientes.

                      \begin{table}[t!]
                       \centering
                      \begin{tabular}{|c|c|c|c|c|c|c|c|}
                      \hline
                      \multirow{2}{1cm}{\hspace*{0.1cm} $\Delta x$} 
                      & \multirow{2}{1.0cm}{ $u_{min}$}
                      & \multirow{2}{2.5cm}{ $(x_{min}, y_{min})$} 
                      & \multirow{2}{1.0cm}{ $u_{max}$}                       
                      & \multirow{2}{2.5cm}{ $(x_{max}, y_{max})$}                                                                    
                      & \multirow{2}{1.4cm}{ $Error$}  & \multirow{2}{1cm}{ $Iter$}\\  
                            &         &                        &            &             &  & \\ \hline
                      \hline
$ 0.1 $  & $ 0.2500 $ & $ (0.4000, 0.1000) $ & $ 2.8900 $ & $ (0.9000, 0.8000) $ & $ 8.8818e-16 $ & $ 94 $ \\                
$ 0.08 $ & $ 0.1600 $ & $ (0.3200, 0.0800) $ & $ 3.3856 $ & $ (0.9600, 0.8800) $ & $ 1.1102e-15 $ & $ 165 $ \\
$ 0.06 $ & $ 0.1296 $ & $ (0.3000, 0.0600) $ & $ 3.4596 $ & $ (0.9600, 0.9000) $ & $ 8.8818e-16 $ & $ 268 $ \\
$ 0.01 $ & $ 0.0144 $ & $ (0.1100, 0.0100) $ & $ 3.8809 $ & $ (0.9900, 0.9800) $ & $ 3.4417e-14 $ & $ 7917 $\\
$ 0.001 $ & $ 0.0011 $& $ (0.0320, 0.0010) $ & $ 3.9880 $ & $ (0.9990, 0.9980) $ & $ 2.6530e-12 $ & $678563$\\                                                                                                                                          
                      \hline
                      \end{tabular}
                      \caption{Resultados num\'{e}ricos de la ecuaci\'{o}n (\ref{Parabola:eq2})
                               para diferentes valores de $ \Delta x = \Delta y $ 
                               sobre $\overline {\Omega}$.}  
                      \label{Parabola:Valor2}  
                      \end{table}
                                                                  
                     
                                          
                                          

                     
                      Conforme $\Delta x $ y $\Delta y $ disminuyen
                      
                      \begin{enumerate}                                                                  
                            \item el dominio $\overline{\Omega}$ y $\underline{\Omega}$
                                  va adquiriendo diferentes formas, pero cada vez
                                  se va aproximando mejor al dominio exacto $\Omega$,
                      
                            \item las posiciones $(x_{min}, y_{min})$ y $(x_{max}, y_{max})$
                                  se aproximan mejor a las posiciones exactas,
                               
                            \item los valores $u_{min}$ y $u_{max}$ cada vez se 
                                  aproximan mejor a los valores exactos,
                                  
                            \item la sucesi\'{o}n generada por $u_{min}$ y $u_{max}$
                                  no es mon\'{o}tona creciente ni decreciente,
                                  
                            \item los valores obtenidos se aproximan cada vez mejor
                                  con los valores exactos,
                                  
                            \item el error disminuye, 
                            
                            \item la sucesi\'{o}n generada por el error
                                  no es mon\'{o}tona creciente ni decreciente,
                                
                            \item el n\'{u}mero de iteraciones aumenta y por ende
                                  el tiempo de ejecuci\'{o}n del programa tambi\'{e}n
                                  aumenta.                                                       
                      \end{enumerate}

                      Algunos de los resultados que obtuvimos se resumen en  
                      las tablas \ref{Parabola:Valor2} y \ref{Parabola:Valor3}, 
                      mientras que en la figura \ref{Parabola:difval} se 
                      pueden observar las gr\'{a}ficas de algunas soluciones 
                      obtenidas sobre el dominio $\overline{\Omega}$ y 
                      $\underline{\Omega}$, en donde con un punto color rojo
                      representamos el valor m\'{a}ximo y con un punto de color
                      azul representamos el valor m\'{i}nimo.
                      
                      N\'{o}tese que tanto en las tablas \ref{Parabola:Valor2}
                      y \ref{Parabola:Valor3} como en las gr\'{a}ficas de la figura
                      \ref{Parabola:difval} se pueden ver reflejadas las 
                      observaciones que acabamos de mencionar.                                                                  
                      
                      Por lo tanto cuando $\Delta x $ y $\Delta y $ tienden a  
                      ser suficientemente peque\~{n}os, entonces el tama\~{n}o 
                      del dominio $\overline{\Omega}$ y $\underline{\Omega}$ se 
                      va ajustando cada vez m\'{a}s y m\'{a}s al dominio $\Omega$,                      
                      \'{e}sto a su vez nos permite aproximarnos cada vez mejor
                      a la soluci\'{o}n con un grado de precisi\'{o}n razonable.

                      \begin{table}[t!]
                       \centering
                      \begin{tabular}{|c|c|c|c|c|c|c|c|}
                      \hline
                      \multirow{2}{1cm}{\hspace*{0.1cm} $\Delta x$} 
                      & \multirow{2}{1.0cm}{ $u_{min}$}
                      & \multirow{2}{2.5cm}{ $(x_{min}, y_{min})$} 
                      & \multirow{2}{1.0cm}{ $u_{max}$}                      
                      & \multirow{2}{2.5cm}{ $(x_{max}, y_{max})$}                                                                     
                      & \multirow{2}{1.4cm}{ $Error$}  & \multirow{2}{1cm}{ $Iter$}\\  
                            &         &                        &            &             &  & \\ \hline
                      \hline
$ 0.1  $  & $ 0.6400 $ & $ (0.6000, 0.2000) $ & $ 1.4400 $ & $ (0.8000, 0.4000) $ & $ 3.3307e-16 $ & $ 34 $ \\                
$ 0.08 $  & $ 0.5184 $ & $ (0.5600, 0.1600) $ & $ 2.0736 $ & $ (0.8800, 0.5600) $ & $ 6.6613e-16 $ & $ 80 $ \\
$ 0.06 $  & $ 0.3600 $ & $ (0.4800, 0.1200) $ & $ 2.4336 $ & $ (0.9000, 0.6600) $ & $ 8.8818e-16 $ & $ 155 $ \\
$ 0.01 $  & $ 0.0400 $ & $ (0.1800, 0.0200) $ & $ 3.6864 $ & $ (0.9800, 0.9400) $ & $ 3.0864e-14 $ & $ 7232 $\\
$ 0.001 $ & $ 0.0032 $ & $ (0.0550, 0.0020) $ & $ 3.9681 $ & $ (0.9980, 0.9940) $ & $ 2.6279e-12 $ & $675528$\\                                                                                                                                          
                      \hline
                      \end{tabular}
                      \caption{Resultados num\'{e}ricos de la ecuaci\'{o}n (\ref{Parabola:eq2})
                               para diferentes valores de $ \Delta x = \Delta y $ sobre
                               $\underline{\Omega}$.}  
                      \label{Parabola:Valor3}  
                      \end{table}
                     




                
                      \begin{figure}[h!]  
                      \centering
                      \begin{tabular}{||c  c||}
                      \hline                                   
                      \includegraphics[width=5.5cm,height=4.0cm]{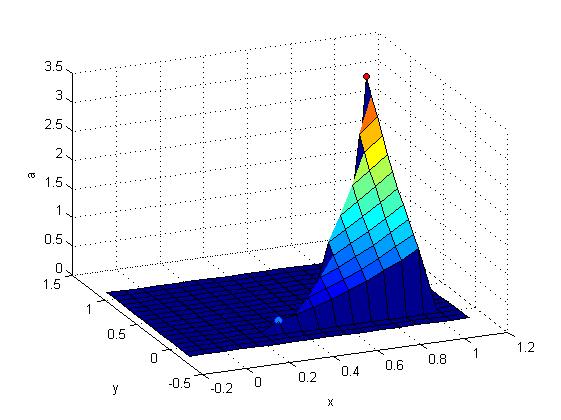} & \includegraphics[width=5.5cm,height=4.0cm]{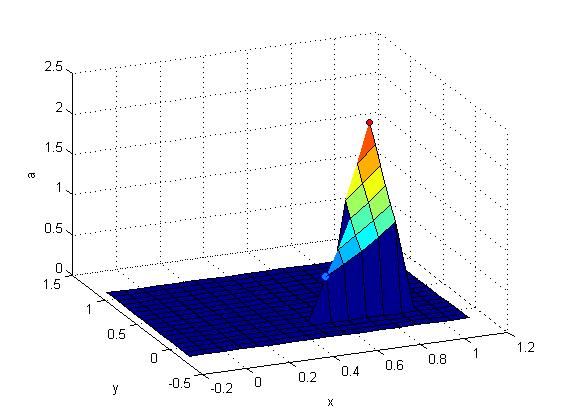}  \\                  
                      \multicolumn{2}{|| c ||}{ $\Delta x = 0.08 $ y $\Delta y = 0.08$ }\\  
                      \includegraphics[width=5.5cm,height=4.0cm]{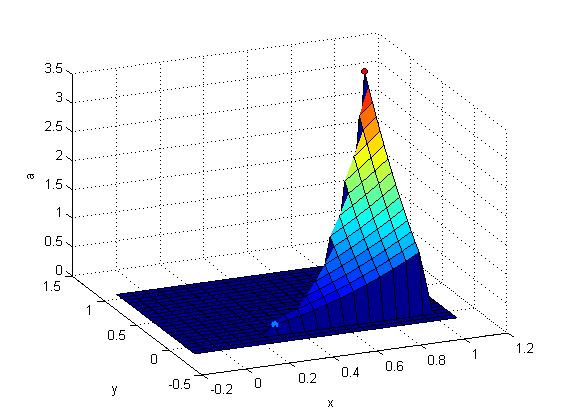} & \includegraphics[width=5.5cm,height=4.0cm]{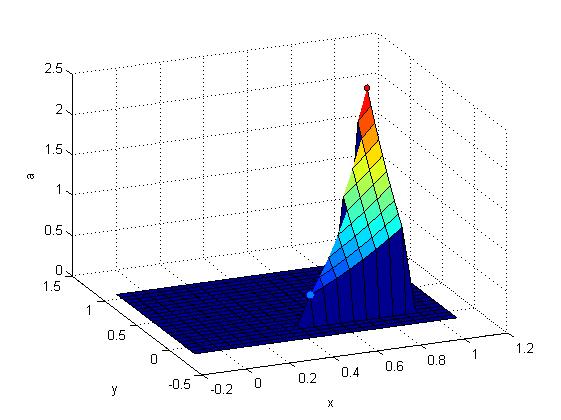} \\                
                      \multicolumn{2}{|| c ||}{ $\Delta x = 0.06 $ y $\Delta y = 0.06$ }\\
                      \includegraphics[width=5.5cm,height=4.0cm]{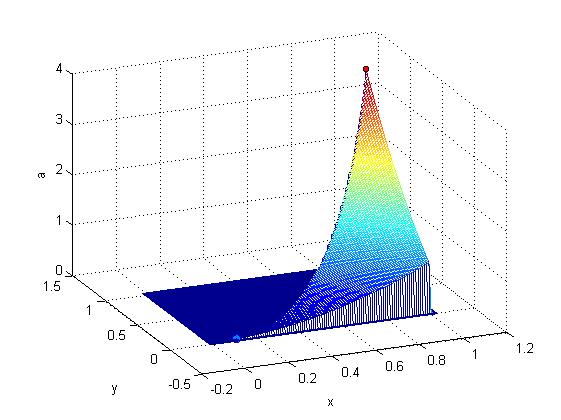} & 
                      \includegraphics[width=5.5cm,height=4.0cm]{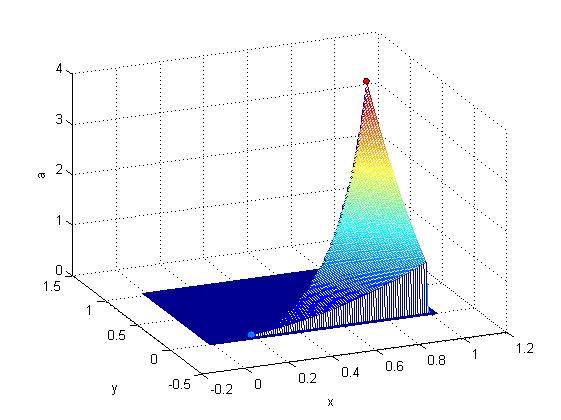} \\                
                      \multicolumn{2}{|| c ||}{ $\Delta x = 0.01 $ y $\Delta y = 0.01$ }\\
                      \includegraphics[width=5.5cm,height=4.0cm]{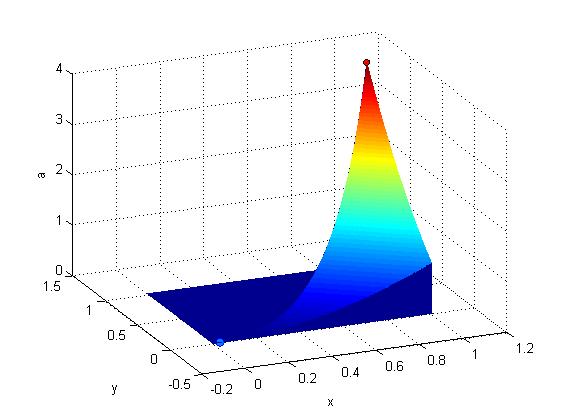} & \includegraphics[width=5.5cm,height=4.0cm]{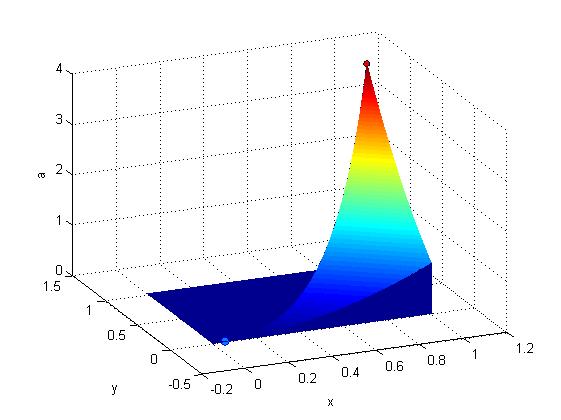} \\                 
                      \multicolumn{2}{|| c ||}{ $\Delta x = 0.001 $ y $\Delta y = 0.001$ }\\
                      \hline 
                      \end{tabular}
                      \caption{Soluci\'{o}n num\'{e}rica de la ecuaci\'{o}n (\ref{Parabola:eq2})
                               para diferentes valores de $\Delta x $ y $\Delta y $
                               sobre $\overline {\Omega}$ y $\underline{\Omega}$, respectivamente.}  
                      \label{Parabola:difval} 
                      \end{figure}

                \end{ejemplo}

\section{Problemas bien planteados}
                
                En este punto quizas el lector debe estar preguntandose
                {\it \textquestiondown qu\'{e} tipo de problemas se pueden 
                resolver con el m\'{e}todo presentado en la secci\'{o}n 
                anterior?} La respuesta depende de varios factores que
                iremos abordando en el desarrollo de este texto. Para
                empezar analicemos el siguiente ejemplo.                                
                
                Consideremos una f\'{a}brica de manufactura, en la que en  
                una etapa de la pro-ducci\'{o}n se deben de cortar placas 
                rectangulares de metal y despu\'{e}s cada placa debe de ser 
                expuesta al calor e ir monitoreando su temperatura en ciertos
                sitios. Una vez que se ha alcanzado la temperatura deseada 
                en dichos sitios, la placa pasar\'{a} a la siguiente etapa.
                
                \begin{figure}[t]
                \centering
                \begin{tabular}{||c||}
                \hline                
                \includegraphics[width=10cm,height=4cm]{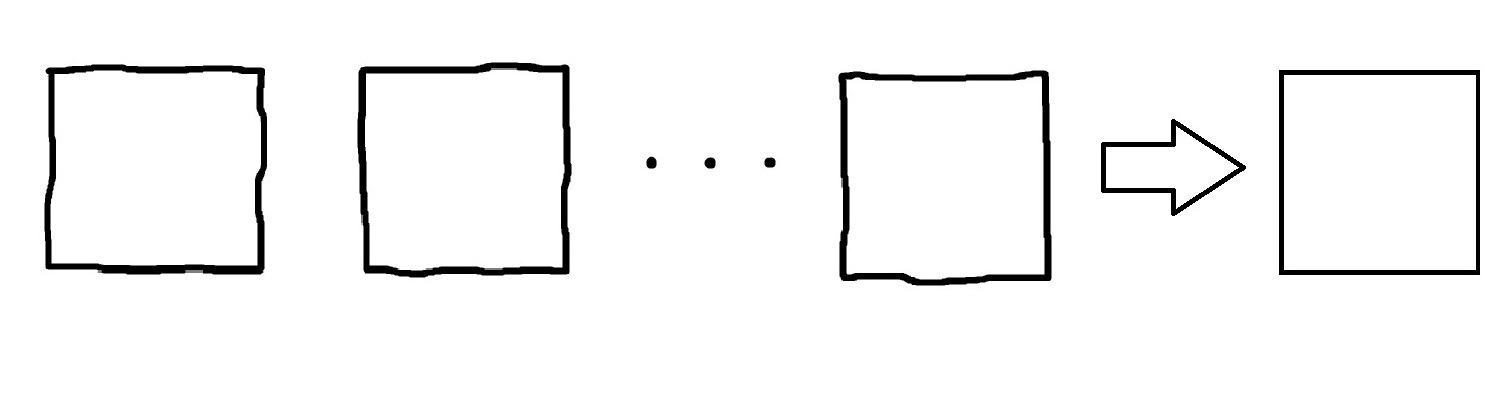} \\             
                \hline
                \end{tabular}
                \caption{La forma de las placas de metal son idealizadas 
                         en forma de un rect\'{a}ngulo.}
                \label{fig:placas}
                \end{figure}  
                
                El proceso de exponer una placa a una fuente de calor 
                puede ser modelado con el fin de c\'{a}lcular 
                \textquestiondown cu\'{a}nto tiempo debe de estar 
                expuesta la placa al calor?, \textquestiondown cu\'{a}ntas
                placas con la temperatura deseada se har\'{a}n en un d\'{i}a?,
                etc. Para dar respuesta a \'{e}stas u otras preguntas que
                puedan surgir, el tama\~{n}o de todas las placas es idealizado 
                en un rect\'{a}ngulo. Por otro lado, sabemos que, f\'{i}sicamente 
                las placas son diferentes entre s\'{i}, debido a que, por 
                ejemplo, los cortes no son exactos como se puede ver en la
                figura \ref{fig:placas}.                                                                                 
                
                M\'{a}s a\'{u}n, al exponer las placas al fuego y e ir 
                monitoreando la temperatura en un mismo sitio, f\'{i}sicamente 
                no es el mismo sitio, pero se considera que s\'{i} lo es, 
                debido a que ni las inexactitudes en el tama\~{n}o ni la 
                dilataci\'{o}n de cada placa afectan considerablemente la 
                posici\'{o}n del sitio donde se tomar\'{a} la temperatura.                                
                
                En este ejemplo resulta claro que a pesar de las diferencias
                entre las placas y de los efectos del proceso que sufren  
                \'{e}stas, el resultado de la producci\'{o}n no se ve afectado.
                M\'{a}s a\'{u}n, el modelo que se aplica para ayudarnos a 
                resolver los problemas ser\'{a} el mismo para todas las placas.                                                                                                                                           
                
                En la vida real los problemas no se presentan en dominios
                idealizados, como se acaba de explicar en el ejemplo anterior.
                Por lo tanto desde el punto de vista matem\'{a}tico tiene 
                sentido considerar problemas sobre dominios no idealizados.
                
\begin{figure}[t]
                \centering
                \begin{tabular}{||c||}
                \hline                
                \includegraphics[width=14cm,height=4cm]{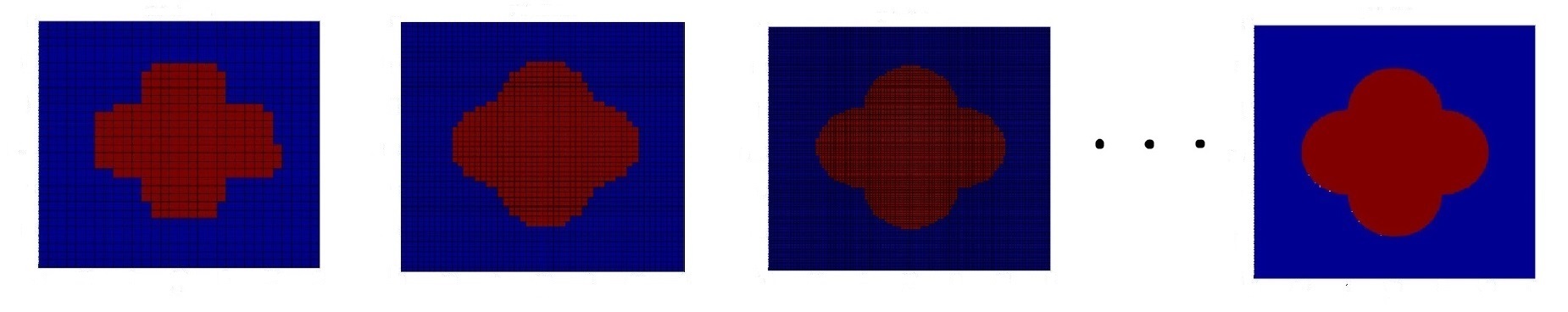} \\             
                \hline
                \end{tabular}
                \caption{Sucesi\'{o}n de dominios tal que convergen a una regi\'{o}n irregular.}
                \label{fig:flor_cua}
                \end{figure}

                Por lo que de aqu\'{i} en adelante vamos a considerar problemas
                sobre dominios idealizados tal que cada dominio idealizado
                tenga asociada o se le pueda asociar una sucesi\'{o}n de dominios
                que se aproxime y converga al mismo dominio idealizado. Esta 
                sucesi\'{o}n de dominios puede establecerse a partir de
                \begin{itemize}
                      \item las propiedades del problema \'{o}
                      \item mediante una estructura matem\'{a}tica,
                \end{itemize}
                como se puede ver en las figuras \ref{fig:placas}
                y \ref{fig:flor_cua}, respectivamente.
                
                Ahora pensemos en el m\'{e}todo presentado en la secci\'{o}n
                \ref{Metodo}, de acuerdo a lo anterior tiene sentido y es 
                razonable la forma en que estamos aplicando diferencias 
                finitas para resolver num\'{e}ricamente un problema sobre  
                un dominio irregular, mediante la construcci\'{o}n de una 
                sucesi\'{o}n de dominios (obtenida a partir de la 
                sucesi\'{o}n de mallas) y de soluciones que aproximen y 
                converjan al dominio original y a la soluci\'{o}n exacta, 
                respectivamente.                                                
                                
                Supongamos que tenemos que resolver un problema sobre un 
                dominio irregular aplicando diferencias finitas, para esto
                vamos a aplicar el m\'{e}todo que se explic\'{o} en la 
                secci\'{o}n \ref{Metodo}. Es decir vamos a sustituir el 
                dominio original por una sucesi\'{o}n de dominios\footnote{En
                algunos casos el dominio irregular se sustituye por un
                rect\'{a}ngulo o un c\'{i}rculo.} tal que esta sucesi\'{o}n 
                converga al dominio irregular y sobre cada uno de \'{e}stos  
                resolveremos el problema para obtener una sucesi\'{o}n de
                soluciones. Entonces es natural preguntarnos \textquestiondown qu\'{e}
                tanto cambiar\'{a} cada una de las soluciones aproximadas
                con respecto a la soluci\'{o}n exacta del problema?, 
                \textquestiondown cu\'{a}l soluci\'{o}n se debe elegir?,
                etc. Si bien, estas preguntas no son f\'{a}ciles de responder, 
                veamos algunos conceptos que debemos de tener en cuenta y                   
                que nos ser\'{a}n de gran utilidad para nuestro prop\'{o}sito.                                 
                
                De acuerdo a \cite{HeathNote} un problema est\'{a} 
                {\it bien planteado} si la soluci\'{o}n                                
                \begin{itemize}
                 \item existe, 
                 \item es \'{u}nica,
                 \item depende continuamente de los datos del problema. 
                \end{itemize}                                             
                
                De aqu\'{i} se sigue que si queremos que la sucesi\'{o}n de 
                soluciones, por extrapolaci\'{o}n e interpolaci\'{o}n, que 
                se obtiene al aplicar el m\'{e}todo la secci\'{o}n \ref{Metodo}
                a un problema no cambie de forma significativa con respecto
                a la soluci\'{o}n exacta, entonces lo recomendable es que el 
                problema dado sea un problema bien planteado. De esta forma
                podr\'{i}amos esperar que cada soluci\'{o}n de la sucesi\'{o}n 
                de soluciones va a ser una buena aproximaci\'{o}n a la 
                soluci\'{o}n exacta.
                
                M\'{a}s a\'{u}n, dado que nuestro objetivo es encontrar una
                aproximaci\'{o}n aceptable y precisa es importante verificar 
                que en el problema un cambio significativo en el dominio, u 
                alg\'{u}n otro cambio no afectan de forma considerable la 
                soluci\'{o}n. Sobre todo considerando el hecho de que vamos 
                a generar la sucesi\'{o}n de dominios y de soluciones mediante 
                algoritmos, como ya se ha visto.
                                                
                En la pr\'{a}ctica s\'{o}lo necesitaremos una buena 
                aproximaci\'{o}n al dominio y con ello generar una buena 
                aproximaci\'{o}n a la soluci\'{o}n del problema.

\section{Extensi\'{o}n y contracci\'{o}n de dominio}

                M\'{a}s all\'{a} de s\'{o}lo justificar desde un punto de
                vista pr\'{a}ctico, respaldado por un poco de teor\'{i}a, 
                la forma en que se aplica el m\'{e}todo de diferencias 
                finitas para resolver problemas sobre dominios irregulares,
                es necesario establecer los fundamentos que nos permitan
                contar con un sustento te\'{o}rico.
                
                \begin{figure}[t] 
                      \centering                
                      \begin{tabular}{||c | c||}     
                      \hline                                                
                      \includegraphics[width=6.0cm,height=5.0cm]{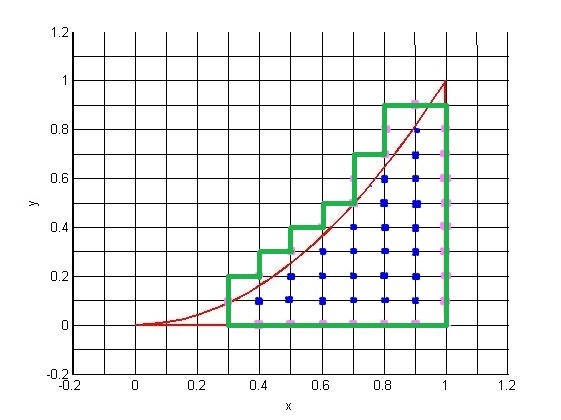} & \includegraphics[width=6.0cm,height=5.0cm]{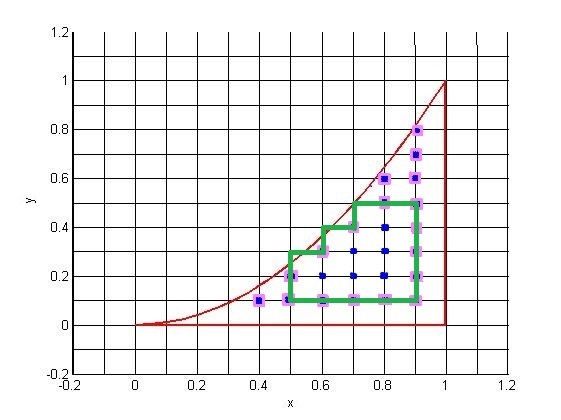} \\                  
                      \hline
                      \end{tabular}
                      \caption{\'{A}rea inducida por $\overline {\Omega}$ 
                               y $\underline{\Omega}$, respectivamente, tal 
                               que aproxima al \'{a}rea del dominio $\Omega$.}
                      \label{fig:areaFrontera}
                \end{figure}   
                
                Para hacer esto, empecemos por analizar el m\'{e}todo 
                presentado en la secci\'{o}n \ref{Metodo}. De este 
                m\'{e}todo sabemos que para construir una malla que 
                aproxime a $\Omega$ primero tenemos que construir un 
                rect\'{a}ngulo $R$ tal que $\Omega \subset R$. Luego
                debemos discretizar $R$ por una malla uniforme $G_R$
                como usualmente se hace e intersectar la malla $G_R$ 
                con el dominio $\Omega$ y a partir de \'{e}sto debemos
                determinar el conjuntos de puntos de $G_R$ que est\'{a}n
                en el interior de $\Omega$. Una vez que ya tenemos el 
                conjunto de puntos interiores, construimos la frontera
                ya sea por extrapolaci\'{o}n o por interpolaci\'{o}n.
                
                Por otra parte recordemos que, el primer punto del m\'{e}todo  
                de exhausti\'{o}n nos dice que de ser necesario, debemos de 
                {\it \textquotedblleft cuadrar\textquotedblright} el \'{a}rea
                mediante una sucesi\'{o}n de figuras cuyas \'{a}reas crecen 
                mon\'{o}tonamente.
                
                Entonces en este punto el lector debe estar pregunt\'{a}ndose                                
                \textquestiondown cu\'{a}l es el \'{a}rea de la sucesi\'{o}n 
                con la qu\'{e} estamos {\it \textquotedblleft cuadrando\textquotedblright} 
                el dominio irregular cuando aplicamos el m\'{e}todo de la 
                secci\'{o}n \ref{Metodo}? La respuesta es sencilla, basta 
                con unir los puntos que forman la frontera de los puntos 
                interiores aproximados de manera adecuada. En este caso 
                {\it podr\'{i}amos decir que estamos hablando de un \'{a}rea 
                inducida por los nodos interiores de un dominio dado}.   
                
                Un ejemplo de \'{e}sto puede verse en la figura 
                \ref{fig:areaFrontera}, en la cual la l\'{i}nea color verde
                determina el \'{a}rea que aproxima al \'{a}rea del dominio 
                $\Omega$ del ejemplo 1. Adem\'{a}s en ambas gr\'{a}ficas 
                podemos observar que puede haber puntos de la frontera del 
                \'{a}rea verde que pueden estar en la frontera del dominio 
                $\Omega$.
                
\begin{figure}[t]
                \centering
                \begin{tabular}{||c||}
                \hline                
                \includegraphics[width=9cm,height=4cm]{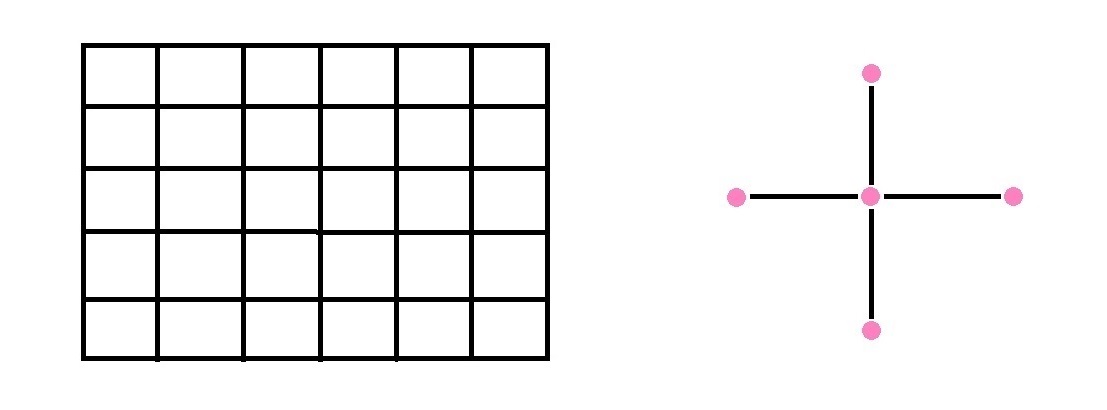} \\             
                \hline
                \end{tabular}
                \caption{Stencil.}
                \label{fig:stencil}
                \end{figure}

                Pero \textquestiondown por qu\'{e} se construye as\'{i} el 
                \'{a}rea? \'{o} lo que es equivalente \textquestiondown por 
                qu\'{e} no se aproxima el \'{a}rea como usualmente se 
                har\'{i}a? Por una parte ya hab\'{i}amos explicado que los
                nodos se eligen en base a la posici\'{o}n con respecto al 
                dominio $\Omega$ para tener mayor exactitud, ver el punto $3$
                del m\'{e}todo de la secci\'{o}n \ref{Metodo}. Pero adem\'{a}s
                es importante se\~{n}alar que la informaci\'{o}n num\'{e}rica
                necesaria para construir una aproximaci\'{o}n a la soluci\'{o}n 
                se encuentra en los puntos que forman el {\it stencil } y no
                en todos los nodos de la discretizaci\'{o}n del dominio $\Omega$.
                Un ejemplo de \'{e}sto se puede ver en la figura \ref{fig:stencil} 
                en donde todos los nodos proporcionan informaci\'{o}n 
                num\'{e}rica excepto {\it las esquinas}.
                
                Por lo que de aqu\'{i} en adelante cuando nos refieramos 
                a $\overline {\Omega}$ y a $\underline{\Omega}$ tambi\'{e}n 
                nos estaremos refiriendo a los puntos que se encuentran en  
                {\it las esquinas}, es decir a aquellos puntos que no 
                contribuyen con ninguna informaci\'{o}n num\'{e}rica, pero 
                que desde el punto de vista te\'{o}rico nos asegurar\'{a}n,
                bajo ciertas condiciones, la continuidad, la convergencia, etc.,
                tanto al dominio como a la soluci\'{o}n.                                                           
                                
                Ahora ya estamos en condiciones de establecer la primera
                hip\'{o}tesis. Para esto consideremos que tenemos un problema
                en un dominio irregular y que deseamos aplicar el m\'{e}todo 
                de diferencias finitas c\'{o}mo se explic\'{o} en la 
                secci\'{o}n \ref{Metodo}, en este caso por cada par de 
                valores $\Delta x $ y $\Delta y $ \'{o} $mx$ y $my$ 
                construiremos $\overline {\Omega}$ y $\underline{\Omega}$,                                                                                 
                \begin{center}
                      {\it de hecho suponemos que
                                    
                       \textbf{ cuando $\Delta x \rightarrow 0 $ 
                                y $\Delta y \rightarrow 0$, entonces
                                $\overline {\Omega}_{m,n} \rightarrow \Omega$ 
                                y $\underline{\Omega}_{m,n} \rightarrow \Omega$}, 
                                     
                                es decir obtenemos el 
                                dominio irregular $\Omega$},
                \end{center}                                
                observemos que impl\'{i}citamente estamos variando el dominio 
                del problema.
                                              
                Bajo esta suposici\'{o}n, que llamaremos {\it Hip\'{o}tesis $1$},
                estamos generando dos sucesiones de mallas uniformes 
                y estructuradas tales que convergen al dominio $\Omega$. De ser 
                as\'{i}, este hecho da lugar a que tendremos dos sucesiones 
                que aproximan a la soluci\'{o}n exacta del problema, como
                veremos a continuaci\'{o}n. 
                              
                A partir de lo anterior, junto con la suposici\'{o}n de la
                {\it hip\'{o}tesis $1$} nos permite establecer {\it la segunda 
                hip\'{o}tesis} de la siguiente manera,
                \begin{center}                              
                      {\it  cuando $\Delta x \rightarrow 0 $ 
                       y $\Delta y \rightarrow 0$, entonces
                                     
                       \textbf{la sucesi\'{o}n de soluciones generada por 
                               interpolaci\'{o}n y por extrapolaci\'{o}n sobre                                     
                               $\overline {\Omega}_{m,n}$ y $\underline{\Omega}_{m,n}$,
                               respectivamente, converge a la soluci\'{o}n exacta 
                               sobre el dominio $\Omega$,}
                                     
                               siempre que $\overline {\Omega}_{m,n} \rightarrow \Omega$ 
                               y $\underline{\Omega}_{m,n} \rightarrow \Omega$.}
                                                                                                            
                \end{center}                
                                                                                
                De aqu\'{i} se sigue que m\'{a}s que construir una
                sucesi\'{o}n de mallas $\overline {\Omega}_{m,n}$ 
                y $\underline{\Omega}_{m,n}$ tales que converjan al 
                dominio $\Omega$ lo ideal es que exista un entorno 
                donde el dominio se pueda extender \'{o} contraer 
                sin que la soluci\'{o}n cambi\'{e} dr\'{a}sticamente. 
                
                Antes de continuar es importante mencionar que aunque
                parece razonable la {\it hip\'{o}tesis $2$} y por tal 
                motivo fue enunciada tan deliberadamente, creemos que
                lo correcto ser\'{i}a afirmar que bajo ciertas condiciones,
                la {\it hip\'{o}tesis $2$} se satisfacer\'{a}. Pero
                creemos que \'{e}sto depender\'{a} del problema que 
                se este abordando.

                Decimos que {\it un problema es de dominio extendible
                \'{o} contra\'{i}ble si las propiedades cualitativas y 
                cuantitativas de la soluci\'{o}n no cambian o varian 
                muy poco}. Pero \textquestiondown qu\'{e} tanto debe
                extenderse o contraerse el dominio de un problema para 
                poder aplicar el m\'{e}todo de la secci\`{o}n \ref{Metodo}?,
                el dominio de un problema se debe extender y contraer 
                al menos lo suficientemente de tal manera que contenga
                cualquier malla $\overline {\Omega}_{m,n}$  
                y $\underline{\Omega}_{m,n}$. 
                \begin{ejemplo} Consideremos nuevamente el problema del primer
                       ejemplo
                      \begin{eqnarray}
                       \label{Parabola:eq3}
                            \nabla^2 a & = & 4,\hspace*{0.3cm}\text{para}\hspace*{0.3cm}(x,y)\in\Omega, \\
                                     a & = & (x+y)^2,\hspace*{0.3cm} \text{para} \hspace*{0.3cm} (x,y) \in   \partial  \Omega , \nonumber
                      \end{eqnarray}    
                      donde $a=a(x,y)$, $\Omega$ es el dominio determinado por 
                      la par\'{a}bola $y=x^2$, la recta $x=1$ y el eje horizontal
                      y $\partial \Omega$ es la frontera de $\Omega$.
                      
                      N\'{o}tese que \'{e}ste es un problema de dominio extendible y 
                      contra\'{i}ble ya que la soluci\'{o}n del problema no cambia
                      dr\'{a}sticamente dada cualquier aproximaci\'{o}n al dominio 
                      $\Omega$. De hecho para encontrar soluciones aproximadas no
                      fue necesario modificar \'{o} cambiar las condiciones de 
                      frontera. \textquestiondown puede el lector explicar este 
                      hecho?                      
                \end{ejemplo}
                
                En otro caso hubieramos tenido que redefinir el problema sobre 
                la frontera extendida, un ejemplo de como hacer \'{e}sto lo 
                podemos ver en \cite{Mattheij}.                             
                                                                                                 
                Por lo tanto cuando estemos considerando extender \'{o} contraer
                el dominio de un problema primero que nada debemos de analizar 
                si tiene sentido hacer \'{e}sto como se vio en la secci\'{o}n 
                anterior. Y al hacerlo debemos de cuidar no caer en singularidades 
                \'{o} en otros problemas que puedan alterar de manera dr\'{a}stica el 
                problema y por ende la soluci\'{o}n que deseamos aproximar. 
                                
                De esta manera tendremos la certeza que al aplicar el m\'{e}todo
                de diferencias finitas como se menciona en la secci\'{o}n \ref{Metodo}
                vamos a obtener una buena aproximaci\'{o}n a la soluci\'{o}n exacta.

\section{Construcci\'{o}n de $R$}                              
               
                En la secci\'{o}n \ref{Metodo}, se present\'{o} una forma de 
                discretizar un dominio irregular $\Omega$, lo cual nos permite
                resolver ecuaciones diferenciales parciales aplicando el 
                m\'{e}todo de diferencias finitas. Para \'{e}sto, lo primero 
                que tenemos que hacer es construir un rect\'{a}ngulo $R$ en el
                plano tal que $\Omega \subset R$, pero \textquestiondown c\'{o}mo
                debemos de construir $R$? Aunque la respuesta a esta pregunta 
                pareciera ser sencilla, es necesario justificarla.                                 
                
                Para esto consideremos una ecuaci\'{o}n \'{o} un sistema de 
                ecuaciones diferenciales parciales sobre $\Omega$ tal que la 
                soluci\'{o}n existe y es \'{u}nica. Definimos {\it la 
                soluci\'{o}n extendida de la ecuaci\'{o}n \'{o} del sistema
                de ecuaciones diferenciales sobre $R$ de la siguiente manera.
                Si $(x,y) \in \Omega$ le asociamos su soluci\'{o}n
                correspondiente y si $(x,y) \in R-\Omega$ le asociamos el cero}.
                
                Esta funci\'{o}n que acabamos de definir, en general no es  
                una funci\'{o}n continua, pero nos ser\'{a} \'{u}til para  
                lo que haremos a continuaci\'{o}n.                                

                Para ver que la soluci\'{o}n de la ecuaci\'{o}n \'{o} del 
                sistema de ecuaciones diferenciales no depende de la 
                elecci\'{o}n de $R$ consideremos dos rect\'{a}ngulos $R_1$ 
                y $R_2$. Entonces tenemos que ver que la soluci\'{o}n de la
                intersecci\'{o}n de $R_1$ y $R_2$ sigue siendo la soluci\'{o}n
                extendida.
                
                Sea $R_3 = R_1 \cap R_2$. Como $\Omega \subset R_1 $ y
                $\Omega \subset R_2 $, entonces para todo $(x,y) \in \Omega $
                se tiene que $(x,y) \in R_1 $ y $(x,y) \in R_2 $. Por lo 
                tanto $R_3$ es distinto del conjunto vac\'{i}o, ya que al 
                menos $\Omega \subset R_3$. 
                
                Ahora bien, como cada elemento de $\Omega$ es tambi\'{e}n 
                elemento de $R_1$ y $R_2$, entonces para todo elemento de $R_3$
                la soluci\'{o}n extendida se sigue satisfaciendo. De \'{e}sto 
                se sigue que 
                                
                \begin{center}
                      {\it                                      
                       \textbf{ la soluci\'{o}n de una ecuaci\'{o}n o de un sistema de
                                ecuaciones diferenciales parciales sobre un dominio
                                irregular no depende de la elecci\'{o}n de $R$.}                                      
                                }
                \end{center}
                
                Ahora que ya sabemos que la soluci\'{o}n de una ecuaci\'{o}n
                \'{o} de un sistema de ecuaciones es independiente de la 
                elecci\'{o}n de $R$ , entonces consideremos cualquier 
                rect\'{a}ngulo $R$ tal que $\Omega \subset R$.
                
                As\'{i} que para aproximarnos a la soluci\'{o}n deseada y                  
                de acuerdo a la secci\'{o}n \ref{Metodo}, tenemos que 
                discretiza $R$, es decir, tenemos que asociarle a $R$ una
                malla uniforme $G_R$, lo cual nos permitir\'{a} construir 
                $\overline {\Omega}_{m,n}$ y $\underline{\Omega}_{m,n}$.
                Por otra parte, hay que recordar que cuando estemos 
                construyendo $\overline {\Omega}_{m,n}$ 
                y $\underline{\Omega}_{m,n}$ necesitaremos extender y contraer
                el dominio $\Omega$ de tal forma que \'{e}ste contenga 
                cualquier $\overline {\Omega}_{m,n}$ y $\underline{\Omega}_{m,n}$.
                
                Por lo tanto el rect\'{a}ngulo $R$ se debe de construir de tal 
                forma que contenga cualquier extensi\'{o}n del dominio $\Omega$,
                \'{e}sto nos permitir\'{a} aplicar de manera adecuada el m\'{e}todo  
                de diferencias finitas \'{o} alg\'{u}n otro m\'{e}todo.                                                
                
                Dado que la soluci\'{o}n de una ecuaci\'{o}n o de un sistema 
                de ecuaciones no depende de la elecci\'{o}n de $R$, quizas,                
                en este punto el lector debe de estar pregunt\'{a}ndose
                \textquestiondown la soluci\'{o}n aproximada se ve afectada 
                por la discretizaci\'{o}n inducida por la elecci\'{o}n de $R$?                  
                Aunque hasta este momento no tenemos una respuesta a esta
                pregunta, es importante mencionar que cualquier intento a 
                \'{e}sta nos ha llevado a establecer nuestra {\it tercera                 
                hip\'{o}tesis.}                                                                                                                  
                
                \begin{center}
                      {\it Bajo ciertas condiciones y dependiendo de las 
                           propiedades cualitativas de la ecuaci\'{o}n \'{o} 
                           del sistema de ecuaciones tenemos que
                           
                       \textbf{la soluci\'{o}n aproximada sobre $\overline {\Omega}$ 
                               y $\underline{\Omega}$ no depende de la elecci\'{o}n 
                               de $R$ ni de la discretizaci\'{o}n inducida por \'{e}ste.}                                      
                                }
                \end{center}

\begin{figure}[h!]  
                \centering
                \begin{tabular}{||c  c||}
                \hline                                   
                \includegraphics[width=6.0cm,height=5.5cm]{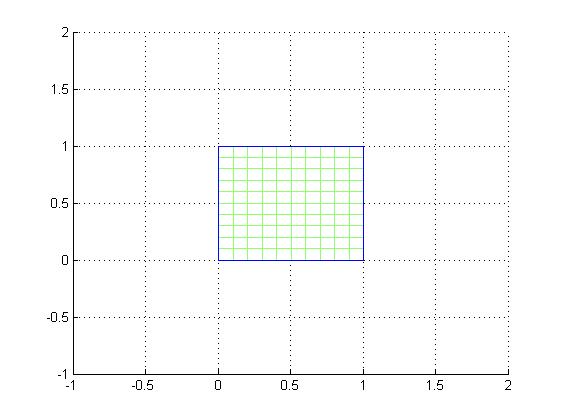} & \includegraphics[width=6.0cm,height=5.5cm]{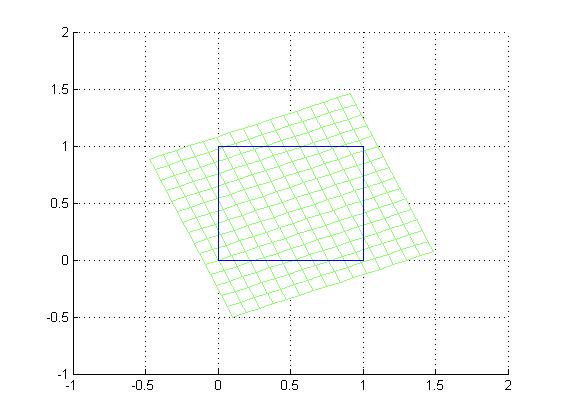} \\ 
                  Malla A & Malla B \\                   
                \includegraphics[width=6.0cm,height=5.5cm]{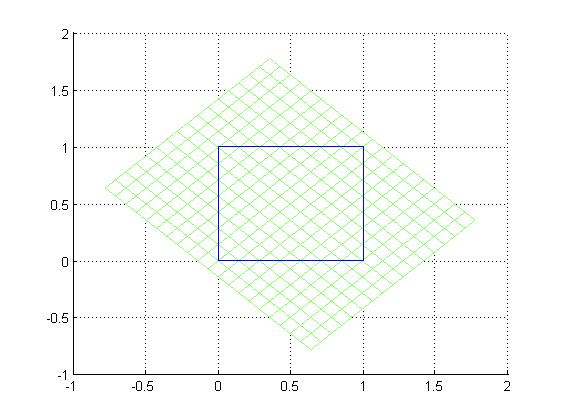} & \includegraphics[width=6.0cm,height=5.5cm]{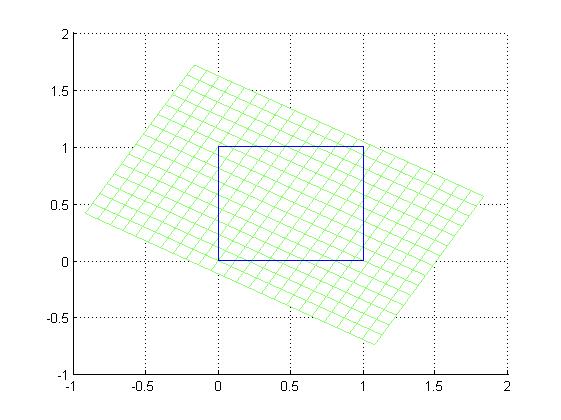} \\ 
                  Malla C & Malla D \\                                            
                \includegraphics[width=6.0cm,height=5.5cm]{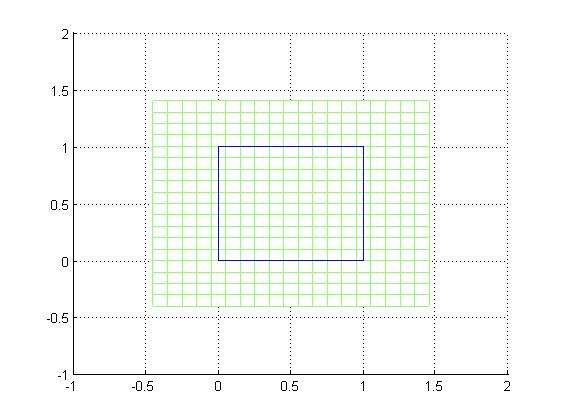} & \includegraphics[width=6.0cm,height=5.5cm]{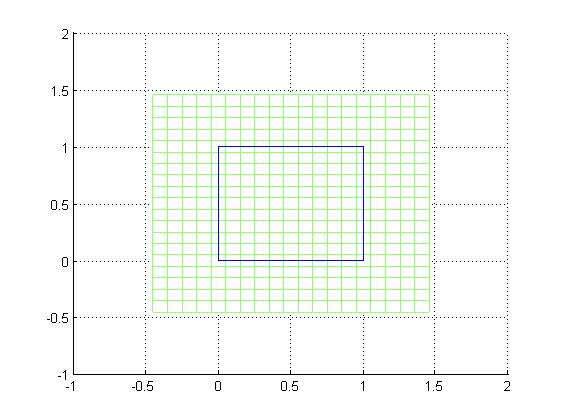} \\ 
                  Malla E & Malla F \\                                                             
                \hline
                \end{tabular}
                \caption{Diferentes elecciones de $R$ que nos permite asociar 
                         diferentes discretizaciones a la regi\'{o}n $\Omega$.}
                \label{mallas:fig1}
                \end{figure}

                \begin{ejemplo} \label{ejemplo4}  
                      Consideremos el caso cuando el dominio $\Omega$ es un cuadrado, 
                      es decir $\Omega=[a,b]\times[a,b]$ y supongamos que queremos
                      aplicar diferencias finitas para resolver un problema dado 
                      sobre $\Omega$. De manera breve, y sin menospreciar todo el 
                      an\'{a}lisis que se tiene que hacer, decimos que lo primero
                      que hacemos es discretizar el cuadrado $\Omega$ como usualmente
                      se acostumbra y despu\'{e}s pasamos a discretizar el problema 
                      aplicando diferencias finitas, esto nos permitir\'{a}  
                      obtener una soluci\'{o}n aproximada.
                      
                      Pero \textquestiondown est\'{a} es la \'{u}nica forma de 
                      asociar una discretizaci\'{o}n a el cuadrado $\Omega$ que 
                      nos permita aplicar diferencias finitas para encontrar una 
                      aproximaci\'{o}n a la soluci\'{o}n de un problema dado? 
                      La respuesta es no, ya que de acuerdo a la secci\'{o}n \ref{Metodo},
                      s\'{o}lo necesitamos construir $R$ tal que $\Omega \subset R$ 
                      y despu\'{e}s discretizar $R$ mediante una malla que llamamos
                      $G_R$. Esto nos permite construir dos discretizaciones 
                      $\overline {\Omega}$ y $\underline{\Omega}$ de $\Omega$. Por 
                      lo tanto los lados de $R$ no necesariamente deben de coincidir 
                      con los lados del cuadrado $\Omega$ para poder asociar una 
                      discretizaci\'{o}n a el cuadrado $\Omega$. Algunos ejemplos 
                      de \'{e}sto se pueden observar en la figura \ref{mallas:fig1}.
                 
                \end{ejemplo}

                Lo anterior nos lleva a preguntarnos \textquestiondown qu\'{e}
                tan bien funciona diferencias finitas para diferentes elecciones
                de $R$, considerando cualquier dominio $\Omega$? La respuesta, 
                como hemos venido diciendo depende del tipo de problema y de 
                sus propiedades.                                 
                                                
                Decimos que {\it un problema es de dominio completamente 
                extendible \'{o} contra\'{i}ble si para cualquier $R$ el
                problema es de dominio extendible \'{o} contra\'{i}ble.                
                Por otra parte decimos que el problema es de dominio
                condicionalmente extendible \'{o} contra\'{i}ble si al 
                menos existe un $R$ en donde el problema es de dominio 
                extendible \'{o} contra\'{i}ble. En otro caso decimos que 
                el problema no es de dominio extendible \'{o} contra\'{i}ble.}
                
                En cuanto a la soluci\'{o}n num\'{e}rica, creemos que en  
                un problema de dominio completamente extendible \'{o} 
                contra\'{i}ble, bajo ciertas condiciones, \'{e}sta no se 
                ve afectada por la elecci\'{o}n de $R$ ni por la 
                discretizaci\'{o}n inducida por \'{e}ste, como se puede  
                apreciar en el siguiente ejemplo.                                
                
                \begin{table}[t!] 
                      \centering                
                      \begin{tabular}{|c|c|c|c|c|c|c|}     
                      \hline    
                       Malla & A & B & C & D & E & F \\ \hline
                       Error & $ 0.0026 $ & $ 0.0026 $ & $ 0.0026 $ & $ 0.0026 $ & $ 0.0026 $ & $ 0.0026 $ \\                  
                      \hline
                      \end{tabular}
                      \caption{Error que se obtiene en cada soluci\'{o}n aproximada de la 
                               figura \ref{sol_mallaD}.}
                      \label{datosMalla}
                \end{table}
                                                                              
\begin{figure}[h!]  
                \centering
                \begin{tabular}{||c  c||}
                \hline                                   
                \includegraphics[width=6.0cm,height=5.5cm]{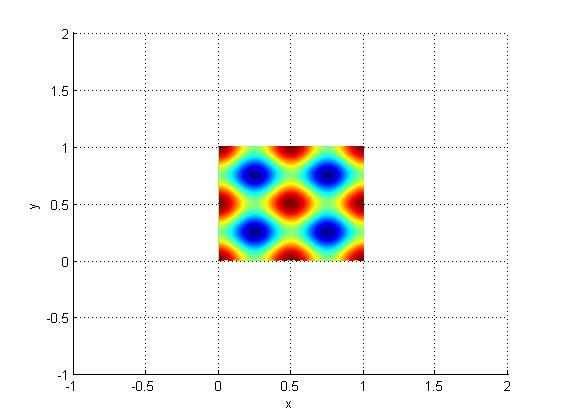} & \includegraphics[width=6.0cm,height=5.5cm]{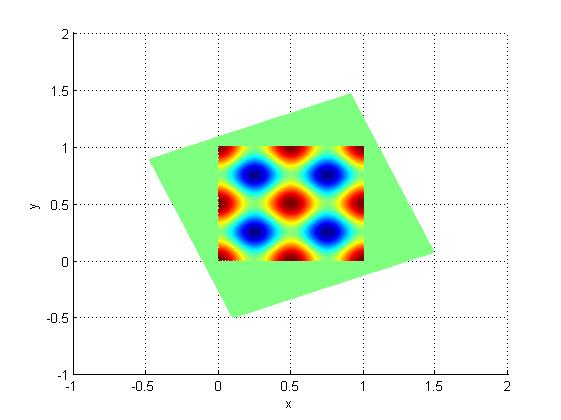} \\ 
                  Malla A & Malla B \\                     
                \includegraphics[width=6.0cm,height=5.5cm]{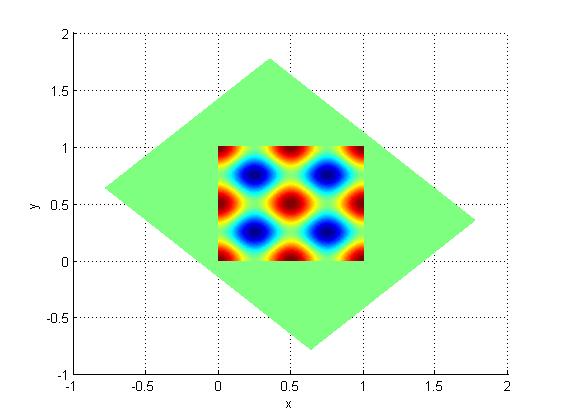} & \includegraphics[width=6.0cm,height=5.5cm]{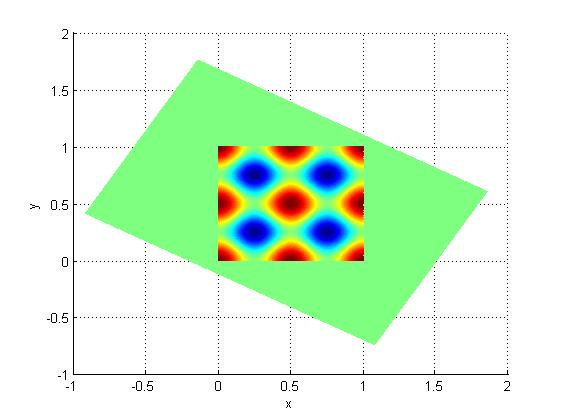} \\ 
                  Malla C & Malla D \\                                              
                \includegraphics[width=6.0cm,height=5.5cm]{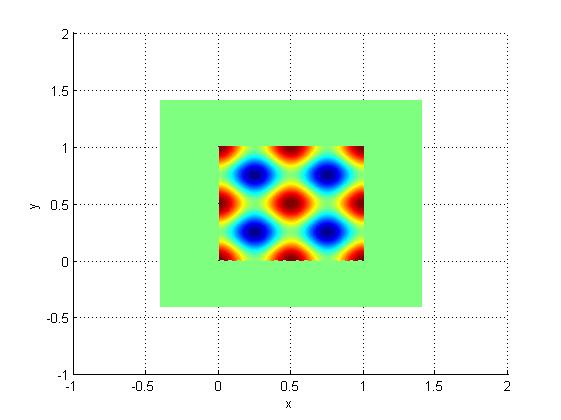} & \includegraphics[width=6.0cm,height=5.5cm]{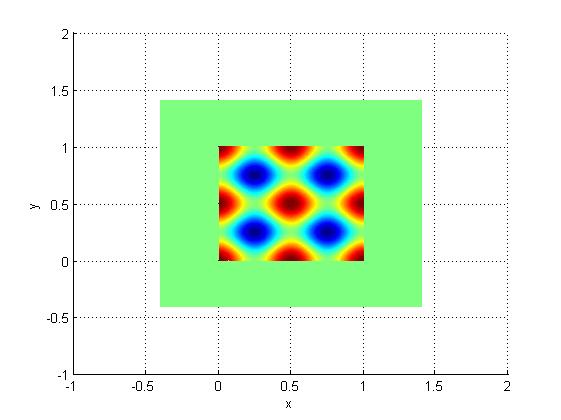} \\ 
                  Malla E & Malla F \\                                                               
                \hline
                \end{tabular}
                \caption{Soluci\'{o}n num\'{e}rica de la ecuaci\'{o}n (\ref{Cuadrado:elip})
                         sobre $\Omega$, con diferentes discretizaciones de $R$.}
                \label{sol_mallaD}
                \end{figure}

                \begin{ejemplo} \label{ejemplo5} Resolver la ecuaci\'{o}n diferencial 
                      con condiciones de Dirichlet en la frontera
                      \begin{eqnarray}
                       \label{Cuadrado:elip}
                            \nabla^2 a & = & -16 \pi ^2 [ cos(4 \pi x) + cos(4 \pi y)],\hspace*{0.3cm}\text{para}\hspace*{0.3cm}(x,y)\in\Omega, \\
                                     a & = & f(x,y),\hspace*{0.3cm} \text{para} \hspace*{0.3cm} (x,y) \in   \partial  \Omega , \nonumber
                      \end{eqnarray}
                      donde $a=a(x,y)$, $f(x,y)= cos(4 \pi x) + cos(4 \pi y)$, $\Omega$ es
                      el cuadrado unitario $(0,1)\times(0,1)$, y $\partial  \Omega$ es la 
                      frontera de $\Omega$. 
                      
                      \begin{enumerate}
                            \item Vamos a resolver este problema aplicando diferencias 
                                  finitas para cada uno de los diferentes rect\'{a}ngulos
                                  $R$ mostrados en el ejemplo \ref{ejemplo4}.
                            \item Despu\'{e}s vamos a calcular el error que se obtiene
                                  entre la soluci\'{o}n aproximada y la soluci\'{o}n 
                                  exacta, la cual es $ a(x,y) = -16 \pi^2 [cos(4\pi x) + cos(4 \pi y)]$.
                      \end{enumerate}
                                                                                                                                                                                                                                         
                      Para esto, primero construimos $\overline {\Omega}$ para cada
                      rect\'{a}ngulo $R$ considerando $\Delta x = 0.001 $ y $\Delta y = 0.001$ 
                      y luego discretizamos el problema (\ref{Cuadrado:elip}). 
                      Esto es, para cada $(i \Delta x,j \Delta y)\in \overline {\Omega}$                       
                      \begin{equation*}
                            \displaystyle\frac{u_{i+1 j} - 2 u_{i j} + u_{i-1 j}}{\Delta x^2} + \displaystyle\frac{u_{i j+1} - 2 u_{i j} + u_{i j-1}}{\Delta y^2}                       
                            = -16 \pi ^2 [ cos(4 \pi i \Delta x) + cos(4 \pi j \Delta y)],               
                      \end{equation*} 
                      mientras que para cada $(i \Delta x, j \Delta y) \in   \partial  \overline {\Omega}$,
                      \begin{equation*}                                                                    
                            u_{i j}  =  cos(4 \pi i \Delta x) + cos(4 \pi j \Delta y),                                       
                      \end{equation*} 
                      donde $ u_{i j} \approx a_{i j}$.
                      
                      Notemos que las expresiones anteriores forman un sistema de
                      ecuaciones lineales, por lo que para resolverlo aplicaremos
                      el m\'{e}todo iterativo de Jacobi.
                                            
                      En los resultados que obtuvimos, los cuales se pueden ver 
                      en la figura \ref{sol_mallaD}, se puede apreciar que cada 
                      soluci\'{o}n aproximada es muy parecida a la soluci\'{o}n 
                      exacta. M\'{a}s a\'{u}n, en la tabla \ref{datosMalla}
                      se puede observar que el error de aproximaci\'{o}n de cada
                      soluci\'{o}n aproximada es el mismo en todos los casos. 
                      
                      Por otra parte, resulta claro que la ecuaci\'{o}n con 
                      condiciones de frontera (\ref{Cuadrado:elip}) es un problema
                      de dominio completamente extendible, caracter\'{i}stica que 
                      reafirma que cada soluci\'{o}n obtenida sobre cada dominio 
                      extendido es una buena aproximaci\'{o}n a la soluci\'{o}n 
                      exacta.                                 
                \end{ejemplo}

                Es importante mencionar que hasta este momento no hemos 
                considerado otros factores que podr\'{i}an afectar la
                obtenci\'{o}n de la soluci\'{o}n num\'{e}rica sobre el 
                dominio extendido, tales como por ejemplo si el m\'{e}todo 
                de discretizaci\'{o}n es expl\'{i}cito \'{o} impl\'{i}cito.

\begin{figure}[h!]  
                \centering
                \begin{tabular}{||c  c||}
                \hline                                   
                \includegraphics[width=5.0cm,height=4.0cm]{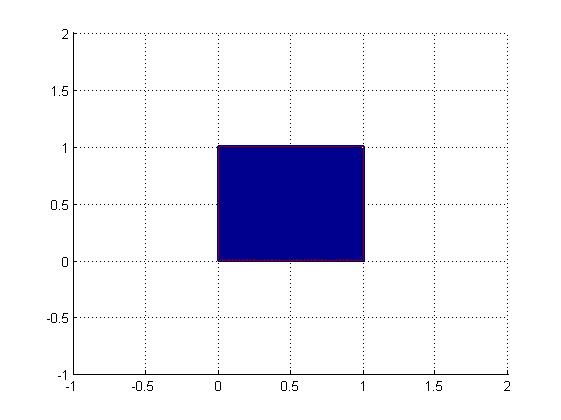}  & \includegraphics[width=5.0cm,height=4.0cm]{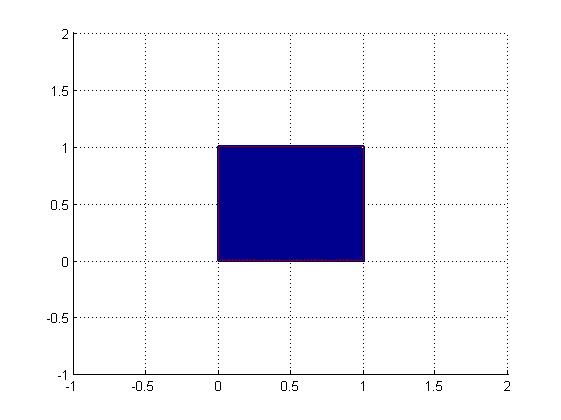} \\                                      
                  \multicolumn{2}{ ||c|| } {Condiciones iniciales} \\                  
                \includegraphics[width=5.0cm,height=4.0cm]{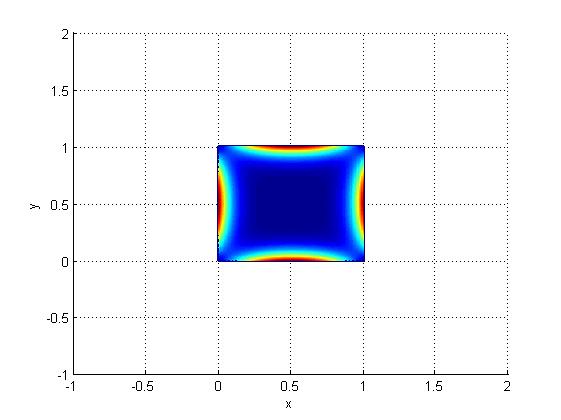}  & \includegraphics[width=5.0cm,height=4.0cm]{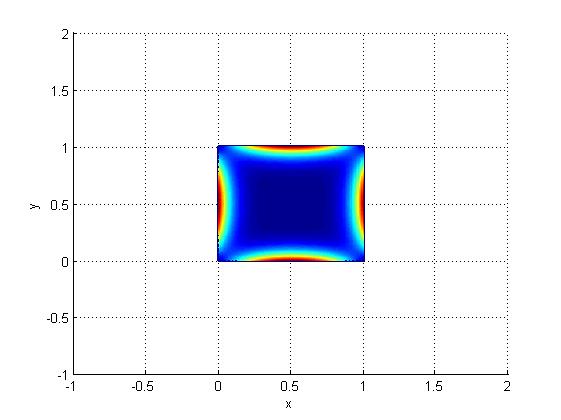} \\    
                  \multicolumn{2}{ ||c|| } { $ t = 0.0096 $ } \\                                                              
                \includegraphics[width=5.0cm,height=4.0cm]{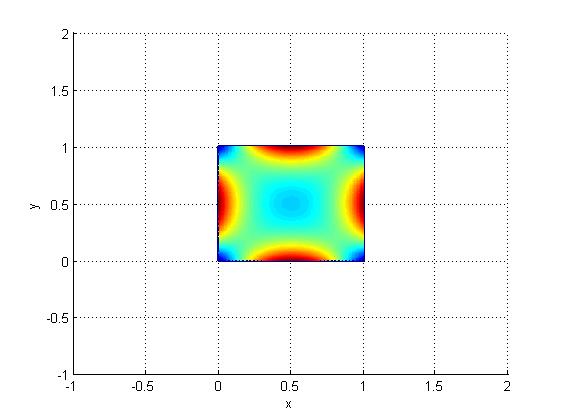}  & \includegraphics[width=5.0cm,height=4.0cm]{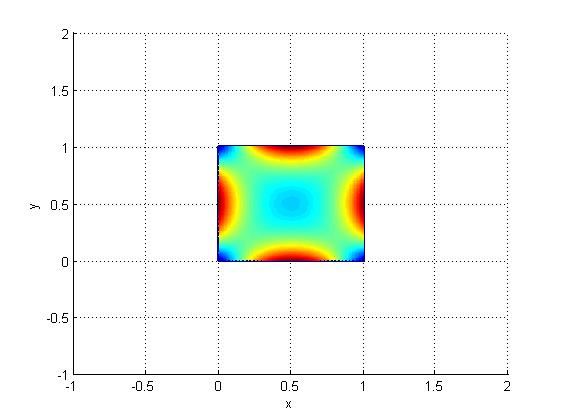} \\ 
                  \multicolumn{2}{ ||c|| } { $ t = 0.1052 $ } \\                    
                \includegraphics[width=5.0cm,height=4.0cm]{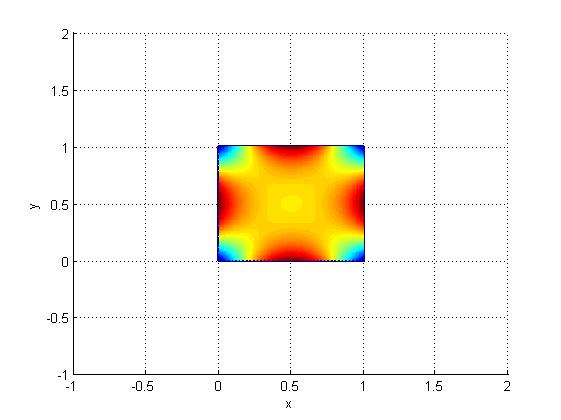}  & \includegraphics[width=5.0cm,height=4.0cm]{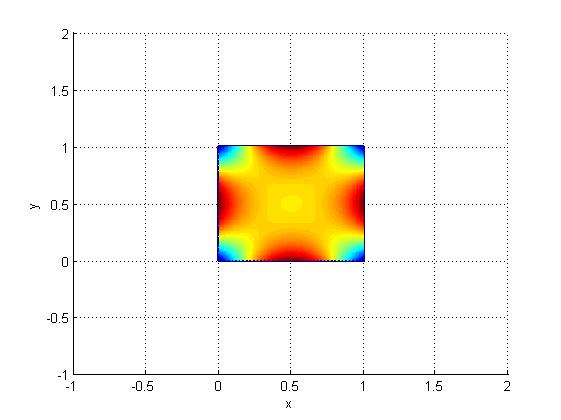} \\   
                  \multicolumn{2}{ ||c|| } { $ t = 0.3548 $ } \\                  
                \hline
                \end{tabular}
                \caption{Soluci\'{o}n num\'{e}rica de la ecuaci\'{o}n (\ref{Cuadrado:tiemLi})
                         sobre $\Omega$, con el m\'{e}todo expl\'{i}cito e impl\'{i}cito, de 
                         izquierda a derecha, respectivamente.}
                \label{sol_tiem_m1}
                \end{figure}

\begin{figure}[h!]  
                \centering
                \begin{tabular}{||c  c||}
                \hline                                   
                \includegraphics[width=5.0cm,height=4.0cm]{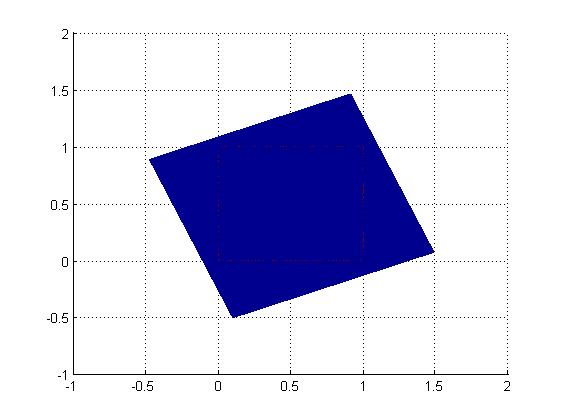} & \includegraphics[width=5.0cm,height=4.0cm]{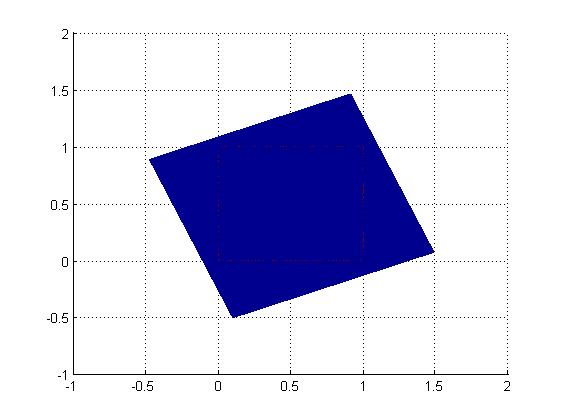} \\ 
                  \multicolumn{2}{ ||c|| } {Condiciones iniciales} \\                   
                \includegraphics[width=5.0cm,height=4.0cm]{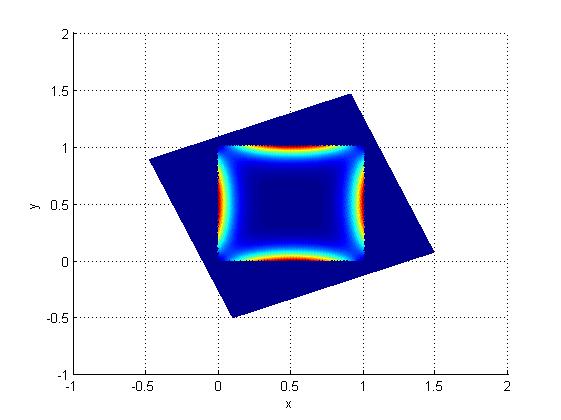} & \includegraphics[width=5.0cm,height=4.0cm]{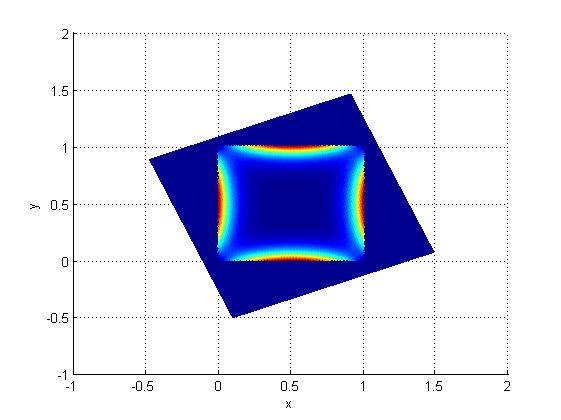} \\ 
                  \multicolumn{2}{ ||c|| } { $ t = 0.0096 $ } \\                                                             
                \includegraphics[width=5.0cm,height=4.0cm]{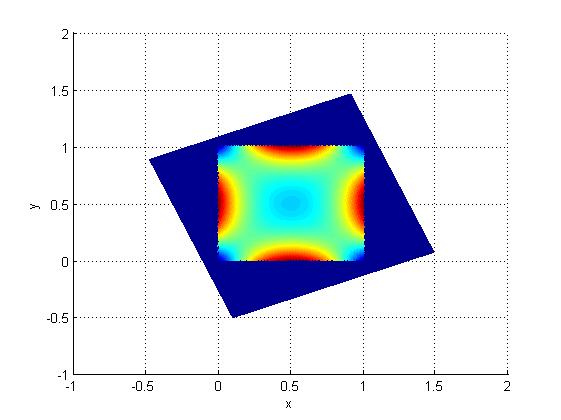} & \includegraphics[width=5.0cm,height=4.0cm]{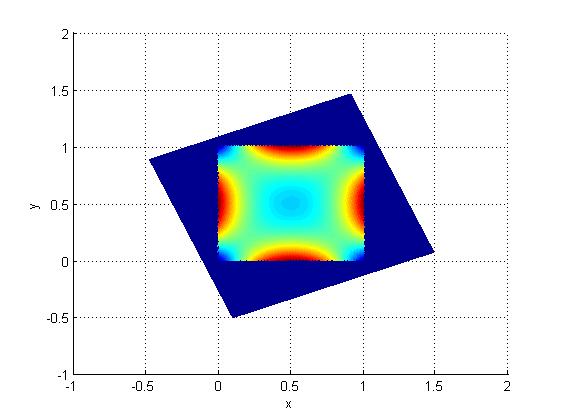} \\ 
                  \multicolumn{2}{ ||c|| } { $ t = 0.1052 $ } \\                    
                \includegraphics[width=5.0cm,height=4.0cm]{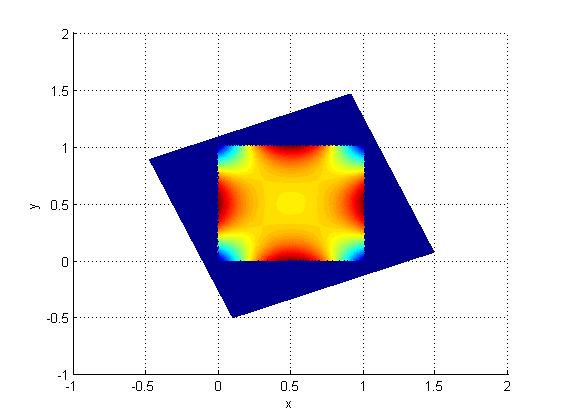} & \includegraphics[width=5.0cm,height=4.0cm]{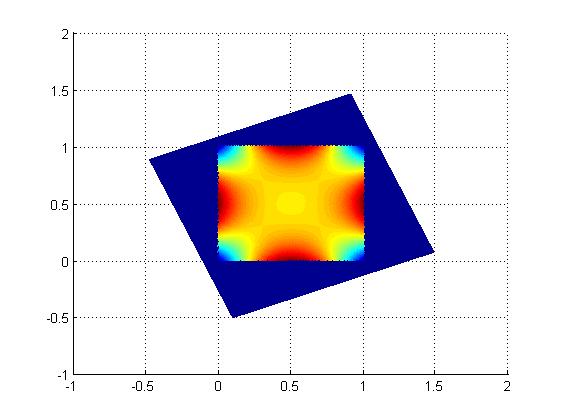} \\ 
                  \multicolumn{2}{ ||c|| } { $ t = 0.3548 $ } \\                    
                \hline
                \end{tabular}
                \caption{Soluci\'{o}n num\'{e}rica de la ecuaci\'{o}n (\ref{Cuadrado:tiemLi})
                         sobre $\Omega$, con el m\'{e}todo expl\'{i}cito e impl\'{i}cito, de 
                         izquierda a derecha, respectivamente.}
                \label{sol_tiem_m2}
                \end{figure}

\begin{figure}[h!]  
                \centering
                \begin{tabular}{||c  c||}
                \hline                                   
                \includegraphics[width=5.0cm,height=4.0cm]{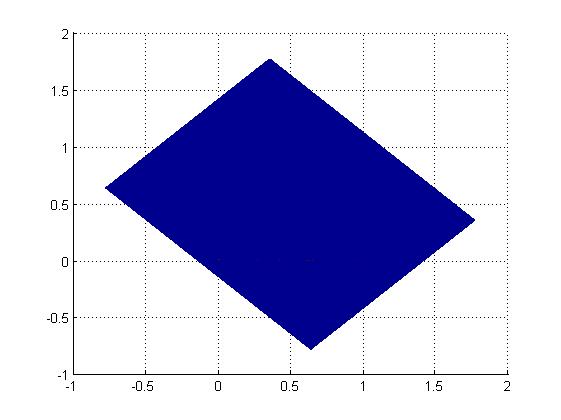} & \includegraphics[width=5.0cm,height=4.0cm]{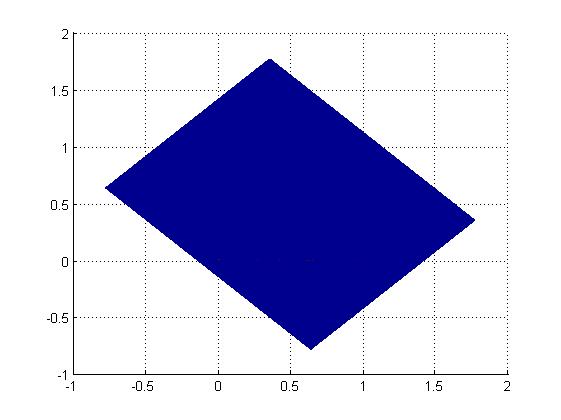} \\ 
                  \multicolumn{2}{ ||c|| } {Condiciones iniciales} \\                   
                \includegraphics[width=5.0cm,height=4.0cm]{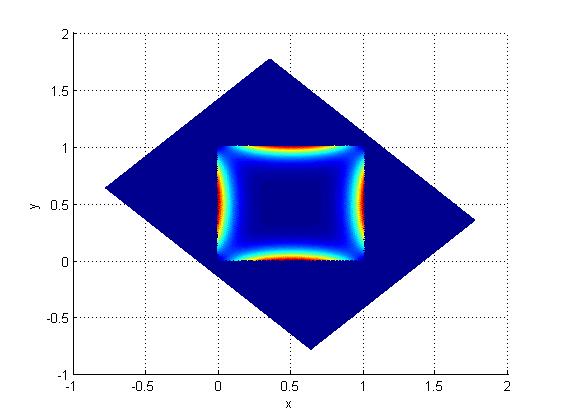} & \includegraphics[width=5.0cm,height=4.0cm]{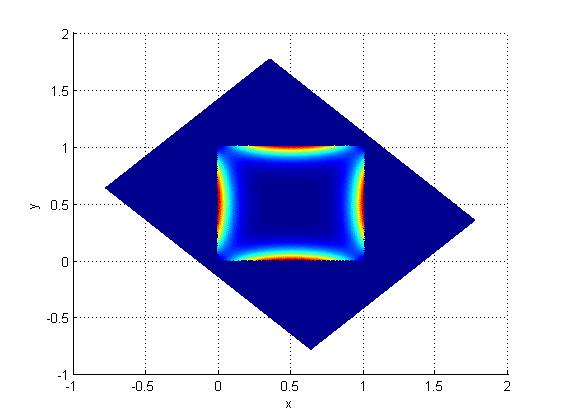} \\ 
                  \multicolumn{2}{ ||c|| } { $ t = 0.0096 $ } \\                                                             
                \includegraphics[width=5.0cm,height=4.0cm]{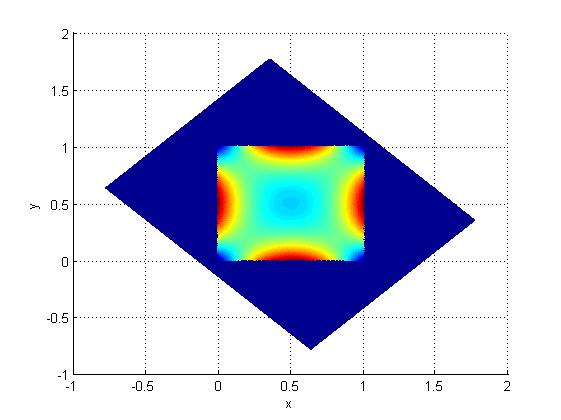} & \includegraphics[width=5.0cm,height=4.0cm]{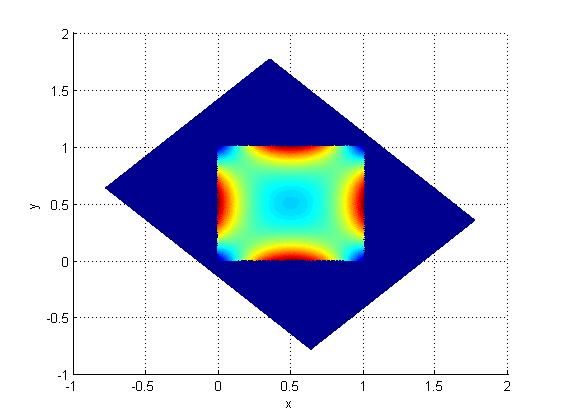} \\ 
                  \multicolumn{2}{ ||c|| } { $ t = 0.1052 $ } \\                    
                \includegraphics[width=5.0cm,height=4.0cm]{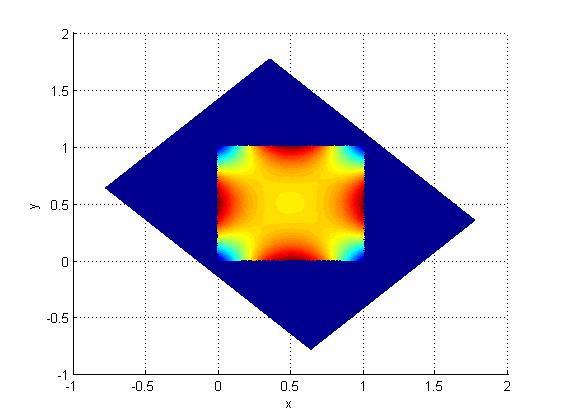} & \includegraphics[width=5.0cm,height=4.0cm]{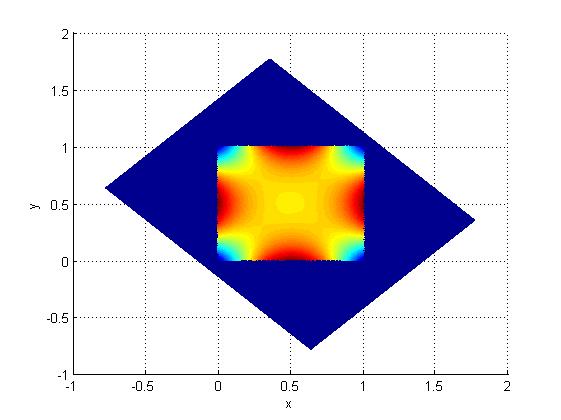} \\ 
                  \multicolumn{2}{ ||c|| } { $ t = 0.3548 $ } \\                    
                \hline
                \end{tabular}
                \caption{Soluci\'{o}n num\'{e}rica de la ecuaci\'{o}n (\ref{Cuadrado:tiemLi})
                         sobre $\Omega$, con el m\'{e}todo expl\'{i}cito e impl\'{i}cito, de 
                         izquierda a derecha, respectivamente.}
                \label{sol_tiem_m3}
                \end{figure}

\begin{figure}[h!]  
                \centering
                \begin{tabular}{||c  c||}
                \hline                                   
                \includegraphics[width=5.0cm,height=4.0cm]{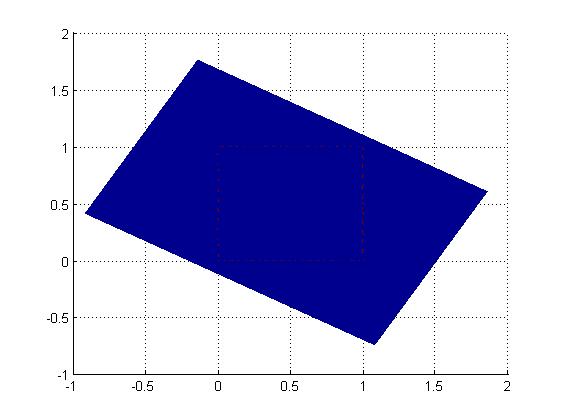} & \includegraphics[width=5.0cm,height=4.0cm]{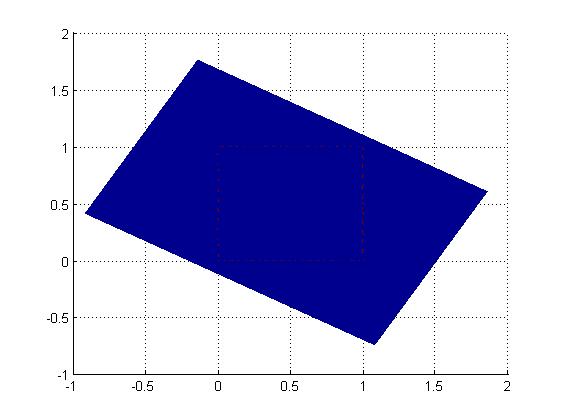} \\ 
                  \multicolumn{2}{ ||c|| } {Condiciones iniciales} \\                   
                \includegraphics[width=5.0cm,height=4.0cm]{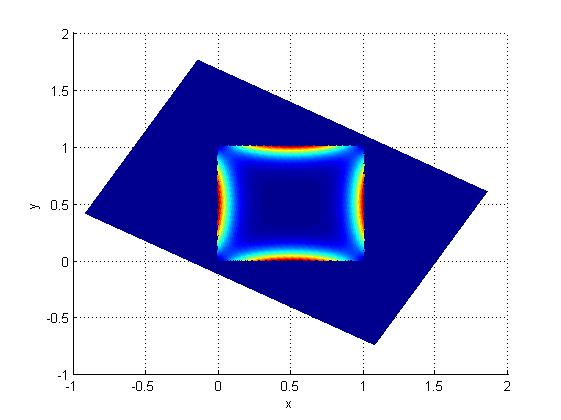} & \includegraphics[width=5.0cm,height=4.0cm]{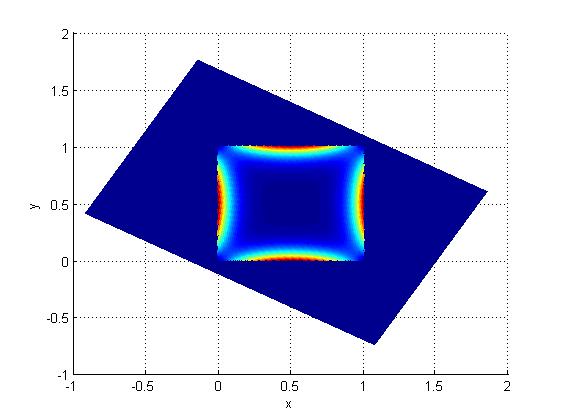} \\ 
                  \multicolumn{2}{ ||c|| } { $ t = 0.0096 $ } \\                                                             
                \includegraphics[width=5.0cm,height=4.0cm]{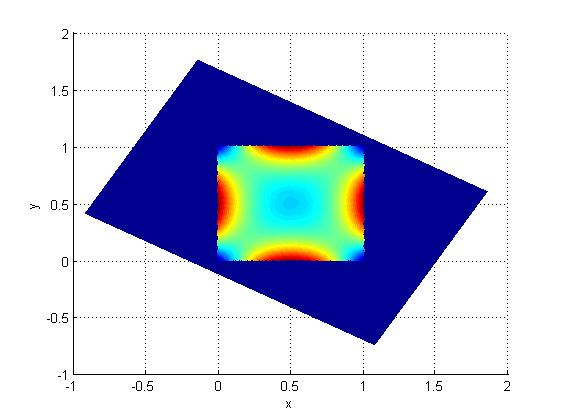} & \includegraphics[width=5.0cm,height=4.0cm]{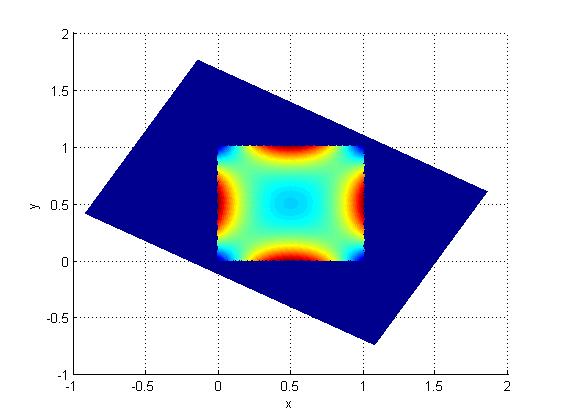} \\ 
                  \multicolumn{2}{ ||c|| } { $ t = 0.1052 $ } \\                    
                \includegraphics[width=5.0cm,height=4.0cm]{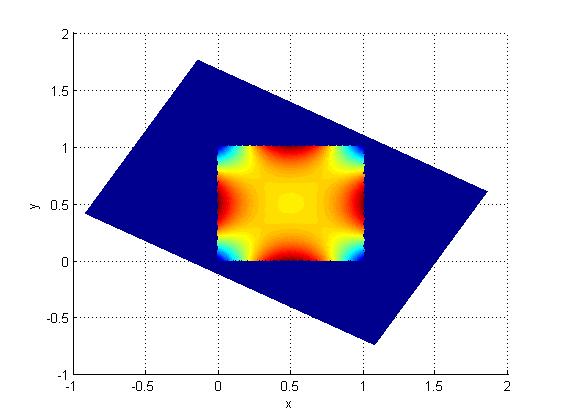} & \includegraphics[width=5.0cm,height=4.0cm]{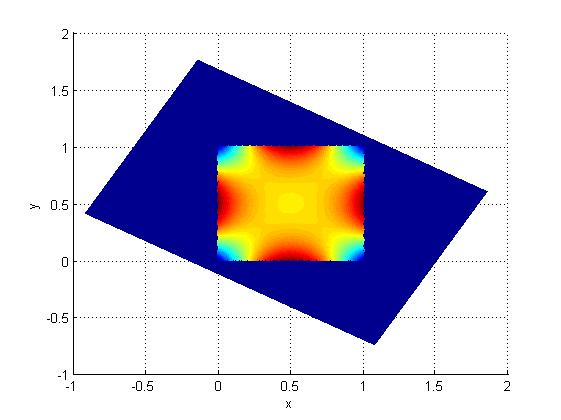} \\ 
                  \multicolumn{2}{ ||c|| } { $ t = 0.3548 $ } \\                    
                \hline
                \end{tabular}
                \caption{Soluci\'{o}n num\'{e}rica de la ecuaci\'{o}n (\ref{Cuadrado:tiemLi})
                         sobre $\Omega$, con el m\'{e}todo expl\'{i}cito e impl\'{i}cito, de 
                         izquierda a derecha, respectivamente.}
                \label{sol_tiem_m4}
                \end{figure}

\begin{figure}[h!]  
                \centering
                \begin{tabular}{||c  c||}
                \hline                                   
                \includegraphics[width=5.0cm,height=4.0cm]{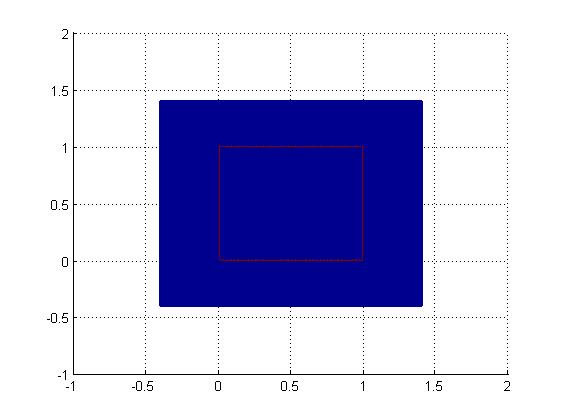} & \includegraphics[width=5.0cm,height=4.0cm]{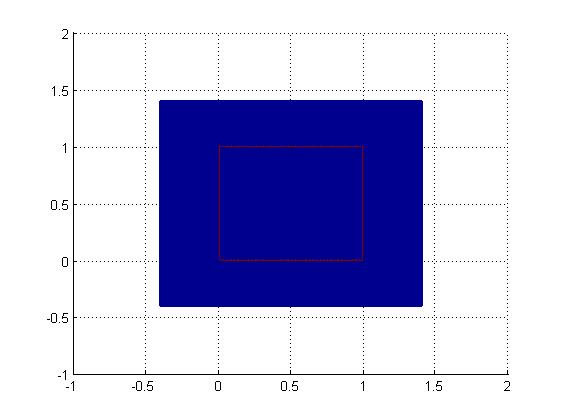} \\ 
                  \multicolumn{2}{ ||c|| } {Condiciones iniciales} \\                   
                \includegraphics[width=5.0cm,height=4.0cm]{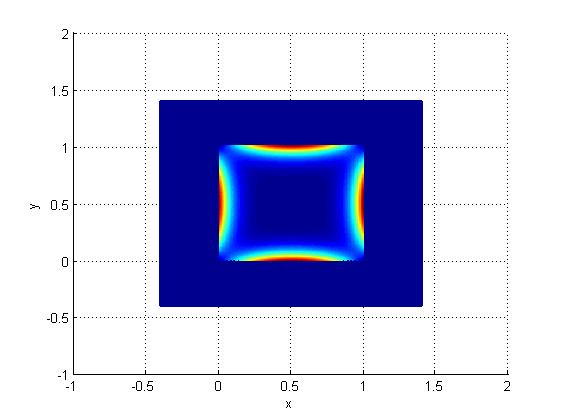} & \includegraphics[width=5.0cm,height=4.0cm]{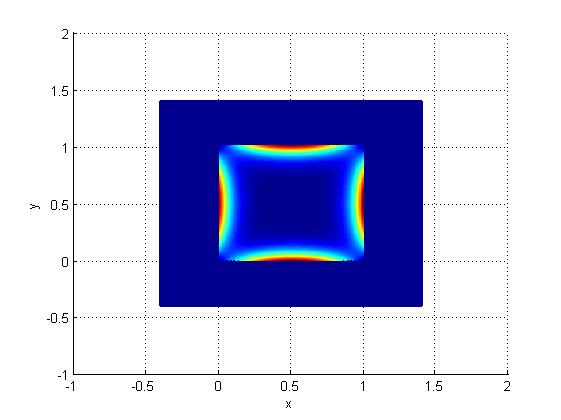} \\ 
                  \multicolumn{2}{ ||c|| } { $ t = 0.0096 $ } \\                                                             
                \includegraphics[width=5.0cm,height=4.0cm]{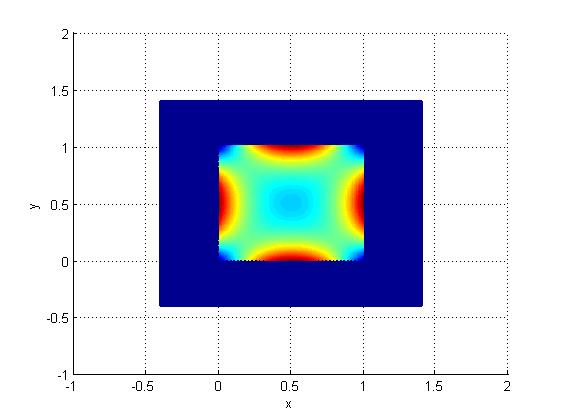} & \includegraphics[width=5.0cm,height=4.0cm]{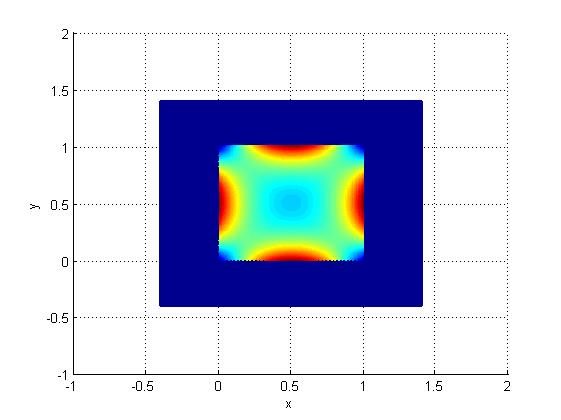} \\ 
                  \multicolumn{2}{ ||c|| } { $ t = 0.1052 $ } \\                    
                \includegraphics[width=5.0cm,height=4.0cm]{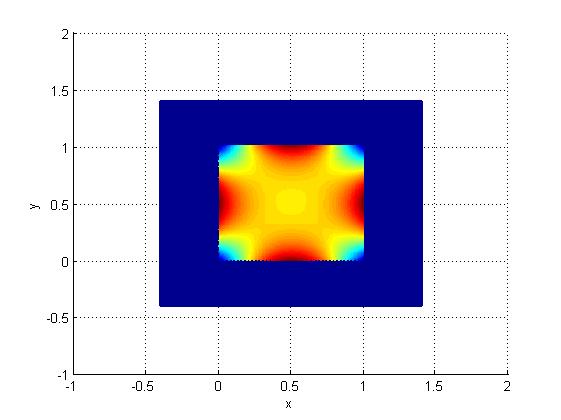} & \includegraphics[width=5.0cm,height=4.0cm]{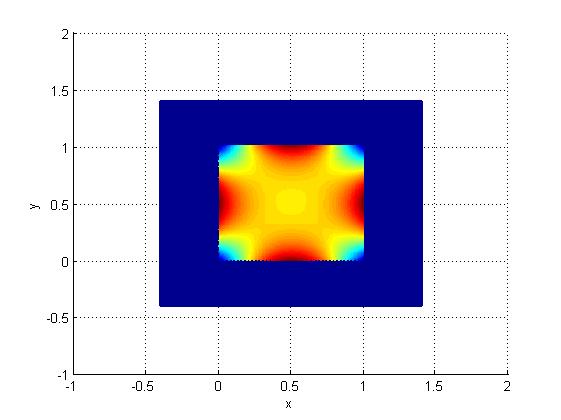} \\ 
                  \multicolumn{2}{ ||c|| } { $ t = 0.3548 $ } \\                    
                \hline
                \end{tabular}
                \caption{Soluci\'{o}n num\'{e}rica de la ecuaci\'{o}n (\ref{Cuadrado:tiemLi})
                         sobre $\Omega$, con el m\'{e}todo expl\'{i}cito e impl\'{i}cito, de 
                         izquierda a derecha, respectivamente.}
                \label{sol_tiem_m5}
                \end{figure}

                \begin{ejemplo} \label{ejemplo6} Consideremos el siguiente problema
                      con condiciones de Dirichlet en la frontera
                      \begin{eqnarray}
                       \label{Cuadrado:tiemLi}
                            a_t & = & \nabla^2 a,\hspace*{0.3cm}\text{para}\hspace*{0.3cm}(x,y)\in\Omega,\hspace*{0.3cm} t>t_0, \\
                                     a(x,y,t) & = & g(x,y,t),\hspace*{0.3cm} \text{para} \hspace*{0.3cm} (x,y,t) \in   \partial  \Omega,\hspace*{0.3cm} t>t_0 , \nonumber \\
                                     a(x,y,t_0) & = & f(x,y,t_0) ,\hspace*{0.3cm} \text{para} \hspace*{0.3cm} (x,y) \in \Omega \cup \partial  \Omega, \nonumber
                      \end{eqnarray}                      
                      donde $a=a(x,y,t)$, $g(x,y,t)=sin(t)[sin( \pi x) + sin(\pi y)]$, $f(x,y,t_0) =0$,
                      $t_0=0$, $\Omega$ es el cuadrado unitario $(0,1)\times(0,1)$ y $\partial  \Omega$ 
                      es la frontera de $\Omega$. 
                        
                      \begin{enumerate}                       
                       \item Al igual que en el ejemplo anterior, vamos a resolver
                             este problema usando diferencias finitas para cada uno 
                             de los diferentes $R$ dados en el ejemplo \ref{ejemplo4}. 
                             Para \'{e}sto vamos a aplicar el m\'{e}todo expl\'{i}cito  
                             y el m\'{e}todo impl\'{i}cito.
                       \item Despu\'{e}s vamos a comparar las soluciones que se obtienen
                             con ambos m\'{e}todos.                             
                      \end{enumerate} 
                      
                      Para resolver el problema (\ref{Cuadrado:tiemLi}) vamos a 
                      considerar $\overline {\Omega}$ para cada $R$ y a denotar
                      con $u_{i j}^n$ la soluci\'{o}n aproximada de $a_{i j}^n$ en el
                      punto $(i \Delta x,j \Delta y, n \Delta t)\in \overline {\Omega} \times t_n$,
                      donde $t_n$ es la discretizaci\'{o}n de la variable temporal $t$.
                      
                      A continuaci\'{o}n vamos a discretizar el problema (\ref{Cuadrado:tiemLi}),
                      primero lo haremos con el m\'{e}todo expl\'{i}cito y despu\'{e}s
                      con el m\'{e}todo impl\'{i}cito.
                                                                                                                                    
                      De acuerdo al m\'{e}todo expl\'{i}cito, primero vamos a 
                      aproximar de manera adecuada las derivadas parciales y 
                      despu\'{e}s vamos a sustituir estas aproximaciones en la
                      ecuaci\'{o}n dada en (\ref{Cuadrado:tiemLi}). Esto da 
                      como resultado que para cada  
                      $(i \Delta x,j \Delta y, n \Delta t)\in \overline {\Omega} \times t_n$
                      \begin{equation*}
                            \displaystyle\frac{u_{i j}^{n+1} - u_{i j}^n}{\Delta t} =
                            \displaystyle\frac{u_{i+1 j}^n - 2 u_{i j}^n + u_{i-1 j}^n}{\Delta x^2} + \displaystyle\frac{u_{i j+1}^n - 2 u_{i j}^n + u_{i j-1}^n}{\Delta y^2},                                                                  
                      \end{equation*} 
                      y como estamos interesado en calcular el nuevo tiempo $u_{i j}^{n+1}$
                      a partir de los valores del tiempo anterior, entonces resolvemos 
                      esta igualdad para $u_{i j}^{n+1}$ y obtenemos que para cada  
                      $(i \Delta x,j\Delta y, n\Delta t)\in \overline {\Omega} \times t_n$,
                      \begin{equation*}
                            u_{i j}^{n+1} = u_{i j}^n +
                            \displaystyle\frac{\Delta t}{\Delta x^2} ( u_{i+1 j}^n - 2 u_{i j}^n + u_{i-1 j}^n ) + \displaystyle\frac{\Delta t}{\Delta y^2} ( u_{i j+1}^n - 2 u_{i j}^n + u_{i j-1}^n ),                                                                 
                      \end{equation*}                       
                      mientras que la discretizaci\'{o}n de la frontera queda de la
                      siguiente manera, 
                      \begin{equation*}                                                                    
                            u_{i j}^n  = g(i \Delta x, j \Delta y, n \Delta t), \text{ para cada }                       
                            (i\Delta x, j\Delta y, n\Delta t) \in \partial \overline {\Omega}\times t_n,                                      
                      \end{equation*} 
                      y la discretizaci\'{o}n de las condiciones iniciales son
                      \begin{equation*}                                                                    
                            u_{i j}^0  = f(i \Delta x, j \Delta y, t_0), \text{ para cada }
                            (i\Delta x, j\Delta y) \in \Omega \cup \partial \overline {\Omega}.                         
                      \end{equation*} 
                                                                                                                                                          
                      Las soluciones aproximadas que obtuvimos con $\Delta x = 0.0040$,
                      $\Delta y = 0.0040$ y $\Delta t = 0.000002$ para diferentes tiempos $t=t_i$ 
                      se pueden apreciar a la izquierda de las figuras \ref{sol_tiem_m1},                      
                      \ref{sol_tiem_m2}, \ref{sol_tiem_m3}, \ref{sol_tiem_m4}, \ref{sol_tiem_m5},
                      y \ref{sol_tiem_m6}.
                                            
                      Ahora vamos a discretizar el problema (\ref{Cuadrado:tiemLi})
                      aplicando el m\'{e}todo impl\'{i}cito. Para \'{e}sto, una vez
                      m\'{a}s vamos a aproximar de las derivadas parciales de manera 
                      adecuada y a sustituir estas aproximaciones en la ecuaci\'{o}n
                      dada en (\ref{Cuadrado:tiemLi}). Una vez hecho lo anterior
                      obtenemos que para cada  
                      $(i \Delta x,j \Delta y, n \Delta t)\in \overline {\Omega} \times t_n$,
                      \begin{equation*}
                            \displaystyle\frac{u_{i j}^{n+1} - u_{i j}^n}{\Delta t} =
                            \displaystyle\frac{u_{i+1 j}^{n+1} - 2 u_{i j}^{n+1} + u_{i-1 j}^{n+1}}{\Delta x^2} + \displaystyle\frac{u_{i j+1}^{n+1} - 2 u_{i j}^{n+1} + u_{i j-1}^{n+1}}{\Delta y^2},                                                                  
                      \end{equation*} 
                      dado que estamos interesados en calcular los valores del nuevo 
                      tiempo a partir de los valores del tiempo anterior, entonces                      
                      resulta claro que esta expresi\'{o}n es equivalente a
                      \begin{equation*}
                            u_{i j}^{n+1} -
                            \displaystyle\frac{\Delta t}{\Delta x^2} (u_{i+1 j}^{n+1} - 2 u_{i j}^{n+1} + u_{i-1 j}^{n+1}) - \displaystyle\frac{\Delta t}{\Delta y^2} (u_{i j+1}^{n+1} - 2 u_{i j}^{n+1} + u_{i j-1}^{n+1}) = u_{i j}^n.                                                                                
                      \end{equation*}
                      para cada  
                      $(i \Delta x,j \Delta y, n \Delta t)\in \overline {\Omega} \times t_n$.
                      
                      Si consideramos estas ecuaciones junto con la discretizaci\'{o}n 
                      de las condiciones de frontera,  
                      \begin{equation*}                                                                    
                            u_{i j}^n  = g(i \Delta x, j \Delta y, n \Delta t), \text{ para cada }                       
                            (i\Delta x, j\Delta y, n\Delta t) \in \partial \overline {\Omega}\times t_n,                                       
                      \end{equation*} 
                      obtenemos un sistema de $m$ ecuaciones lineales con $m$ 
                      incognitas de la forma
                      \begin{equation*}
                            Au^{n+1} = u^n, \text{ para } n=0,1,\dots, 
                      \end{equation*}
                      donde $u^0$ esta dado por las condiciones iniciales discretizadas,
                      esto es                         
                      \begin{equation*}                                                                    
                            u_{i j}^0  = f(i \Delta x, j \Delta y, t_0), \text{ para cada }
                            (i\Delta x, j\Delta y) \in \Omega \cup \partial \overline {\Omega}.                        
                      \end{equation*}.
                                                                                                          
                      Los resultados num\'{e}ricos que obtuvimos con $\Delta x = 0.0040$,
                      $\Delta y = 0.0040$ y $\Delta t = 0.0002$ para diferentes tiempos $t=t_i$  
                      se pueden observar a la derecha de las figuras \ref{sol_tiem_m1},                      
                      \ref{sol_tiem_m2}, \ref{sol_tiem_m3}, \ref{sol_tiem_m4}, \ref{sol_tiem_m5}
                      y \ref{sol_tiem_m6}. Es importante mencionar que para resolver 
                      cada sistema de ecuaciones lineales aplicamos el m\'{e}todo 
                      iterativo de Jacobi.                      
                      
                      Al comparar las soluciones num\'{e}ricas para los mismo 
                      tiempos $t=t_i$, podemos observar que \'{e}stas son muy
                      parecidas entre ellas. A pesar de la diferencias que existen
                      entre los dos m\'{e}todos y del error introducido por la 
                      estructura de la discretizaci\'{o}n inducida por cada $R$. 
                      
                      De hecho, despu\'{e}s de realizar varias pruebas observamos
                      que para cada $R$ se obtienen los mismos resultados considerando 
                      diferentes valores adecuados para $\Delta x$, $\Delta y$ 
                      y $\Delta t$. En particular notamos que cuando resolvimos 
                      el problema con el m\'{e}todo expl\'{i}cito los valores
                      propuestos deben de satisfacer la condici\'{o}n que se
                      propone en \cite{Jwth}, mientras que con el m\'{e}todo
                      impl\'{i}cito no fue as\'{i}.                      
                \end{ejemplo}

                \begin{figure}[h!]  
                \centering
                \begin{tabular}{||c  c||}
                \hline                                   
                \includegraphics[width=5.0cm,height=4.0cm]{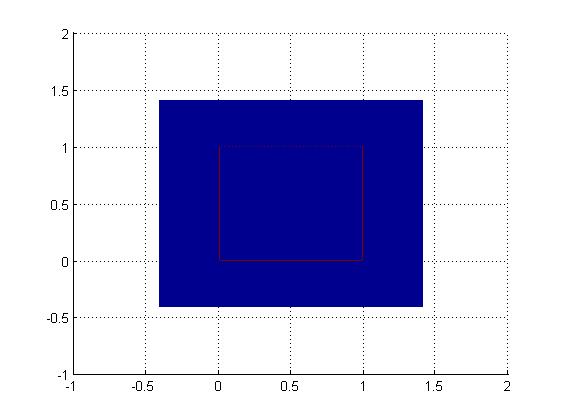} & \includegraphics[width=5.0cm,height=4.0cm]{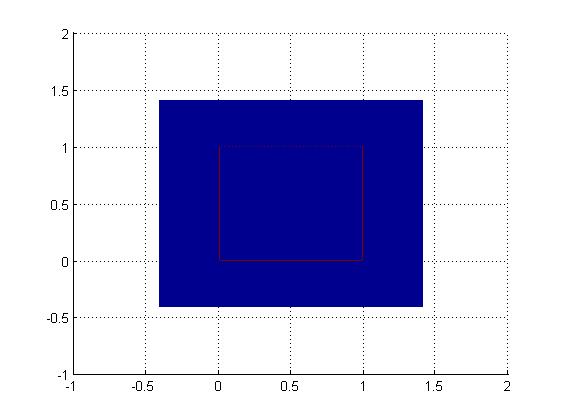} \\ 
                  \multicolumn{2}{ ||c|| } {Condiciones iniciales} \\                   
                \includegraphics[width=5.0cm,height=4.0cm]{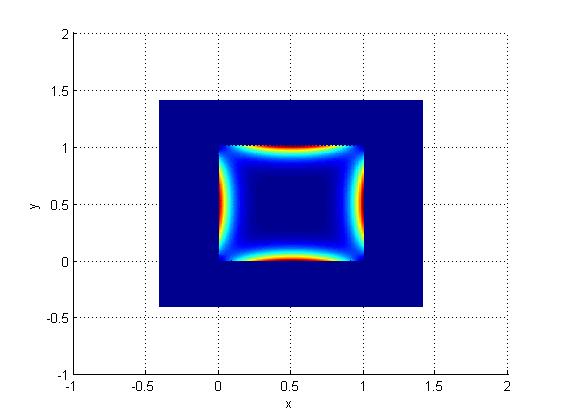} & \includegraphics[width=5.0cm,height=4.0cm]{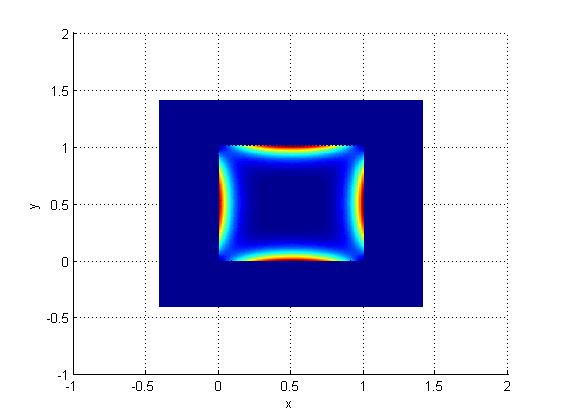} \\ 
                  \multicolumn{2}{ ||c|| } { $ t = 0.0096 $ } \\                                                             
                \includegraphics[width=5.0cm,height=4.0cm]{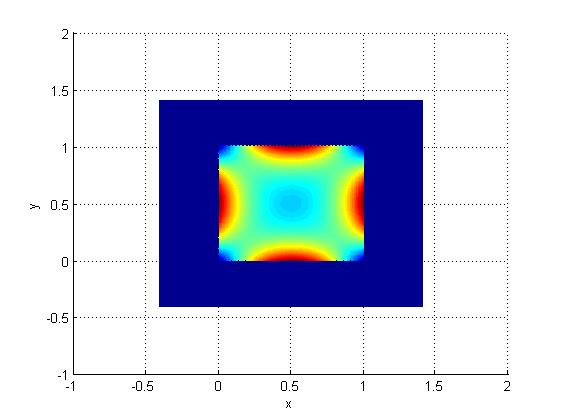} & \includegraphics[width=5.0cm,height=4.0cm]{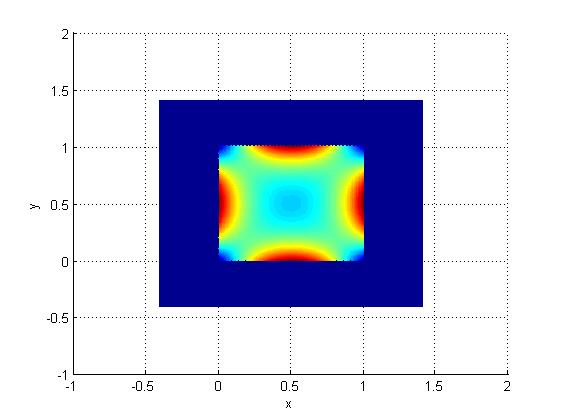} \\ 
                  \multicolumn{2}{ ||c|| } { $ t = 0.1052 $ } \\                    
                \includegraphics[width=5.0cm,height=4.0cm]{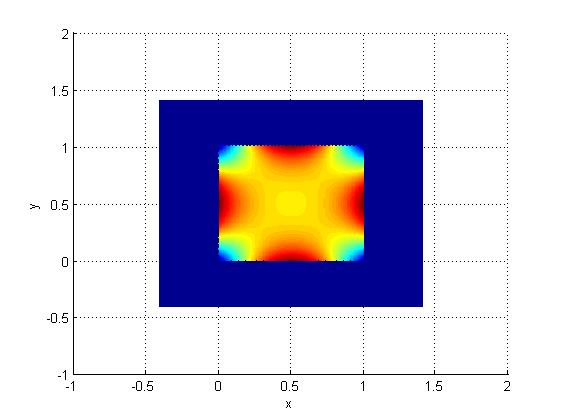} & \includegraphics[width=5.0cm,height=4.0cm]{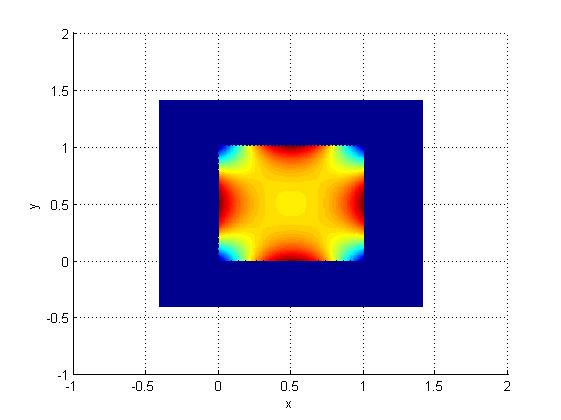} \\ 
                  \multicolumn{2}{ ||c|| } { $ t = 0.3548 $ } \\                    
                \hline
                \end{tabular}
                \caption{Soluci\'{o}n num\'{e}rica de la ecuaci\'{o}n (\ref{Cuadrado:tiemLi})
                         sobre $\Omega$, con el m\'{e}todo expl\'{i}cito e impl\'{i}cito, de 
                         izquierda a derecha, respectivamente.}
                \label{sol_tiem_m6}
                \end{figure}
                
                Ahora bien, si nuestro prop\'{o}sito es resolver un problema 
                sobre un dominio irregular aplicando diferencias finitas, ver  
                la secci\'{o}n \ref{Metodo}, lo ideal es no elegir tan 
                deliberadamente $R$, ya que una elecci\'{o}n adecuada de $R$ 
                nos puede ayudar a disminuir \'{o} incluso anular los errores 
                que se introducen cuando consideramos el dominio extendido 
                de $\Omega$.                                                                                
                
                Esto nos lleva a que nos preguntemos \textquestiondown c\'{o}mo 
                debemos de construir $R$ tal que la discretizaci\'{o}n,
                $\overline {\Omega}_{m,n}$ y $\underline{\Omega}_{m,n}$, inducida 
                por $G_R$ sea \'{o}ptima para $\Omega$? o equivalentemente                
                
                \begin{center}
                      {\it 
                           
                       \textbf{*** \textquestiondown cu\'{a}l es la forma m\'{a}s 
                               \'{o}ptima de ocupar un dominio $\Omega$ con el 
                               stencil?***}                                      
                                }
                \end{center}

                {\it
                Decimos que una discretizaci\'{o}n $\overline {\Omega}_{m,n}$ 
                \'{o} $\underline{\Omega}_{m,n}$ inducida por $G_R$ es una 
                buena discretizaci\'{o}n del dominio $\Omega$ si
                \begin{itemize}
                      \item el \'{a}rea inducida por $\overline {\Omega}_{m,n}$
                            \'{o} $\underline{\Omega}_{m,n}$ coincide con el
                            \'{a}rea usual de $\Omega$. Pero, \textquestiondown
                            a qu\'{e} nos estamos refiriendo con esta afirmaci\'{o}n?
                            
                            Antes de contestar a esta pregunta, es importante 
                            mencionar que en t\'{e}rminos generales vamos a  
                            tratar el \'{a}rea inducida de la misma manera que 
                            el \'{a}rea usual, ya que estamos interesados en 
                            considerar esta \'{a}rea como usualmente comprendemos.
                            De hecho, hasta este momento no hay razones para 
                            considerar de forma distinta el \'{a}rea inducida
                            del \'{a}rea usual.
                            
                            Ahora pasemos a responder a nuestra pregunta, para
                            esto notemos que por construcci\'{o}n el \'{a}rea 
                            inducida y el \'{a}rea usual estan superpuestas, 
                            por lo que en nuestra afirmaci\'{o}n nos estamos 
                            refiriendo a que si la forma geom\'{e}trica de las
                            \'{a}reas de los dominios son iguales puntualmente,
                            sin considerar rotaciones, simetr\'{i}as, traslaciones,
                            etc. As\'{i} mismo debe de quedar claro que no nos
                            estamos refiriendo a que si el c\'{a}lculo de las
                            \'{a}reas de \'{e}stos dominios son iguales.                                                                                                                                                                                                                                                                                               
                            
                            Por ejemplo, el \'{a}rea inducida 
                            por $\overline {\Omega}_{m,n}$ de la malla $A$ de la 
                            figura \ref{mallas:fig1} es un ejemplo de \'{e}sto.
                      \item todos los nodos de $\overline {\Omega}_{m,n}$ 
                            \'{o} $\underline{\Omega}_{m,n}$ involucrados en el 
                            c\'{a}lculo de la soluci\'{o}n pertenecen a $\Omega$.
                            Un ejemplo de \'{e}sto se puede observar en la malla
                            $C$ de la figura \ref{mallas:fig1}.
                \end{itemize}}

                En el \'{u}ltimo caso hay que recordar que hay nodos que no aportan
                informaci\'{o}n al calculo num\'{e}rico  de la soluci\'{o}n y es 
                por esta raz\'{o}n que, desde el punto de vista pr\'{a}ctico,
                consideramos que esta \'{u}ltima forma de asociar una discretizaci\'{o}n
                a un dominio es correcta.

                \begin{center}
                      {\it En t\'{e}rminos m\'{a}s generales estamos buscando 
                           construir $R$ tal que                            
                           
                       \textbf{para toda $n$ y $m$ la discretizaci\'{o}n, 
                               $\overline {\Omega}_{m,n}$ \'{o} $\underline{\Omega}_{m,n}$,
                               inducida por $G_R$ sea \'{o}ptima para $\Omega$.}
                               }                                                            
                \end{center}                                                                                                                                  
                            
                \begin{ejemplo} \label{ejemplo7} 
                      Consideremos nuevamente el caso cuando el dominio $\Omega$ 
                      es un cuadrado, es decir $\Omega=[a,b]\times[a,b]$ y supongamos 
                      de nueva cuenta que queremos aplicar diferencias finitas de 
                      acuerdo a la secci\'{o}n \ref{Metodo} para resolver un problema 
                      sobre $\Omega$. 
                      
                      Para este dominio, resulta claro que la \'{o}ptima 
                      discretizaci\'{o}n $\overline {\Omega}_{m,n}$ inducida 
                      por $G_R$ para toda $n$ y $m$, es aquella en la que 
                      todos los lados de $R$ coinciden con todos los lados
                      del cuadrado $\Omega$. Adem\'{a}s tenemos que para 
                      toda $n$ y $m$ se satisfacen las siguientes afirmaciones.
                      
                      \begin{itemize}                      
                            \item La discretizaci\'{o}n usual de $\Omega$ coincide con 
                                  la discretizaci\'{o}n dada por $\overline {\Omega}_{m,n}$.
                            
                            \item El \'{a}rea usual de $\Omega$ coincide con el
                                  \'{a}rea inducida por $\overline {\Omega}_{m,n}$.                                                                                                                                                         
                      \end{itemize}
                      
                      Por otra parte si $R$ es el cuadrado m\'{a}s peque\~{n}o girado 
                      $45^\circ$ tal que $\Omega \subset R$, entonces
                      
                      \begin{itemize}                      
                            \item Existe una discretizaci\'{o}n de $R$ tal que 
                                  $\underline{\Omega}_{m,n}$ es una buena 
                                  discretizaci\'{o}n. 
                                  
                                  Para ver esto basta con usar como referencia 
                                  las mitad de las dos diagonales de $R$ para 
                                  definir los valores de la discretizaci\'{o}n  
                                  de los intervalos. Y prolongar esta misma
                                  discretizaci\'{o}n a lo largo de las dos 
                                  diagonales.
                      \end{itemize}
                                      
                \end{ejemplo}
                                
                En este punto es importante notar que no siempre sucede que 
                cuando el stencil ocupa todo el dominio $\Omega$ se tiene 
                que el \'{a}rea inducida coincide con el \'{a}rea usual del
                dominio. Por lo que a continuaci\'{o}n vamos a {\it explicar 
                y a definir como medir el error que ocurre.}                                                                 
                
                {\it Dada una discretizaci\'{o}n $\overline {\Omega}_{m,n}$ asociada
                al dominio $\Omega$, definimos el error $\overline {E}(m,n)$  
                como                              
                                
                \begin{equation}
                      \overline {E}(m,n) = \left \{ \begin{matrix} 0 & \mbox{si las \'{a}reas $\overline {\Omega}_{m,n}$ y $\Omega$ coinciden,} \\
                                            | \overline{\Omega}_{m,n} - \Omega | & \mbox{si las \'{a}reas  $\overline {\Omega}_{m,n}$ y $\Omega$ no coinciden,}\end{matrix} \right.                              
                \end{equation}
                y de manera an\'{a}loga se define $\underline{E}(m,n)$.} 
                
                De manera inmediata notamos que el resultado de $\overline {E}(m,n)$
                y de $\underline{E}(m,n)$ s\'{o}lo puede ser cero \'{o} ser
                un dominio distinto del vac\'{i}o. De hecho este nuevo dominio
                no necesariamente tendr\'{a} las mismas propiedades que los
                dominios $\overline {\Omega}_{m,n}$ y $\Omega$ \'{o} que los
                dominios $\underline {\Omega}_{m,n}$ y $\Omega$, respectivamente.

                As\'{i} que debido a la 
                naturaleza de $\overline {E}(m,n)$ y $\underline{E}(m,n)$ 
                necesitamos establecer una {\it medida que nos permita hablar 
                acerca de su magnitud.} Por lo que de manera natural resulta
                evidente que el {\it c\'{a}lculo del \'{a}rea} puede usarse para \'{e}ste fin.
                Esto es,
                
                \begin{equation}
                       \overline {A}(m,n) = \left \{ \begin{matrix} 0 & \mbox{ si } \overline {E}(m,n) =0, \\
                                            Area ( | \overline{\Omega}_{m,n} - \Omega |) & \mbox{en otro caso,}\end{matrix} \right.                              
                \end{equation}
                y de manera an\'{a}loga se define $\underline{A}(m,n)$.
                
                Ahora bien, hay que recordar que por la forma en que construimos 
                $\overline {\Omega}$ y $\underline{\Omega}$ tenemos que    
                cuando $\Delta x \rightarrow 0 $ y $\Delta y \rightarrow 0$, 
                entonces $\overline {\Omega}_{m,n} \rightarrow \Omega$ 
                y $\underline{\Omega}_{m,n} \rightarrow \Omega$, por lo que 
                los siguientes l\'{i}mites est\'{a}n bien definidos
                \begin{equation}
                 \label{limiteError}
                       lim_{m \rightarrow \infty , \hspace{0.2cm} n \rightarrow \infty } \hspace{0.2cm} 
                                 \overline{E}(m,n) =0, \hspace{1cm} 
                       lim_{m \rightarrow \infty , \hspace{0.2cm} n \rightarrow \infty } \hspace{0.2cm} 
                                 \underline{E}(m,n) =0.
                \end{equation}
                este hecho implica que los siguientes l\'{i}mites tambi\'{e}n
                est\'{e}n bien definidos 
                \begin{equation}
                 \label{limiteAError}
                       lim_{m \rightarrow \infty , \hspace{0.2cm} n \rightarrow \infty } \hspace{0.2cm} 
                                 \overline{A}(m,n) =0, \hspace{0.8cm}
                       lim_{m \rightarrow \infty , \hspace{0.2cm} n \rightarrow \infty } \hspace{0.2cm} 
                                 \overline{A}(m,n) =0.
                \end{equation}   
                
                Por lo tanto las sucesiones generadas por $\overline{E}(m,n)$, 
                $\underline{E}(m,n)$, $\overline{A}(m,n)$ y $\underline{A}(m,n)$ 
                est\'{a}n bien definidas y convergen a cero.                                                
            
                \begin{itemize}                                           
                                                                                                                
                      \item Notemos que el \'{a}rea $( | \underline{\Omega}_{m,n} - \Omega |)$  
                            es diferente de $| \text{\it{\'{a}}}rea (\underline{\Omega}_{m,n})-  \text{\it{\'{a}}}rea(\Omega) |$, ya que $| \text{\it{\'{a}}}rea(\underline{\Omega}_{m,n})-  \text{\it{\'{a}}}rea(\Omega) |$ 
                            puede ser igual a cero y no necesariamente los dominios 
                            $\underline{\Omega}_{m,n}$ y $\Omega$ coinciden.                                                  
                      
                \end{itemize}

                Antes de finalizar esta secci\'{o}n es importante mencionar que
                dado que la definici\'{o}n de $\underline{\Omega}$ y $\overline {\Omega}$ 
                no depende de la elecci\'{o}n de $R$ ni de la discretizaci\'{o}n 
                inducida por \'{e}ste, entonces $R$ puede ser construido y depender
                de $\Delta x$ y $\Delta y$. En otras palabras el tama\~{n}o de $R$
                tambi\'{e}n puede variar cuando $\Delta x$ y $\Delta y$ varian, 
                pero con la condici\'{o}n de que cualquier extensi\'{o}n 
                de $\Omega$ siempre sea un subconjunto completamente contenido 
                en $R$. Algunos ejemplos de \'{e}sto se muestran a continuaci\'{o}n.
                
                \begin{ejemplo} \label{ejemplo8}
                      Consideremos que se tienen los siguientes dominios
                      \begin{enumerate}
                       \item Si el dominio es un cuadrado, es decir $\Omega=[a,b]\times[a,b]$, 
                             entonces $R=[a-\Delta x,b+\Delta x]\times[a-\Delta y,b+\Delta y]$
                       \item Si el dominio es un rect\'{a}ngulo, es decir $\Omega=[a,b]\times[c,d]$, 
                             entonces $R=[a-\Delta x,b+\Delta x]\times[c-\Delta y,d+\Delta y]$  
                       \item Si el dominio $\Omega$ es irregular y $min(x)$, $min(y)$,
                             $max(x)$ y $max(y)$ son los valores m\'{i}ninos y m\'{a}ximos,
                             respectivamente, del dominio $\Omega$, entonces
                             $R=[min(x)-\Delta x,max(x)+\Delta x]\times[min(y)-\Delta y,max(y)+\Delta y]$.
                      \end{enumerate}                      
                \end{ejemplo}
                
                De hecho, para llevar a cabo los c\'{a}lculos num\'{e}ricos, 
                se recomienda en algunos casos que $R$ dependa de $\Delta x$
                y $\Delta y$, ya que \'{e}sto nos ayudar\'{a} a disminuir el 
                tiempo de c\'{o}mputo.

\section{Propiedades}

                Resulta claro que la forma en que estamos aplicando 
                diferencias finitas, de acuerdo al m\'{e}todo descrito
                en la secci\'{o}n \ref{Metodo}, para resolver problemas 
                que involucran ecuaciones diferenciales parciales 
                sobre dominios planos irregurales difiere de cualquier 
                m\'{e}todo de diferencias finitas que usualmente 
                conocemos para resolver este tipo de problemas.
                
                De hecho los m\'{e}todos de diferencias finitas que se
                conocen para resolver este tipo de problemas funcionan 
                \'{u}nicamente cuando el dominio del problema cumple 
                ciertas caracter\'{i}sticas muy espec\'{i}ficas.
                
                Sin embargo, como el m\'{e}todo de diferencias finitas 
                de la secci\'{o}n \ref{Metodo} est\'{a} basado en los
                met\'{o}dos de diferencias finitas que usualmente 
                conocemos, entonces
                \begin{enumerate}
                      \item los m\'{e}todos de diferencias finitas
                            resultan ser un caso particular del 
                            m\'{e}todo de diferencias finitas de 
                            la secci\'{o}n \ref{Metodo}.
                
                      \item creemos que todas las propiedades que 
                            conocemos de los met\'{o}dos de diferencias 
                            finitas tambi\'{e}n las satisface, bajo ciertas
                            condiciones, el m\'{e}todo de diferencias
                            finitas de la secci\'{o}n \ref{Metodo}.
                \end{enumerate}

                \begin{ejemplo} \label{ejemplo9} Consideremos el siguiente problema
                      con condiciones de Dirichlet en la frontera y con condiciones
                      iniciales, ver \cite{Jwth},
                      \begin{eqnarray}
                       \label{Cuadrado:par}
                            a_t & = & \nu \nabla^2 a,\hspace*{0.3cm}\text{para}\hspace*{0.3cm}(x,y)\in\Omega, \hspace*{0.3cm} t>0, \nonumber \\ 
                            a(x,y,t) & = & 0,\hspace*{0.3cm} \text{para} \hspace*{0.3cm} (x,y) \in   \partial  \Omega , \hspace*{0.3cm} t>0,\\ 
                            a(x,y,0) & = & sin(\pi x)sin(2\pi y),\hspace*{0.3cm} \text{para} \hspace*{0.3cm} (x,y) \in  \Omega \cup \partial  \Omega, \nonumber                                                                 
                      \end{eqnarray}                      
                      donde $a=a(x,y,t)$,
                      $\Omega$ es el cuadrado unitario $(0,1)\times(0,1)$ y $\partial  \Omega$ 
                      es la frontera de $\Omega$.

                      \begin{enumerate}                       
                       \item Ahora vamos a resolver este problema usando diferencias 
                             finitas para cada uno de los diferentes $R$ dados en el
                             ejemplo \ref{ejemplo4}. Para \'{e}sto vamos a aplicar 
                             el m\'{e}todo expl\'{i}cito. 
                             
                       \item Despu\'{e}s vamos a calcular el que se obtiene entre las
                             soluciones aproximadas y la soluci\'{o}n conocida.       
                      \end{enumerate}

                      Vamos a resolver este problema usando diferencias finitas 
                                                                  
                      Los resultados num\'{e}ricos que obtuvimos para diferentes 
                      tiempos $t=t_i$ con el m\'{e}todo expl\'{i}cito sobre las mallas
                      del ejemplo \ref{ejemplo4} se pueden apreciar en la figuras. 
                      Es importante mencionar que no obtuvimos una soluci\'{o}n 
                      num\'{e}rica usando la malla $C$.

                \end{ejemplo}

                \begin{figure}[h!]  
                \centering
                \begin{tabular}{||c  c||}
                \hline                                   
                \includegraphics[width=5.0cm,height=4.0cm]{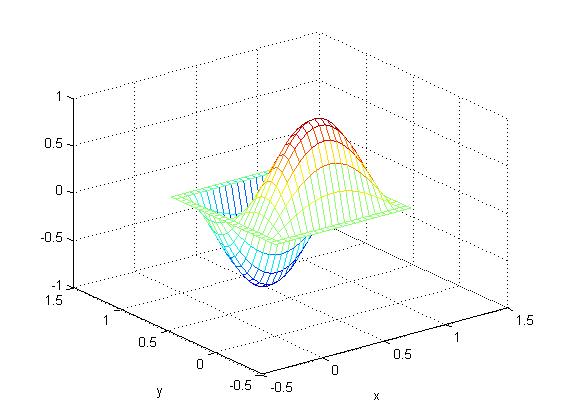} & \includegraphics[width=5.0cm,height=4.0cm]{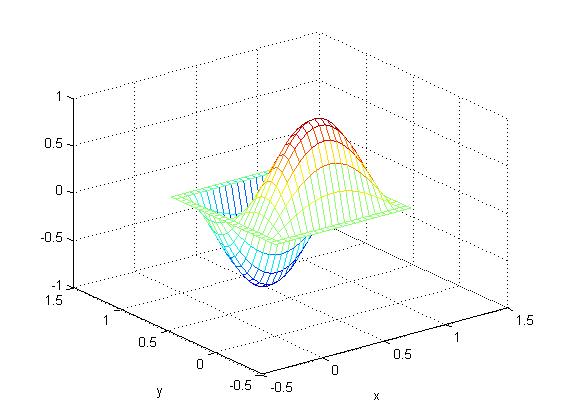} \\ 
                  \multicolumn{2}{ ||c|| } {Condiciones iniciales} \\                   
                \includegraphics[width=5.0cm,height=4.0cm]{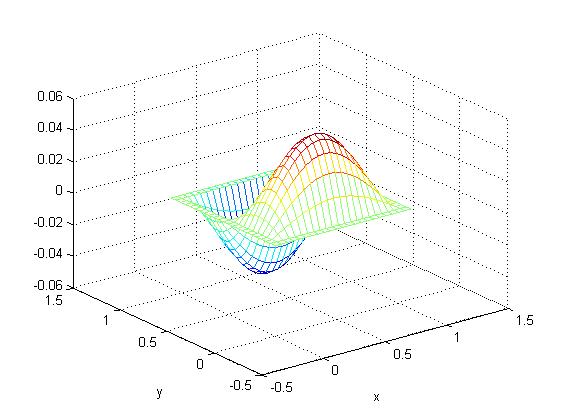} & \includegraphics[width=5.0cm,height=4.0cm]{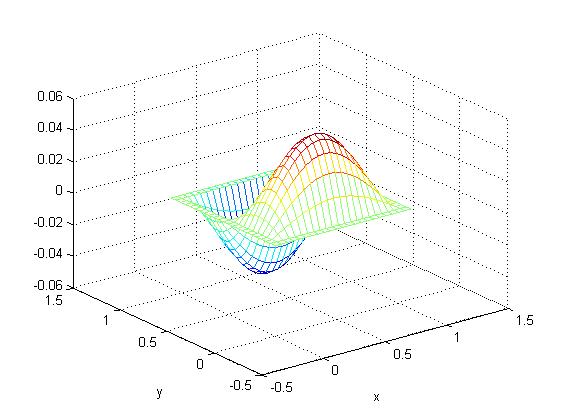} \\ 
                  \multicolumn{2}{ ||c|| } { $ t = 0.06 $ } \\                                                             
                \includegraphics[width=5.0cm,height=4.0cm]{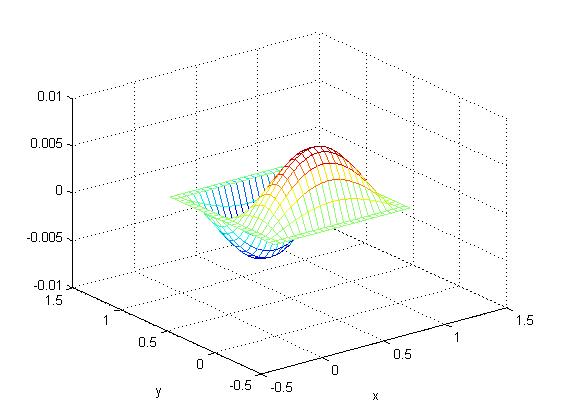} & \includegraphics[width=5.0cm,height=4.0cm]{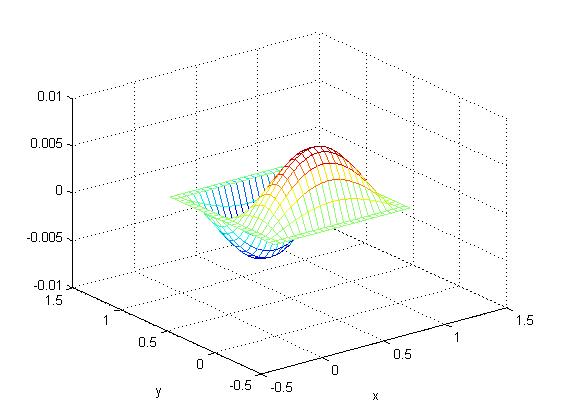} \\ 
                  \multicolumn{2}{ ||c|| } { $ t = 0.1 $ } \\                    
                \includegraphics[width=5.0cm,height=4.0cm]{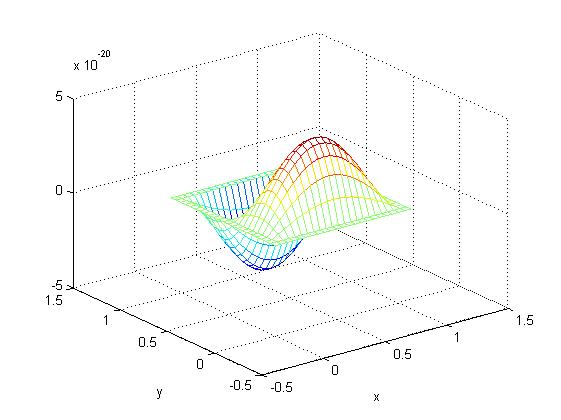} & \includegraphics[width=5.0cm,height=4.0cm]{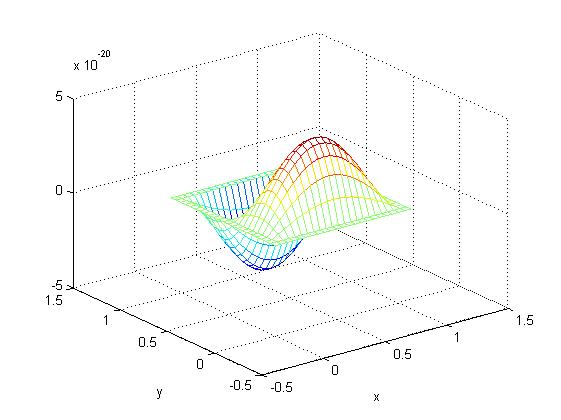} \\ 
                  \multicolumn{2}{ ||c|| } { $ t = 0.9 $ } \\                    
                \hline
                \end{tabular}
                \caption{Soluci\'{o}n num\'{e}rica de la ecuaci\'{o}n (\ref{Cuadrado:par})
                         sobre $\Omega$, con el m\'{e}todo expl\'{i}cito e impl\'{i}cito, de 
                         izquierda a derecha, respectivamente.}
                \label{sol_tiem_m6}
                \end{figure}

                      \begin{figure}[t]
                      \centering
                      \begin{tabular}{||c||}
                      \hline                
                      \includegraphics[width=6cm,height=5cm]{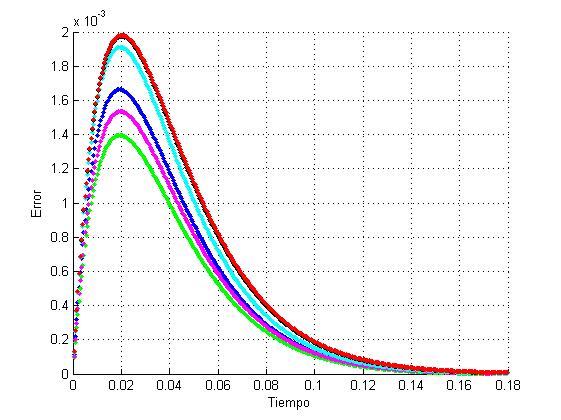} \\             
                      \hline
                      \end{tabular}
                      \caption{Grafica de los Errores.}
                      \label{fig:puntosInterioresA}
                      \end{figure}

\section{Convergencia, Consistencia y Estabilidad}                 
                
                Hasta este momento s\'{o}lo hemos visto como construir una 
                sucesi\'{o}n de dominios y de soluciones que aproximan
                a la soluci\'{o}n de un problema dado.
                
                De hecho, esta sucesi\'{o}n de soluciones la obtenemos a 
                trav\'{e}s de aplicar de manera adecuada el m\'{e}todo de 
                diferencias finitas a una sucesi\'{o}n de mallas que aproximan
                a una discretizaci\'{o}n del dominio irregular. 
                
                Esto es, en primer lugar estamos {\it inscribiendo }
                una sucesi\'{o}n de mallas $\underline{\Omega}$ crecientes
                que converge a $\Omega$ cuando $\Delta x \rightarrow 0 $ 
                y $\Delta y \rightarrow 0$. En segundo lugar, estamos 
                {\it circunscribiendo} una sucesi\'{o}n de mallas 
                $\overline {\Omega}$ crecientes que converge a $\Omega$ 
                cuando $\Delta x \rightarrow 0 $ y $\Delta y \rightarrow 0$.

                Por otra parte, debe de quedar claro que no estamos intentando
                definir lo que se entiende por {\it cuadrar una regi\'{o}n},
                lo \'{u}nico que estamos haciendo es construir sucesiones de
                mallas tales que \'{e}stas convergan a la regi\'{o}n $\Omega$.
                
                Informacion adicional sobre la naturaleza del problema
                puede usarse para desarrollar algoritmos que produzcan
                sucesiones que converjan a la solucion.

                Entonces tiene sentido preguntarnos de qu\'{e} tipo de 
                convergencia estamos hablando. Por lo que a continuaci\'{o}n
                pasamos a formalizar \'{e}ste y otros conceptos.
                
                M\'{a}s a\'{u}n, cuando se lleve a cabo el an\'{a}lisis del error entre 
                la soluci\'{o}n exacta y la soluci\'{o}n aproximada hay que tener
                en cuenta un nuevo t\'{e}rmino de error: el de las mismas diferencias 
                finitas introducido por la forma en que se discretiz\'{o} la 
                regi\'{o}n.

                \begin{ejemplo} Consideremos el siguiente problema tomado
                    del libro Kelley                
                      \begin{eqnarray}
                             - \nabla ( cos(x) \nabla a ) & = & f(x,y),\hspace*{0.3cm}\text{para}\hspace*{0.3cm}(x,y)\in\Omega, \\
                                     a & = & 0,\hspace*{0.3cm} \text{para} \hspace*{0.3cm} (x,y) \in   \partial  \Omega , \nonumber
                      \end{eqnarray}    
                      donde $a=a(x,y)$,$f(x,y)$ es tal que la soluci\'{o}n es 
                      $10xy(1-x)(1-y)exp(x^{4.5})$, $\Omega$ es el dominio 
                      determinado por la par\'{a}bola $y=-x(x-1)$
                      y el eje horizontal. Vamos a resolver este problema mediante
                      el m\'{e}todo de la secci\'{o}n \ref{Metodo}.

                \end{ejemplo}

                Para cada $\Delta x $ y $\Delta y $ obtenemos una discretizaci\'{o}n 
                de un rect\'{a}ngulo que contiene a la regi\'{o}n $\Omega$ y 
                construimos
                \begin{eqnarray}
                    \overline {\Omega}_1 ;\hspace{0.2cm} \overline {\Omega}_2; \hspace{0.2cm} \overline {\Omega}_3 \hspace{0.2cm} \cdots \\
                    \underline{\Omega}_1 ; \hspace{0.2cm} \underline{\Omega}_2; \hspace{0.2cm} \underline{\Omega}_3 \hspace{0.2cm} \cdots                 
                \end{eqnarray}
                de tal manera que cuando $\Delta x \rightarrow 0 $ y $\Delta y \rightarrow 0 $
                las sucesiones convergan a $\Omega$, esto es
                \begin{equation}
                      lim_{\Delta x \rightarrow 0, \hspace{0.2cm} \Delta y \rightarrow 0 } \hspace{0.2cm} \overline {\Omega} = \Omega \hspace{0.5cm} \text{ y } \hspace{0.5cm}  lim_{ \Delta x \rightarrow 0, \hspace{0.2cm} \Delta y \rightarrow 0 } \hspace{0.2cm}   \underline{\Omega} = \Omega,    
                \end{equation}                                                                               
                as\'{i} que, por cada $\underline{\Omega}$ y $\overline {\Omega}$
                definimos,                                                                
                \begin{eqnarray*}
                       \underline {\omega} & = & \{ \text{es el conjunto de todos los rect\'{a}ngulos}   \\
                                           &   &  \text{ que est\'{a}n en el interior de } {\Omega} \}, \\                                           
                       \overline {\omega}  & = &  \{ \text{es el conjunto de todos los rect\'{a}ngulos }   \\
                                           &   & \text{ que est\'{a}n en el interior y en la frontera de }  {\Omega}  \}.
                \end{eqnarray*}  
                
                Hasta aqu\'{i}, debe de estar claro que cada discretizaci\'{o}n $\underline{\Omega}$ 
                y $\overline {\Omega}$ tiene asociada su propia regi\'{o}n $\underline {\omega}$ 
                y $\overline {\omega}$, respectivemente. Ya que el conjunto $\underline {\omega}$
                y  $\overline {\omega}$ es un conjunto de rect\'{a}ngulos que poseen ciertas
                caracter\'{i}sticas, mientras que el conjunto $\underline {\Omega}$ y
                $\overline {\Omega}$ est\'{a} conformado por los v\'{e}rtices de estos mismos 
                rect\'{a}ngulos, respectivamente. En particular,
                \begin{equation}
                 \underline{\Omega} \text{ } \cap \text{ } \underline {\omega} = \underline{\Omega}  
                \text{ y } \overline {\Omega} \text{ } \cap \text{ } \overline {\omega} = \overline {\Omega}.
                \end{equation}                
                
                Dado que la definici\'{o}n de $\underline{\Omega}$ y $\overline {\Omega}$ no  
                depende de la elecci\'{o}n de $R$ ni de la discretizaci\'{o}n de \'{e}ste. 
                Entonces definimos,                
                \begin{eqnarray*}
                       \underline{\Omega}^* & = & sup \{\underline {\Omega} \cap G_R | G_R  \text{
                                  es cualquier discretizaci\'{o}n de }  R  \}, \\
                       \overline{\Omega}^* & = & inf \{\overline{\Omega} \cap G_R   |  G_R \text{ 
                                 es cualquier discretizaci\'{o}n de } R  \},
                \end{eqnarray*}   
                donde $sup$ e $inf$ es el {\it supremo} y el {\it infimo}, respectivamente.
                
                Esta \'{u}ltima definici\'{o}n coincide con la definici\'{o}n del {\it 
                contenido de Jordan}, lo cual tiene sentido, ya que nos estamos basando
                en esas propiedades. Adem\'{a}s, sabemos que $\Omega$ es el contenido
                de Jordan si
                \begin{equation}
                      \underline{\Omega}^*= \Omega =\overline{\Omega}^* .
                \end{equation}

                \'{E}sto nos permite afirmar que ambas discretizaciones
                se aproximan de manera adecuada a $\Omega$.
                
                De manera an\'{a}loga definimos 
                \begin{eqnarray*}
                       \underline {\omega}^* & = & sup  \{ \text{la regi\'{o}n limitada por 
                                 la frontera de } \underline {\Omega}  \}, \\
                       \overline {\omega}^* & = & inf  \{ \text{la regi\'{o}n limitada por 
                                 la frontera de } \overline {\Omega}  \}.
                \end{eqnarray*}   
                
                En esta definici\'{o}n nos referimos a el {\it supremo} y al {\it infimo}
                pero con respecto al \'{a}rea de $ \underline {\Omega} $ y de
                $ \overline {\Omega} $, respectivamente. De aqu\'{i} se sigue que
                \begin{equation}
                      \underline{\omega}^*= \Omega =\overline{\omega}^*.
                \end{equation}

                Ahora que ya contamos dos sucesiones $\underline {\omega}$ y
                $\overline {\omega}$ que aproximan a la regi\'{o}n $\Omega$,
                ya estamos en condiciones de explicar como nos vamos a aproximar
                a la soluci\'{o}n del problema. Para esto, consideremos                 
                una ecuaci\'{o}n diferencial parcial 
                \begin{equation}
                       \mathcal{L} v= F,
                \end{equation}
                donde $v$ y $F$ son funciones y $\mathcal{L} $ es el operador 
                diferencial.                
                
                Denotemos con
                \begin{equation}
                      I(\Omega, \mathcal{L} v=F ) , \text{ } I(\underline{\omega}, 
                \mathcal{L} v=F )  \text{ y }  I (\overline { \omega },  \mathcal{L} v=F )
                \end{equation}                            
                la soluci\'{o}n exacta de la ecuaci\'{o}n diferencial parcial 
                sobre las regiones $\omega$, $\underline{\omega}$ y $\overline{\omega}$, 
                respectivamente.
                
                Decimos que $I(\underline{\omega},\mathcal{L} v=F )$ y 
                $ I (\overline { \omega }$, $\mathcal{L} v=F )$ son convergentes 
                si
                \begin{eqnarray}                
                       \displaystyle\lim_{ \Delta x, \Delta y \rightarrow 0}{I(\underline{ \Omega}, \mathcal{L}v=F)} & = &  I(\Omega, \mathcal{L}v=F ) , \\
                       \displaystyle\lim_{ \Delta x, \Delta y \rightarrow 0}{I(\overline { \Omega},  \mathcal{L} v=F )} & = & I(\Omega, \mathcal{L} v=F ), \nonumber
                \end{eqnarray}
                
                Decimos que $I(\underline{\omega},\mathcal{L} v=F )$ y 
                $  I (\overline { \omega },  \mathcal{L} v=F )$ aproximan
                a $ I(\Omega, \mathcal{L} v=F )$                
                \begin{eqnarray}                
                       \displaystyle\lim_{ \Delta x, \Delta y \rightarrow 0}{I(\underline{ \Omega}, \mathcal{L}v=F)} & = &  I(\Omega, \mathcal{L}v=F ) , \\
                       \displaystyle\lim_{ \Delta x, \Delta y \rightarrow 0}{I(\overline { \Omega},  \mathcal{L} v=F )} & = & I(\Omega, \mathcal{L} v=F ), \nonumber
                \end{eqnarray}
                
                Como $\underline{\omega}$ y $\overline{\omega}$ son dos sucesiones,
                generadas junto con $\underline{\Omega}$ y $\overline{\Omega}$,
                que aproximan a $\Omega$.  
                
                donde $\textbardbl . \textbardbl$ es una norma.
                $\mathcal{L} v=F$ es ``integrable o soluble''  si se satisface

                En este sentido, si ambas sucesiones convergen al poblema original y la
                forma en que se da esta convergencia se mantiene {\it bien portada}, 
                entonces decimos que el problema $\mathcal{L} v=F$ es de dominio
                estable si peque\~{n}os cambios en la regi\'{o}n $\Omega$ producen
                peque\~{n}os cambios en la soluci\'{o}n del problema, es decir 
                los errores no crecen exponencialmente.
                
                Es f\'{a}cil ver que una de las condiciones para que $\mathcal{L} v=F$  
                sea de dominio estable es que las condiciones de frontera deben de 
                ser al menos diferenciables y continuas en un entorno de $\Omega$,
                tal que contiene a $\underline{ \Omega}$ y $\overline { \Omega}$. 
                
                Desde el punto de vista num\'{e}rico, este hecho implica,
                que el error num\'{e}rico de aproximar\footnote{Hay que tener presente que
                el c\'{a}lculo de $\Omega$ se hace usando aritm\'{e}tica de punto flotante.}
                $\Omega$ no estar\'{i}a afectando de manera muy grave el c\'{a}lculo 
                num\'{e}rico de $\mathcal{L} v=F$.                                 
                
                Ahora supongamos que tenemos un problema  $\mathcal{L} v=F$ que es ``integrable o soluble'', 
                entonces $ I(\underline{\Omega}, \mathcal{L} v=F ) $ 
                y $I (\overline {\Omega},  \mathcal{L} v=F )$ pueden ser aproximados aplicando 
                el m\'{e}todo de diferencias finitas.
                                
                Esto es f\'{a}cil de ver, ya que como estamos cubriendo la regi\'{o}n $\Omega$ 
                con rect\'{a}ngulos y por la construcci\'{o}n de \'{e}stos y de $ \underline{\Omega}$ 
                y $\overline {\Omega}$ se sigue que cada punto interior de $ \underline{\Omega}$ 
                y $\overline {\Omega}$ tiene cuatro puntos equidistantes por lo que de manera 
                natural el m\'{e}todo de diferencias finitas es aplicable.
                De hecho                       
                \begin{eqnarray}                            
                       \displaystyle\lim_{ \Delta x, \Delta y \rightarrow 0 } {\underline{\Omega}}=\Omega, \hspace{0.5cm} \text{y} \hspace{0.5cm}
                       \displaystyle\lim_{ \Delta x, \Delta y \rightarrow 0 } {\overline {\Omega}}=\Omega,                                      
                 \end{eqnarray}                                                                   
                
                Denotemos con $ I(\Omega, L^n_k u^n_k = G^n_k ) $ ,
                $ I(\underline{\Omega}, L^n_k u^n_k = G^n_k ) $ 
                y $ I(\overline {\Omega}, L^n_k u^n_k = G^n_k )$ las soluciones 
                aproximadas de $ I(\Omega, \mathcal{L} v=F ) $,
                $ I(\underline{\Omega}, \mathcal{L} v=F ) $ y
                $ I (\overline {\Omega},  \mathcal{L} v=F )$, respectivamente. 
                
                \begin{enumerate}
                        \item \textbf{Convergencia} Decimos que un esquema $ I(\underline{\Omega}, L^n_k u^n_k = G^n_k ) $
                                que aproxima a la ecuaci\'{o}n diferencial $ I(\underline{\Omega}, \mathcal{L} v=F ) $ 
                                que a su vez aproxima a $ I(\Omega, \mathcal{L} v=F ) $  
                                es un esquema puntualmente convergente si para cualquier $x$, $y$ y $t$, 
                                $ (j \Delta x, k \Delta y, n \Delta t) \in \underline{\Omega} \times \mathbb{R_+}$ converge
                                a $(x,y,t) \in \Omega \times \mathbb{R} $, $u^n_k \in \underline{\Omega} \times \mathbb{R_+}$ converge a
                                $v(x,y,t) \in \Omega \times \mathbb{R_+} $ cuando $ \Delta x, \Delta y, \Delta t \rightarrow 0$
                                         
                        \item \textbf{Consistencia} Decimos que un esquema $ I(\underline{\Omega}, L^n_k u^n_k = G^n_k ) $
                                es puntualmente consistente con la ecuaci\'{o}n diferencial $ I(\underline{\Omega}, \mathcal{L} v=F ) $
                                que aproxima a $ I(\Omega, \mathcal{L} v=F ) $ 
                                en $ ( x, y, t) $ si para cualquier funci\'{o}n $\phi = \phi(x,y,t)$
                                \begin{equation}
                                            (\mathcal{L} \phi - F ) \mid^n_k   - [ L^n_k \phi (j \Delta x, k \Delta y, n \Delta t) - G^n_k] \rightarrow 0
                                \end{equation}
                                cuando $ \Delta x, \Delta y, \Delta t \rightarrow 0$ 
                                y $ (j \Delta x, k \Delta y, (n+1) \Delta t) \in \underline{\Omega} \times \mathbb{R_+} \rightarrow  (x,y,t) \in \Omega \times \mathbb{R_+} $
                                
                        \item \textbf{ Estabilidad } El esquema $ I(\underline{\Omega}, L^n_k u^n_k = G^n_k ) $  
                               se dice que es estable con respecto a una norma $ \left \| \cdot \right \|$ si existen 
                               constantes positivas $\Delta x_0, \Delta y_0, \Delta t_0$ y constantes no negativas $k$
                               y $\beta$ tal que 
                               \begin{equation}
                                         \left \|  u^{n+1} \right \| \leq k e^{\beta t} \left \| u^0 \right \|,
                               \end{equation}
                               para todo $0 \leq t= (n+1) \Delta t$, $0 \leq \Delta x \leq \Delta x_0$, $0 \leq \Delta y \leq \Delta y_0$
                               y $0 \leq \Delta t \leq \Delta t_0$.
                \end{enumerate}
                
                De manera an\'{a}loga se define la \textbf{Convergencia}, \textbf{Consistencia}
                y \textbf{Estabilidad} para el esquema $ I(\overline {\Omega}, L^n_k u^n_k = G^n_k )$.
                
                Ahora bien, si los esquemas $ I(\underline{\Omega}, L^n_k u^n_k = G^n_k ) $ 
                y $ I(\overline {\Omega}, L^n_k u^n_k = G^n_k )$ satisfacen las propiedades 
                antes mencionadas, entonces es posible obtener una buena aproxi-maci\'{o}n a 
                la soluci\'{o}n de $ I(\Omega, \mathcal{L} v=F ) $ mediante el m\'{e}todo de
                diferencias finitas. 
                
                Sea $ \Omega= \{ (x,y) \mid x_1  \leq x \leq x_2  , y_1 \leq y \leq y_2 \}$
                una regi\'{o}n rect\'{a}ngular. Si $R=\Omega$, entonces la discretizaci\'{o}n de $ \Omega $
                coincide con $\overline {\Omega}$ y adem\'{a}s 
                \begin{equation}
                       I(\overline {\Omega}, L^n_k u^n_k = G^n_k ) = I(\Omega, L^n_k u^n_k = G^n_k )=  I(\underline{\Omega}, L^n_k u^n_k = G^n_k ), 
                \end{equation} 
                cuando $ \Delta x, \Delta y \rightarrow 0 $.
                
                De hecho es posible encontar una gran cantidad de ejemplos de este tipo dentro
                de la literatura matem\'{a}tica, ver por ejemplo \cite{Jwth} y \cite{Jwt2}.

\bibliographystyle{apalike}
\bibliography{bibliografia}

\end{document}